\documentclass[12pt,a4paper]{article}
\usepackage[english]{babel}
\usepackage{verbatim}
\usepackage{bm}

\usepackage{amsmath}
\usepackage{amssymb}
\usepackage{amsfonts}
\usepackage{amsbsy}

\usepackage{graphicx}
\usepackage{graphics}
\usepackage{subfigure}
\usepackage{epsfig}
\usepackage{color}
\usepackage{array}
\usepackage{booktabs}
\usepackage{caption}
\usepackage{tabularx}

\begin{document}

\title{The TR-BDF2 method for second order\\
 problems in structural mechanics}

\author{Luca Bonaventura$^{(1)}$,\ \  Macarena G\'omez M\'armol$^{(2)}$}

\maketitle

\begin{center}
{\small
$^{(1)}$ MOX -- Modelling and Scientific Computing, \\
Dipartimento di Matematica ``F. Brioschi'', Politecnico di Milano \\
Via Bonardi 9, 20133 Milano, Italy\\
{\tt luca.bonaventura@polimi.it}
}
\end{center}

\begin{center}
{\small
$^{(2)}$
Departamento
de Ecuaciones Diferenciales y An\'alisis Num\'erico, \\
Universidad de Sevilla\\
Apdo.\,de correos 1160,  41080 Sevilla, Spain\\
{\tt macarena@us.es}
}
\end{center}

\date{}


\vspace*{0.5cm}


\vspace*{0.5cm}

\pagebreak

\abstract{The application of the TR-BDF2 method to second order problems typical of structural mechanics and seismic engineering is discussed. A  reformulation of this method is presented, that only requires the solution of algebraic
systems of size equal to the number of  displacement  degrees of freedom. A linear analysis and
 numerical experiments on relevant benchmarks show that the TR-BDF2  method is superior in terms of accuracy and efficiency
 to the classical Newmark method  and to  its generalizations.}

\pagebreak

 \section{Introduction}
 \label{intro} \indent
 
 The Newmark method \cite{newmark:1959} is an implicit time discretization technique for second order ordinary differential
 equations that is very widely used in structural dynamics applications,  when the stiffness of the system causes standard
 explicit methods such as the St\"ormer-Verlet method \cite{hairer:1993} to be inefficient.
The so called Generalized $\alpha - $ methods \cite{chung:1993} are extensions of
  the Newmark method  using different time averaging parameters for different forcing terms. 
 A full reference list and discussion of the properties
 of these methods are reported in \cite{erlicher:2002}, along with an analysis of their behaviour in the limit of arbitrarily large frequencies.
 Generalized  $\alpha-$methods (shortly denoted as G($\alpha$) methods in the following) introduce
  numerical dissipation in the approximation of the highest frequency modes of the solution, thus allowing to achieve unconditional stability   and, in some regimes, to avoid  overshoots, as discussed in \cite{erlicher:2002}. However,
  in fully nonlinear regimes spurious oscillations may still be present. Furthermore, the stability and accuracy are dependent
  on the values of the numerical parameters that define the method, which have to be tuned for each specific application.
  For these reasons, alternatives to the G($\alpha$) methods have been sought in \cite{owren:1995}, \cite{piche:1995}
  in the class of  Singly Diagonally Implicit Runge Kutta (SDIRK)  and Rosenbrock methods, respectively.
  These investigations, which date back more than 25 years by now, do not seem to have changed
  the attitude of practitioners, who apparently have kept on using G($\alpha$) methods for this kind of problems,
  apart from a limited number of exceptions, see e.g. \cite{bursi:2008}, \cite{bursi:2011}, \cite{hamkar:2012}, \cite{hartmann:2007}, 
  \cite{lamarche:2009}, \cite{meijaard:2003}. One possible reason of this preference is  that, in the straightforward application of Runge-Kutta type methods to structural dynamics, the size of the system to be solved is  twice that required by G($\alpha$)
  methods.
  
  In this work, we extend the analysis and comparison of \cite{owren:1995} to the  TR-BDF2 method
   \cite{bank:1985}, a second order accurate method whose remarkable accuracy, stability and efficiency
   properties  have been fully analyzed in \cite{hosea:1996}. In particular, this L-stable method is endowed with an embedded third order
   method that permits effective  time step adaptation. Unconditionally monotonic extensions of the  TR-BDF2 method
   have been recently derived in \cite{bonaventura:2017} and its multirate version has been shown in \cite{bonaventura:2020a}
   to be quite effective in reducing the computational cost of the time discretization for multirate hyperbolic problems.
   Furthermore, the combination of this method with high order, adaptive discontinuous finite element discretizations
   for model problems represented by first order equations in time has been shown in \cite{tumolo:2015} to be extremely
   effective in reducing the computational cost of these high order methods.
   
   For these reasons, it is of interest to derive a form of this method that allows its application to structural mechanics problems.
   More specifically,
   we present a  reformulation of the TR-BDF2 method for structural mechanics problems that only implies the solution of nonlinear systems
  of the same size as the number of displacement  degrees of freedom.
  We also show how 
  the velocity  degrees of freedom  do not have to be stored explicitly and can be recomputed whenever needed,
   thus avoiding excessive memory requirements for large scale applications. 
   A similar reformulation was presented in \cite{piche:1995} for a two stage Rosenbrock method. 
   An analysis of the accuracy and dissipation properties
   of the TR-BDF2 method   is carried out in the linear regime, extending the classical analyses to consider also the approximation of the damping   terms. Both this analysis and  a number of numerical experiments on relevant benchmarks show that 
    TR-BDF2 is superior  in terms of accuracy and efficiency
 to the classical  G($\alpha$) methods.

 The paper is organized as follows.  In section \ref{trbdf2_struct} we show
 how the TR-BDF2 method can be applied to second order problems and  reformulated so as to avoid solving systems
 of dimension larger than the number of displacement degrees of freedom, while at the same time not requiring
 extra storage of the velocity degrees of freedom.
 In section \ref{analysis}, the dissipation properties of the method are compared with those of the G($\alpha$) methods.
 In section \ref{tests} results of several numerical simulations are presented, which highlight the
 accuracy and efficiency of the method for structural mechanics applications. Some conclusions and perspectives
 for future work are presented in section \ref{conclu}.

 \section{The  TR-BDF2 method for second order problems}
 \label{trbdf2_struct} 
 \indent
 
  We first consider the generic Cauchy problem $ {\bf y}^{\prime}={\bf f}({\bf y},t) $
  on the time interval $t\in [0,T].$ 
  Even though all the results presented in the following also hold for variable time steps, 
  we present for simplicity the TR-BDF2 method as
 employing  a constant time step $h=T/N  $   and we follow in this presentation the notation and conventions
 in \cite{bonaventura:2017}. In its
 original formulation  \cite{bank:1985}, the TR-BDF2 method is defined by the two  stages:
\begin{eqnarray}\label{trbdf2}
  {\bf u}^{n\gamma} - \frac{\gamma h}2  {\bf f}({\bf u}^{n\gamma},t_{n}+\gamma h)
   &=& {\bf u}^n + \frac{\gamma h}2 {\bf f}({\bf u}^{n},t_{n}), \nonumber \\
    && \ \ \nonumber \\
  {\bf u}^{n+1} - \gamma_2  h  {\bf f}({\bf u}^{n+1},t_{n+1}) &=& (1-\gamma_3 ){\bf u}^n +\gamma_3 {\bf u}^{n+\gamma}.
\end{eqnarray}
Here,  ${\bf u}^{n}$ denotes the numerical approximation of the solution at time level $t_n=nh,$
$\gamma \in (0,1) $ is an implicitness parameter and $\gamma_2 =  (1-\gamma)/(2-\gamma), 
 \gamma_3  =1/\gamma(2-\gamma).$
The first stage  of \eqref{trbdf2} is nothing but the application of the trapezoidal rule
 (or Crank-Nicolson method) over the interval 
 $[t_n,t_{n}+\gamma h].$
 The outcome of this stage is then used as input for the BDF2 implicit method.
 The resulting combination yields a method with several interesting accuracy and stability
 properties.  A detailed analysis of these properties is given in \cite{hosea:1996}, where it
 is also shown that TR-BDF2  is equivalent to a three stage  Diagonally Implicit Runge Kutta (DIRK) method
defined by the stages
\begin{eqnarray}
{\bf k}_1& =&  {\bf f}\left ({\bf u}^n, t_{n}\right) \nonumber \\
 {\bf k}_2 &=& {\bf f}\left ({\bf u}^n +\frac{\gamma h}2  {\bf k}_1  +\frac{\gamma h}2  {\bf k}_2,   t_{n}+\gamma  h\right)\nonumber \\
 {\bf k}_3 &=& {\bf f}\left ({\bf u}^n +\frac{\gamma_3h}2 {\bf k}_1+ \frac{\gamma_3h}2{\bf k}_2 
 +\gamma_2 h{\bf k}_3,t_{n+1} \right )
 \nonumber \\ 
 {\bf u}^{n+1} &=& {\bf u}^n +h \left ( \frac{\gamma_3}2{\bf k}_1+ \frac{\gamma_3}2{\bf k}_2 +\gamma_2{\bf k}_3 \right ).
\label{tr_dirk}
\end{eqnarray}
Notice that this method is not  Singly Diagonally Implicit Runge Kutta (SDIRK), due to the fully explicit first stage.
As shown in  \cite{hosea:1996},
the TR-BDF2 method is second order accurate and A-stable  for any value of $\gamma.$ 
Furthermore, for $\gamma=2-\sqrt{2} $ it is also  L-stable. Therefore, with this coefficient
value it can be safely applied to 
 problems with eigenvalues whose imaginary part is large, such as typically arise from the discretization
 of highly oscillatory second order  problems.  This is not the case for the standard 
trapezoidal rule (or Crank-Nicolson) implicit method, whose linear stability region is exactly bounded by
the imaginary axis.   Recently, the method   has been applied
   in \cite{tumolo:2015}, \cite{tumolo:2016} to high order accuracy spatial discretization of wave propagation phenomena written as  systems
  of first order in time. Its monotonicity properties have also been studied in \cite{bonaventura:2017}, where it was shown
  that an unconditionally monotonic extension of the method can be derived and that the method does not suffer from order reduction when applied to stiff systems. Several specific variants of the method
  for nonlinear problems have been proposed and analyzed in \cite{edwards:2011}, while an analysis of a multirate implicit method
  based on TR-BDF2 has been presented in \cite{bonaventura:2020a}.

We   now focus on systems of
second order  ordinary differential equations. In particular, as common in structural mechanics applications,
we will consider systems  
 of the form
  \begin{equation}
 \label{2ndcauchypb}
 {\bf M}  {\bf y}^{\prime\prime}=    -{\bf C}  {\bf y}^{\prime}-{\bf K}{\bf y} +{\bf g}( {\bf y}) +{\bf z}(t)
\end{equation}
with the initial conditions ${\bf  y}(0) = {\bf y}_0 \ \ \ \ {\bf  y}^{\prime}(0) = {\bf v}_0.$
Here, $ {\bf M} $ denotes a symmetric and positive definite mass matrix,
 $ {\bf K} $ denotes the stiffness matrix, which is also assumed to be symmetric and positive definite, while
the matrix  ${\bf C}  $ represents friction terms. 
  
 In order to introduce the application of the TR-BDF2 method
to this kind of problems, system \eqref{2ndcauchypb} is rewritten as a first order system
letting $ {\bf v}= {\bf y}^{\prime} $ 
and setting
${\boldsymbol \Gamma}= {\bf M} ^{-1}{\bf C},$  ${\bf B}={\bf M}^{-1} {\bf K} , $ 
and
${\bf f}( {\bf y}) ={\bf M}^{-1}{\bf g}( {\bf y}), $ ${\bf s}(t) ={\bf M}^{-1}{\bf z}(t),$   so as to obtain

\begin{eqnarray}
 \label{structpb}
   {\bf y}^{\prime}&=&   {\bf v}     \nonumber \\
    {\bf v}^{\prime}&=& -{\boldsymbol \Gamma}  {\bf v}-{\bf B}{\bf y }+{\bf f} ({\bf y}) + {\bf s}(t).  
\end{eqnarray}
%
%
%
We now apply method  \eqref{trbdf2} to problem \eqref{structpb}. Denoting by
${\bf u}^n ,{\bf w}^n,  $ respectively, the numerical approximations of ${\bf y}, {\bf v} $ at time level $n, $
 the  trapezoidal rule stage can be written as

\begin{eqnarray}\label{trbdf2_1st}
    {\bf u}^{n+\gamma} &=&  {\bf u}^n + \frac{\gamma  h}2 {\bf w}^{n}+\frac{\gamma  h}2   {\bf w}^{n+\gamma}    \\
 &&   \nonumber \\
  {\bf w}^{n+\gamma} &+&\frac{\gamma  h}2  {\boldsymbol \Gamma} {\bf w}^{n+\gamma} 
+\frac{\gamma  h}2  {\bf B} {\bf u}^{n+\gamma}  - \frac{\gamma  h}2  {\bf f}({\bf u}^{n+\gamma})   \nonumber \\
  &=& {\bf w}^n  -\frac{\gamma  h}2  {\boldsymbol \Gamma} {\bf w}^{n} 
-\frac{\gamma  h}2 {\bf B} {\bf u}^{n} + \frac{\gamma  h}2  {\bf f}({\bf u}^{n})+ \frac{\gamma  h}2{ \bar {\bf s}},  
\end{eqnarray}
where we have set ${ \bar {\bf s}}= {\bf s}(t^{n+\gamma}) + {\bf s}(t^n). $ The BDF2 stage  yields instead

\begin{eqnarray}\label{trbdf2_2nd}
  {\bf u}^{n+1}  &=&  \gamma_2  h  {\bf w}^{n+1} + (1-\gamma_3 ){\bf u}^n +\gamma_3 {\bf u}^{n+\gamma} \\
 && \ \ \nonumber \\
    {\bf w}^{n+1} &+&\gamma_2  h {\boldsymbol \Gamma} {\bf w}^{n+1} 
+\gamma_2 h {\bf B} {\bf u}^{n+1}  - \gamma_2  h  {\bf f}({\bf u}^{n+1})  \nonumber \\
  &=& (1-\gamma_3 ){\bf w}^n +\gamma_3 {\bf w}^{n+\gamma} + \gamma_2 h {\bf s}^{n+1}. \
\end{eqnarray}
Along the lines of what is done for the Newmark and G($\alpha$) methods, each of the two stages is  now
rewritten in terms of a single implicit step for ${\bf u}^{n+\gamma}, {\bf u}^{n+1}, $ respectively.
This amounts to 

\begin{eqnarray}\label{trbdf2_1st_un}
   {\bf w}^{n+\gamma} &=& 2\left (  {\bf u}^{n+\gamma} -{\bf u}^n - \gamma  h {\bf w}^{n}/2  \right)/{\gamma  h} \\
  && \ \ \nonumber \\
  {\bf A}_1{\bf u}^{n+\gamma} &-& \frac{\gamma^2  h^2}4  {\bf f}({\bf u}^{n+\gamma}) ={\bf b}_1 
\end{eqnarray}

\begin{eqnarray}\label{trbdf2_2nd_un}
  {\bf w}^{n+1} &=&   \left [ {\bf u}^{n+1} - (1-\gamma_3 ){\bf u}^n -\gamma_3 {\bf u}^{n+\gamma}  \right]/\gamma_2  h \\
 && \ \ \nonumber \\
   {\bf A}_2{\bf u}^{n+1}&-& \gamma_2^2  h^2  {\bf f}({\bf u}^{n+1}) ={\bf b}_2 
\end{eqnarray}
where we have now set

\begin{eqnarray}
  {\bf A}_1&=&  {\bf I}  + \frac{\gamma  h}2  {\boldsymbol \Gamma} +\frac{\gamma^2  h^2}4 {\bf B}    \\
    {\bf A}_2&=&  {\bf I}  + \gamma_2  h{\boldsymbol \Gamma} +\gamma_2^2  h^2 {\bf B}    \\
   && \ \ \nonumber \\
  {\bf b}_1&=&  \left({\bf I}  + \frac{\gamma  h}2  {\boldsymbol \Gamma} -\frac{\gamma^2  h^2}4 {\bf B}\right ) {\bf u}^{n}  \nonumber \\
                 &+&\frac{\gamma^2  h^2}4 {\bf f}({\bf u}^{n}) +\frac{\gamma^2  h^2}4{ \bar {\bf s}}  +\gamma h {\bf w}^n\\
                  && \ \ \nonumber \\
                 {\bf b}_2&=& \gamma_2^2 h^2 {\bf s}^{n+1} +\gamma_2(1-\gamma_3) h{\bf w}^n   +\gamma_2\gamma_3 h{\bf w}^{n+\gamma} 
                 \nonumber \\
                 &+&(1-\gamma_3 )({\bf I}  + \gamma_2  h{\boldsymbol \Gamma} ) {\bf u}^{n} 
                  +\gamma_3({\bf I}  + \gamma_2  h{\boldsymbol \Gamma}) {\bf u}^{n+\gamma}. 
                  \end{eqnarray}
Fixing from now on the value  $\gamma=2-\sqrt{2},$ so as to achieve L-stability, one obtains
 $\gamma_2=\gamma/2, $ so that ${\bf A}_1={\bf A}_2={\bf A}.$
 Notice that this formulation has a number of advantages. 
 First, it allows in practice to avoid doubling the degrees of freedom of the discrete
 problem, even though a method for first order problems is employed. For the first time step, the value
 of ${\bf w}^n$ would be recovered from the initial datum, while formula \eqref{trbdf2_2nd_un}
 allows the reconstruction of the same term at subsequent time levels. Furthermore,  
both stages of the method are defined in terms of the same matrix. 
 In order to improve the  efficiency of the algebraic solvers, 
  both equations can also be rewritten in terms of the
 increments
 $$ \delta {\bf u}^{n+\gamma}=  {\bf u}^{n+\gamma} -{\bf u}^{n} 
 \ \ \ \ \delta {\bf u}^{n+1}=  {\bf u}^{n+1} -{\bf u}^{n+\gamma}, $$
 so that the scheme can be rewritten as
 \begin{eqnarray}\label{trbdf2_delta}
&&  {\bf A}\delta {\bf u}^{n+\gamma} - \frac{\gamma^2  h^2}4  {\bf f}({\bf u}^{n}+\delta {\bf u}^{n+\gamma}) = {\bf b}_1 -{\bf A}{\bf u}^{n}  \\
&& {\bf A}\delta{\bf u}^{n+1} - \frac{\gamma^2  h^2}4  {\bf f}({\bf u}^{n+\gamma}+\delta{\bf u}^{n+1}) ={\bf b}_2 -{\bf A}{\bf u}^{n+\gamma}. 
\end{eqnarray}
Finally, it is easy to notice that the mass matrix inversion is not necessary in practice. Indeed, each of the previous 
equations can be multiplied by $ {\bf M}, $ thus yielding the two nonlinear systems

\begin{eqnarray}\label{trbdf2_delta_nom}
&&  {\bf A}\delta {\bf u}^{n+\gamma} - \frac{\gamma^2  h^2}4  {\bf g}({\bf u}^{n}+\delta {\bf u}^{n+\gamma}) ={\bf \tilde b}_1 \\
&& {\bf A}\delta{\bf u}^{n+1} - \frac{\gamma^2  h^2}4  {\bf g}({\bf u}^{n+\gamma}+\delta{\bf u}^{n+1}) =\tilde{\bf b}_2 
\end{eqnarray}
where now $  {\bf A} $ has been redefined as
 $ {\bf A}=  {\bf M}  + \gamma  h{\bf C}/2 +\gamma^2  h^2 {\bf K}/4 $
and the corresponding right hand sides as

\begin{eqnarray}
\label{trbdf2_delta_rhs}
\tilde{\bf b}_1&=&  \left({\bf M}  + \frac{\gamma  h}2  {\bf C} -\frac{\gamma^2  h^2}4 {\bf K}\right ) {\bf u}^{n}  \nonumber \\
                 &+&\frac{\gamma^2  h^2}4 {\bf g}({\bf u}^{n}) +\frac{\gamma^2  h^2}4{ \bar {\bf z}} +\gamma h {\bf w}^n -{\bf A}{\bf u}^{n}  \nonumber\\
                  && \ \ \nonumber \\
                 \tilde  {\bf b}_2&=& \gamma_2^2 h^2 {\bf z}^{n+1} +\gamma_2(1-\gamma_3) h{\bf M}{\bf w}^n   +\gamma_2\gamma_3 h{\bf M}{\bf w}^{n+\gamma}  \nonumber \\
                 &+&(1-\gamma_3 )({{\bf M}}  + \gamma_2  h{{\bf C}} ) {\bf u}^{n} 
                  +\gamma_3({\bf M}  + \gamma_2  h{\bf C}) {\bf u}^{n+\gamma} -{\bf A}{\bf u}^{n+\gamma}.  \nonumber
\end{eqnarray}

%

 For completeness, we also present the classical 
 Newmark method and its G($\alpha$) generalizations. Following \cite{erlicher:2002}, 
 we introduce  discrete approximations  $ {\bf u}^n, {\bf w}^n,{\bf a}^n $ of the continuous displacement,
 velocity and acceleration values, respectively:
 
\begin{eqnarray}\label{eq:galfa}
&&{\bf u}^{n+1}={\bf u}^{n}+h{\bf w}^n+ h^2\left [\left (\beta -\frac 12\right ){\bf a}^n+\beta {\bf a}^{n+1}\right] \\
  &&{\bf w}^{n+1}={\bf w}^n+h [  (1-\gamma_N ){\bf a}^n+\gamma_N {\bf a}^{n+1} ]  \\
   &&  {\bf M}   {\bf a}^{n+1-\alpha_m}+{\bf C} {\bf w}^{n+1-\alpha_f}+{\bf K}{\bf u}^{n+1-\alpha_f} 
   \nonumber \\
   && \ \ \ \ \ \ \ \ \ \ \ \ \ \ \ ={\bf g}({\bf u}^{n+1-\alpha_f}) + {\bf z}^{n+1-\alpha_f},
  \end{eqnarray}
where for a generic variable $\phi $ one has $\phi^{n+1-\alpha} = ( 1- \alpha) \phi^{n+1} + \alpha\phi^{n}. $
where  $\gamma_N,  \beta, \alpha_m,\alpha_f $ are  specific averaging parameters.
Notice that different definitions of the $\alpha $ averaging  can also be employed and that
several relationships between  the values of the parameters $\gamma_N, \beta, \alpha_m,\alpha_f $ must be respected
for the method to be second order consistent  in the displacement and unconditionally stable.
 We refer again to \cite{erlicher:2002} for a complete discussion of these issues and of a simple strategy to link 
the values of these parameters to that of the effective dissipation rate of the highest frequencies, denoted
by $\rho_{\infty}.$  Here, we just summarize the key results, which are reported in the following expression

\begin{equation}
\beta= \frac{1}{(1+\rho_{\infty})^2}\ \ \ \ \gamma_N=\frac 12\frac{3-\rho_{\infty}}{1+\rho_{\infty}}
\end{equation}
for the parameters $\beta, \gamma_N $ and
\begin{equation}
\alpha_f=\frac{2\rho_{\infty}-1}{1+\rho_{\infty}} \ \ \ \ \alpha_m=\frac{ \rho_{\infty}}{1+\rho_{\infty}}
\end{equation}
for the parameters that define the G($\alpha$)  method introduced in \cite{chung:1993}, which we consider
in the following  as example of the
broader G($\alpha$)  class.

 \section{Analysis of the TR-BDF2 method}
\label{analysis} \indent

In  \cite{owren:1995}, several SDIRK methods were compared to some
G($\alpha$) methods with respect to their dissipation properties. The
behaviour of the TR-BDF2 method is very similar to the two stage SDIRK method
analyzed in Section 2.2 of  \cite{owren:1995}. However, the analysis will be
repeated here for convenience and extended to the case in which also dissipative terms are taken into account.
More specifically, we will consider the methods introduced in section \ref{trbdf2_struct} as applied to the linear scalar problem

\begin{equation}
\label{eq:osci}
y^{\prime \prime}= -\frac cm  {y}^{\prime}-{\frac km}{y},
\end{equation}
where $k,m>0, c\geq 0.$
We will only consider the   case in which $ c^2< 4km,$ for which it is well known that
the analytic solution has the form $$y(t)=A_+\exp{(\alpha_+ t)}+A_-\exp{(\alpha_- t)}, $$ where one has
$$\alpha_{\pm}=-\tilde \zeta \pm i\tilde \omega \ \ \
{\rm with} \ \ \  \tilde \zeta=c/2m \geq 0,  \ \ \ \ \tilde \omega= \sqrt{k/m - (c/2m)^2}>0. $$
As a consequence, one can write for this case the exact evolution operator

\begin{equation}\label{eq:ex_evol}
\mathbf{y}^{n+1}
=\mathbf{E}\mathbf{y}^{n},
\end{equation}
where we have set $\mathbf{y}^{n}=[y(t^{n}), hy^{\prime}(t^{n})]^T, $ $\zeta= h \tilde \zeta, $ $\omega= h \tilde \omega$  and
  \begin{equation}
 \label{eq:exmat1}
\mathbf{E}
=     \left[
\begin{array}{cc}
\left(\cos{\omega} +\frac{\zeta}{\omega}\sin{\omega}\right){\rm e}^{-\zeta}&  \frac{\sin{\omega}}{\omega}{\rm e}^{-\zeta}\\
-(\zeta^2+\omega^2)\frac{\sin{\omega}}{\omega} {\rm e}^{-\zeta}&  \left(\cos{\omega} -\frac{\zeta}{\omega}\sin{\omega}\right) {\rm e}^{-\zeta}\\
\end{array}
  \right]. 
\end{equation}
Notice that the velocity variable has been rescaled (see e.g.  \cite{bathe:1973}), so as to obtain
an evolution operator which only includes dimensionally homogeneous quantities. Defining then
${\bf x}^{n}=[u^n, z^n]^T, $ where we set $z^n=hw^n, $ the  discrete evolution operator corresponding to the TR-BDF2 method
can be reconstructed from the specific forms of the TR and BDF2 stages, which can be written as
\begin{equation}\label{eq:trbdf2_lin}
\mathbf{S}_1\mathbf{x}^{n+\gamma}=\mathbf{T}_1\mathbf{x}^{n} \ \ \ \ \mathbf{S}_2\mathbf{x}^{n+1}
=\mathbf{T}_2\mathbf{x}^{n+\gamma} +\mathbf{U}_2\mathbf{x}^{n},
\end{equation}
where one defines
\begin{eqnarray}
\label{eq:trbdf2_lin_mat}
&&\mathbf{S}_1=  \left[
\begin{array}{cc}
 \left(1+\gamma\zeta +(\omega^2+\zeta^2)\frac{\gamma^2}4 \right)&  0 \\
1 &  -\frac{\gamma}2\\
\end{array}
  \right]   
  \nonumber \\
  && 
   \mathbf{T}_1=  \left[
\begin{array}{cc}
 \left(1+\gamma\zeta -(\omega^2+\zeta^2)\frac{\gamma^2}4 \right) & \gamma \\
 1& \frac{\gamma}2 \\
\end{array}
  \right]  \nonumber \\
&& \mathbf{S}_2=  \left[
\begin{array}{cc}
 \left(1+\gamma\zeta +(\omega^2+\zeta^2)\frac{\gamma^2}4 \right) &  0 \\
 1 &  -\frac{\gamma}2 \\
\end{array}
  \right] 
   \nonumber \\
  && 
   \mathbf{T}_2=  \left[
\begin{array}{cc}
 \frac{1}{\gamma(2-\gamma)} \left(1+\frac{\gamma}{2}\zeta\right)& \frac{1}{2(2-\gamma)}  \\
\frac{1}{\gamma(2-\gamma)} &  0\\
\end{array}
  \right]    \nonumber \\
  &&\mathbf{U}_2=  \left[
\begin{array}{cc}
 -\frac{(1-\gamma)^2}{\gamma(2-\gamma)}\left(1+\frac{\gamma}{2}\zeta\right)&  -\frac{(1-\gamma)^2}{2(2-\gamma)} \\
-\frac{(1-\gamma)^2}{\gamma(2-\gamma)} &  0\\
\end{array}
  \right].  
\end{eqnarray}
By straightforward algebraic manipulations, one obtains 

\begin{equation}\label{eq:trbdf2_lin2}
\mathbf{x}^{n+1}
=\mathbf{S}_2^{-1}\left [ \mathbf{T}_2 \left( \mathbf{S}_2^{-1} \mathbf{T}_2\right )     \mathbf{U}_2\right]\mathbf{x}^{n},
\end{equation}
which can then be compared to \eqref{eq:ex_evol} in order to assess the properties of the TR-BDF2 method.
A similar  comparison can be carried out for the G($\alpha $)  methods,
taking into account that these methods also employ an approximation of the second derivative.
Therefore, one has to consider the
exact evolution operator derived from the relationship

\begin{equation}\label{eq:ex_evol_acc}
\left [\begin{array}{cc}
\mathbf{I} & 0 \\
\boldsymbol{\xi} & 1 \\
\end{array} \right ]\left [\begin{array}{c}
\mathbf{y}^{n+1}\\
h^2y^{\prime\prime}(t^{n+1})
\end{array} \right ]
=\left [\begin{array}{cc}
\mathbf{E} & 0 \\
\mathbf{0}  & 0 \\
\end{array} \right ]
 \left [\begin{array}{c}
\mathbf{y}^{n}\\
h^2y^{\prime\prime}(t^{n})
\end{array} \right ],
\end{equation}
where $\boldsymbol{\xi}=[\omega^2+\zeta^2,\zeta]^T  $ and equation \eqref{eq:osci} has been used to introduce
the relationship between acceleration, velocity and displacement values.
Analogously, for the generic G($\alpha $)  method one has

\begin{eqnarray}\label{eq:galfa_evol}
&&\left [\begin{array}{ccc}
1 & 0 & -\beta h^2\\
0 &1 & -\gamma_N \\
(1-\alpha_f) \omega^2& (1-\alpha_f) \zeta &  (1-\alpha_m) \
\end{array} 
\right ]
\left [\begin{array}{c}
u^{n+1}\\
hw^{n+1}\\
h^2a^{n+1}\\
\end{array} \right ] \nonumber \\
&&=\left [\begin{array}{ccc}
1 & 1 & (\beta -1/2) \\
0 & 1 & (1-\gamma_N) \\
\alpha_f\omega^2   & \alpha_f\zeta & \alpha_m \\
\end{array} \right ]
\left [\begin{array}{c}
u^{n}\\
hw^{n}\\
h^2a^{n}\\
\end{array} \right ],
\end{eqnarray}
so that the comparison is between the exact operator

$$  \left [\begin{array}{cc}
\mathbf{I} & 0 \\
\boldsymbol{\xi} & 1 \\
\end{array} \right ]^{-1}\left [\begin{array}{cc}
\mathbf{E} & 0 \\
\mathbf{0}  & 0 \\
\end{array} \right ]$$
and the approximate operator
 
 $$
 \left [\begin{array}{ccc}
1 & 0 & -\beta h^2\\
0 &1 & -\gamma_N \\
(1-\alpha_f) \omega^2& (1-\alpha_f) \zeta &  (1-\alpha_m) \
\end{array} 
\right ]^{-1}\left [\begin{array}{ccc}
1 & 1 & (\beta -1/2) \\
0 & 1 & (1-\gamma_N) \\
\alpha_f\omega^2   & \alpha_f\zeta & \alpha_m \\
\end{array} \right ].
 $$
 
\begin{figure}
 \begin{center}
	\includegraphics[width=0.45\textwidth]{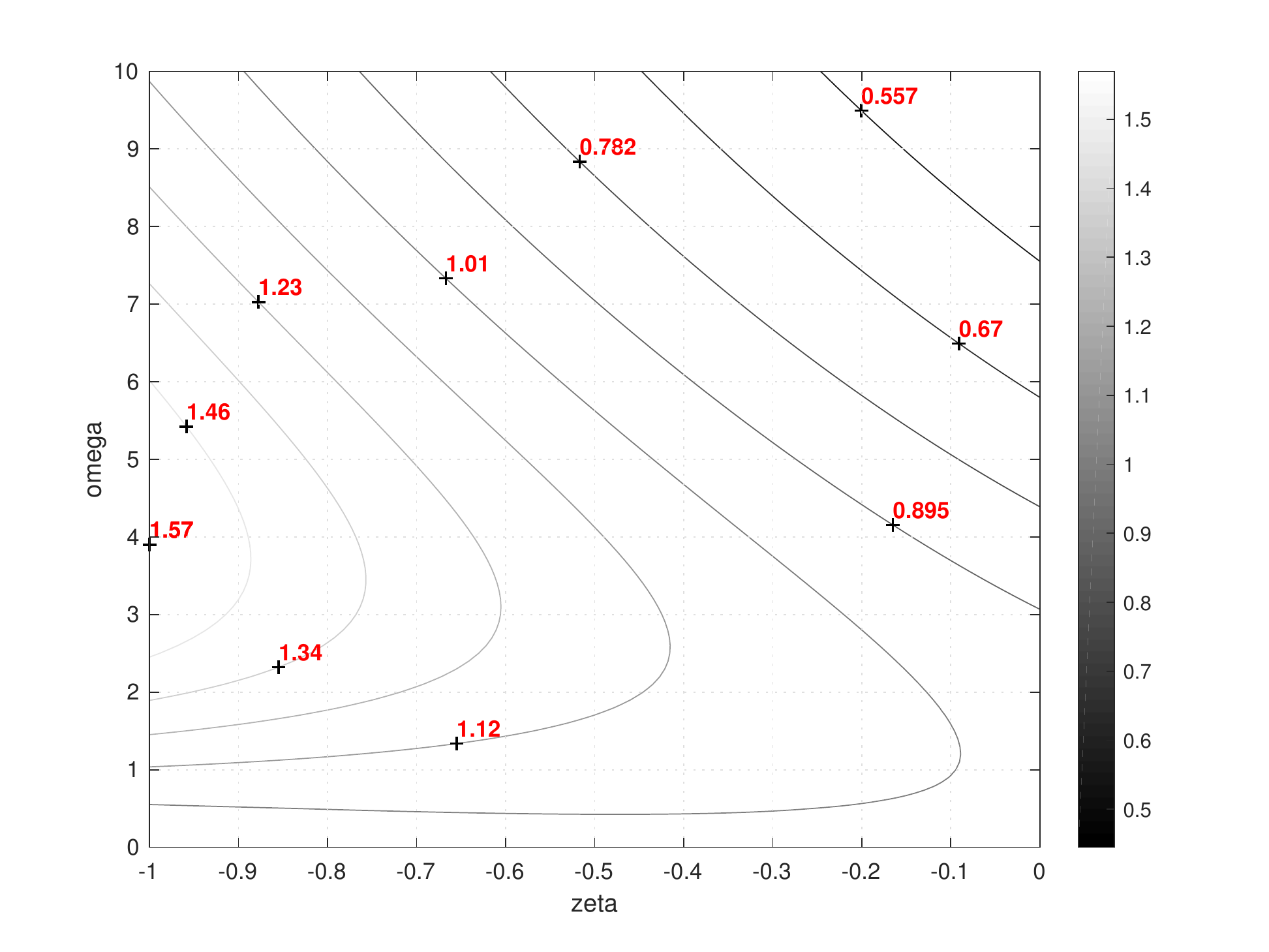}a)
	\includegraphics[width=0.45\textwidth]{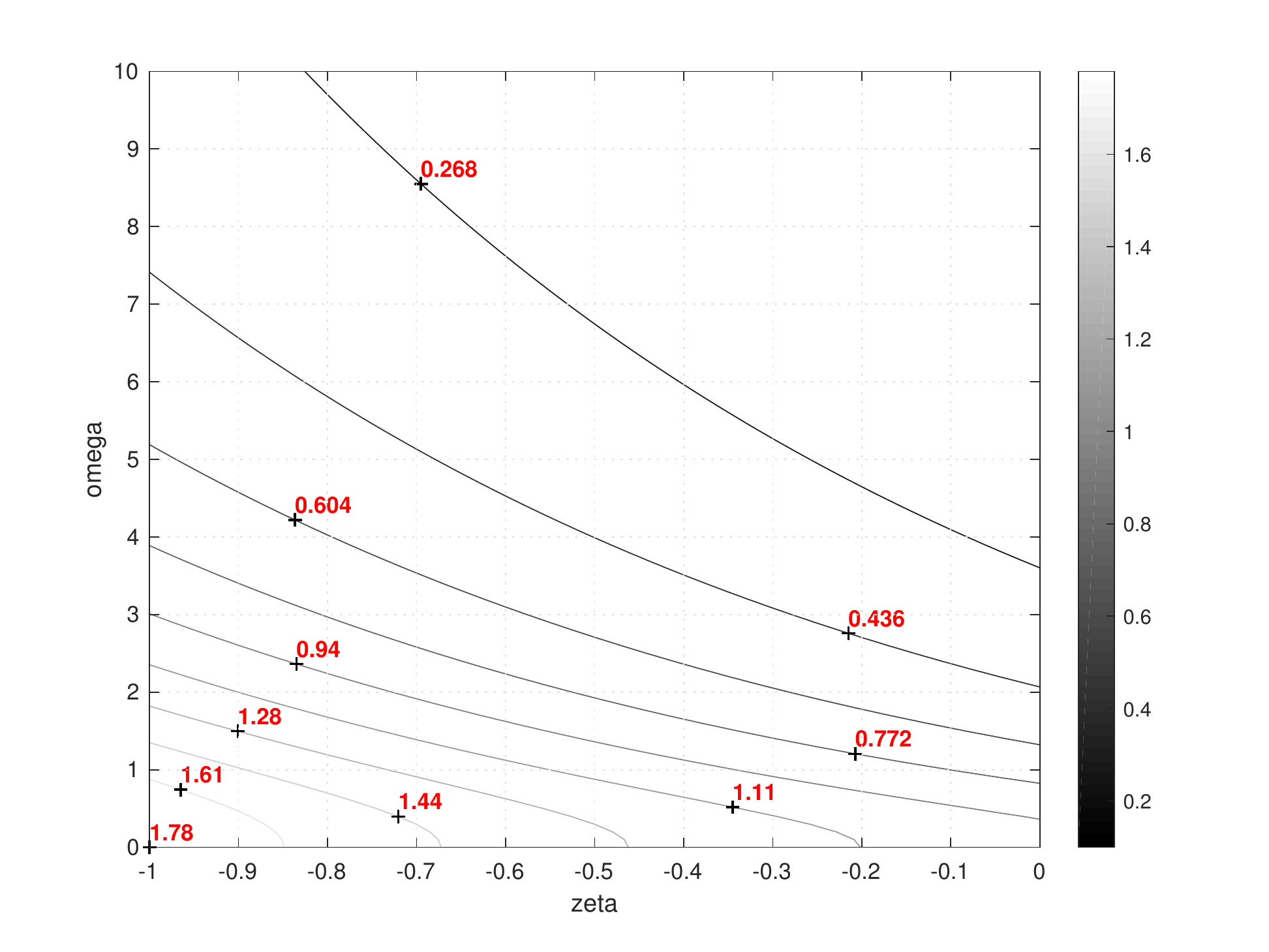}b)
	\includegraphics[width=0.45\textwidth]{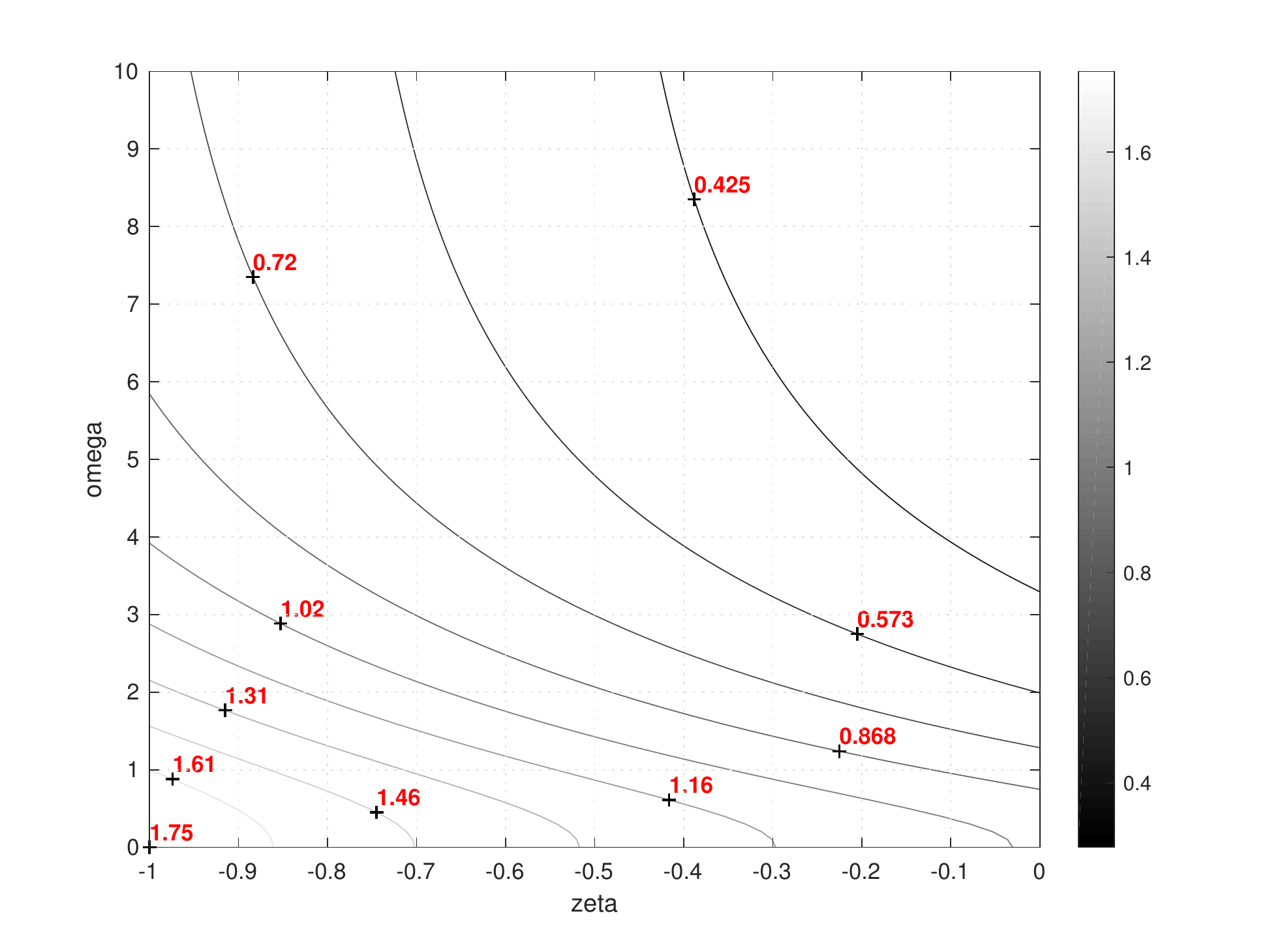}c)
	\includegraphics[width=0.45\textwidth]{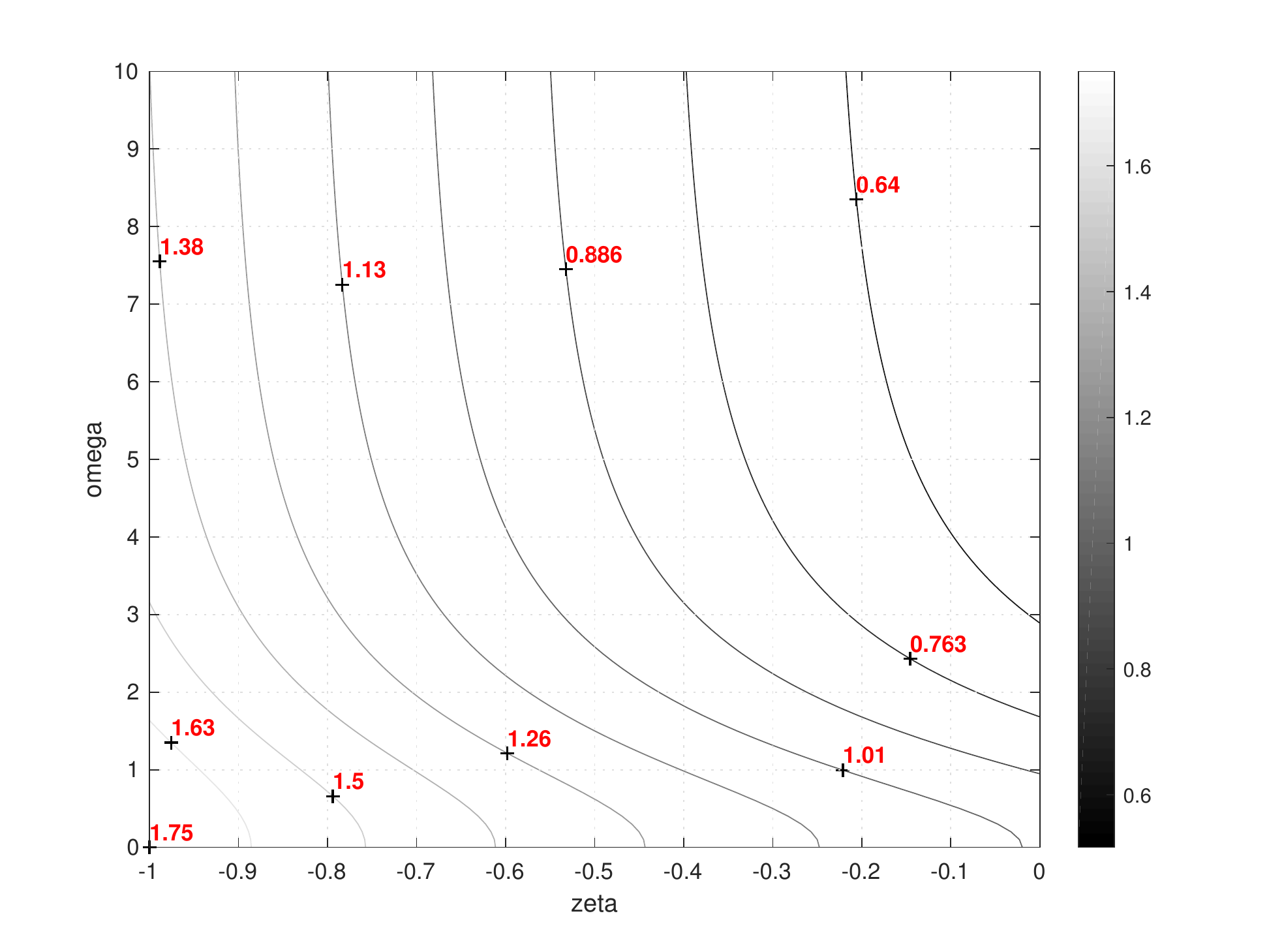}d)
	\includegraphics[width=0.45\textwidth]{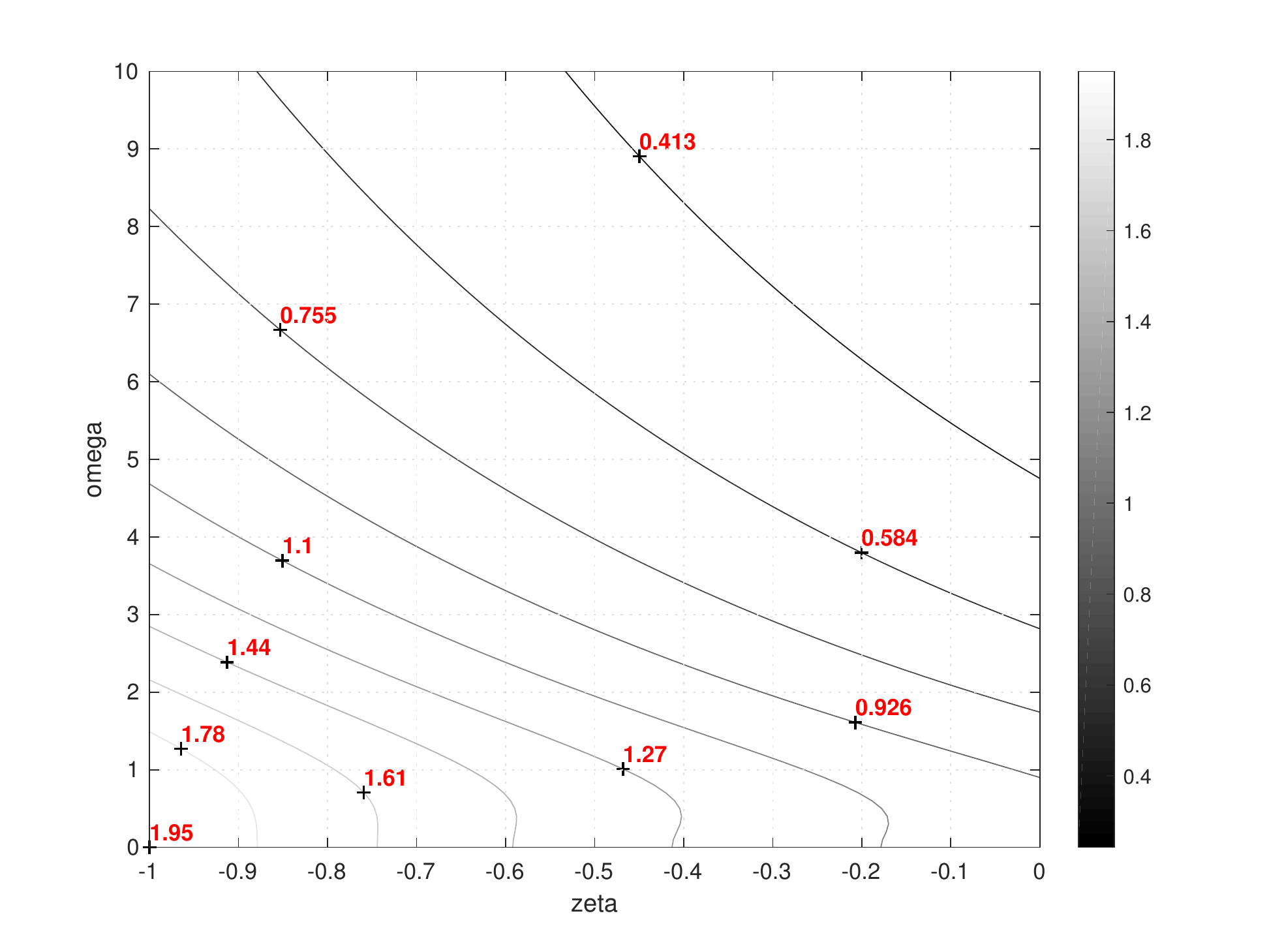}e)
	\includegraphics[width=0.45\textwidth]{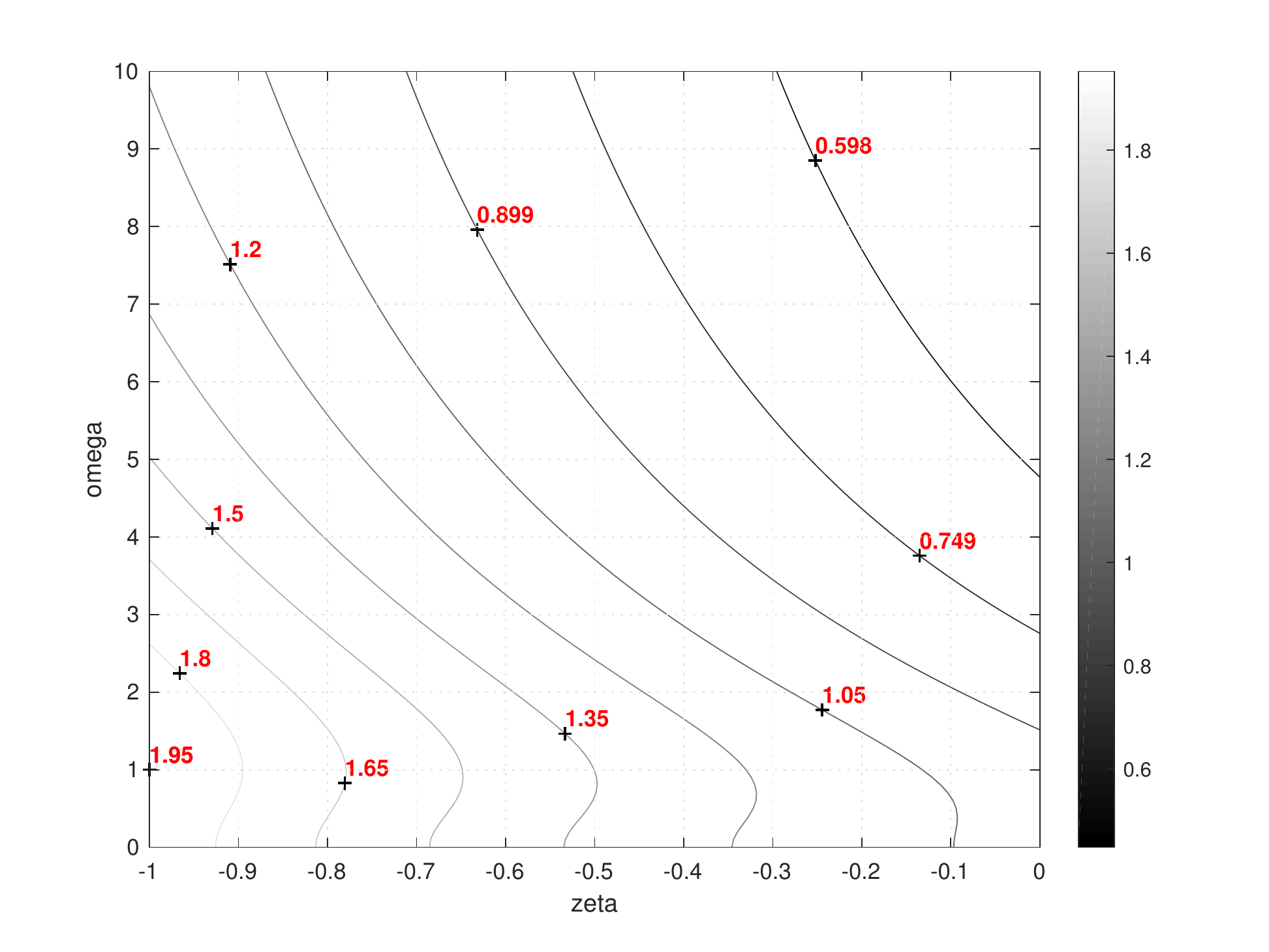}f)
	\end{center}
	\caption{Ratio of the spectral norm of the discrete evolution operators to that of the exact evolution operator
	for a) TR-BDF2 method, b) the Newmark method with $\rho_{\infty}=0,  $  c) the Newmark method with $\rho_{\infty}=0.25, $ 
	 d) the Newmark method with $\rho_{\infty}=0.5 $  e) the CH-$\alpha$ method with $\rho_{\infty}=0 $ 
	  e) the CH-$\alpha$ method with $\rho_{\infty}=0.25.$ }
	\label{fig:ratio_methods}
\end{figure}

\begin{figure}
  \begin{center}
	\includegraphics[width=0.45\textwidth]{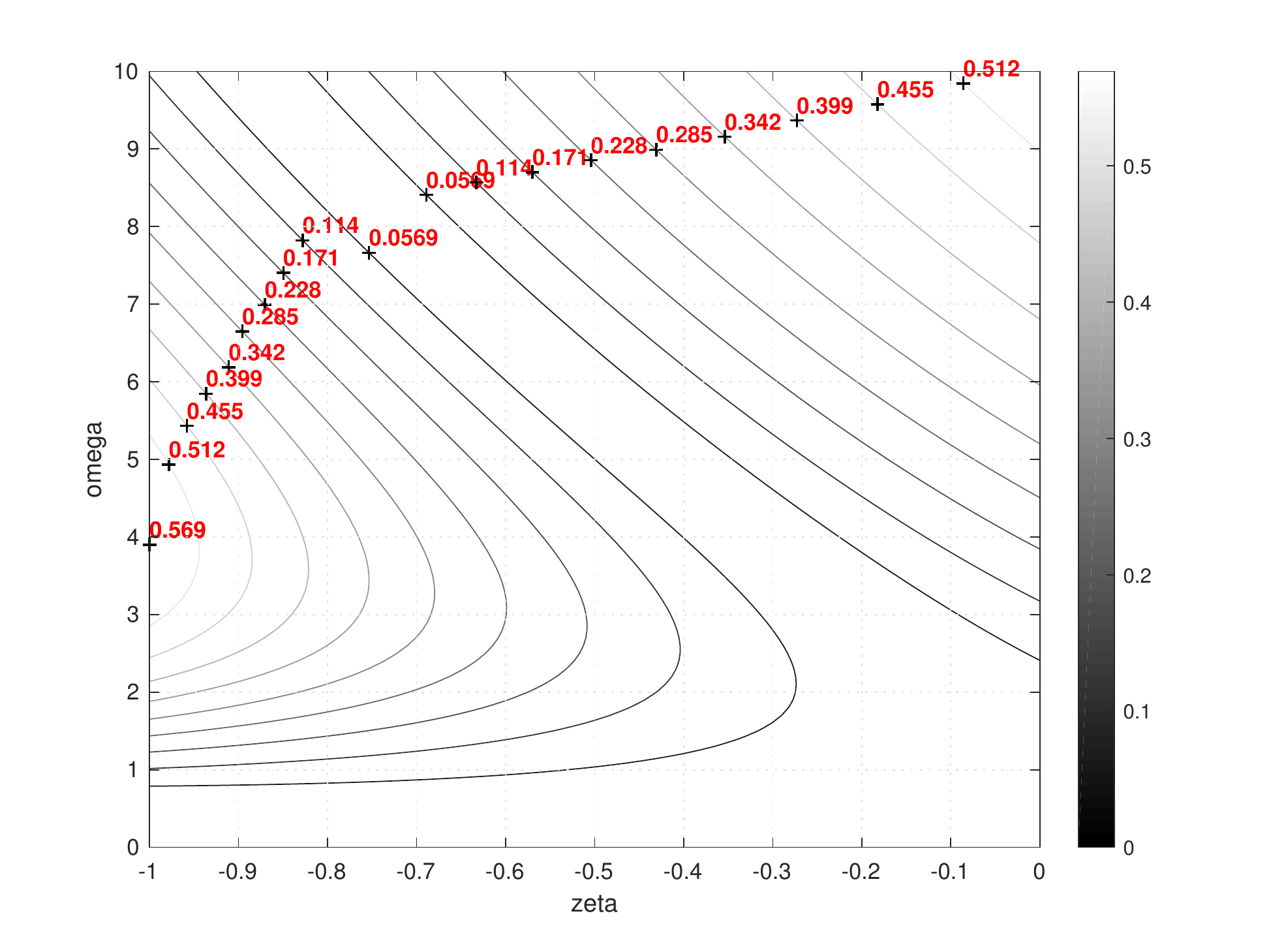}a)
	\includegraphics[width=0.45\textwidth]{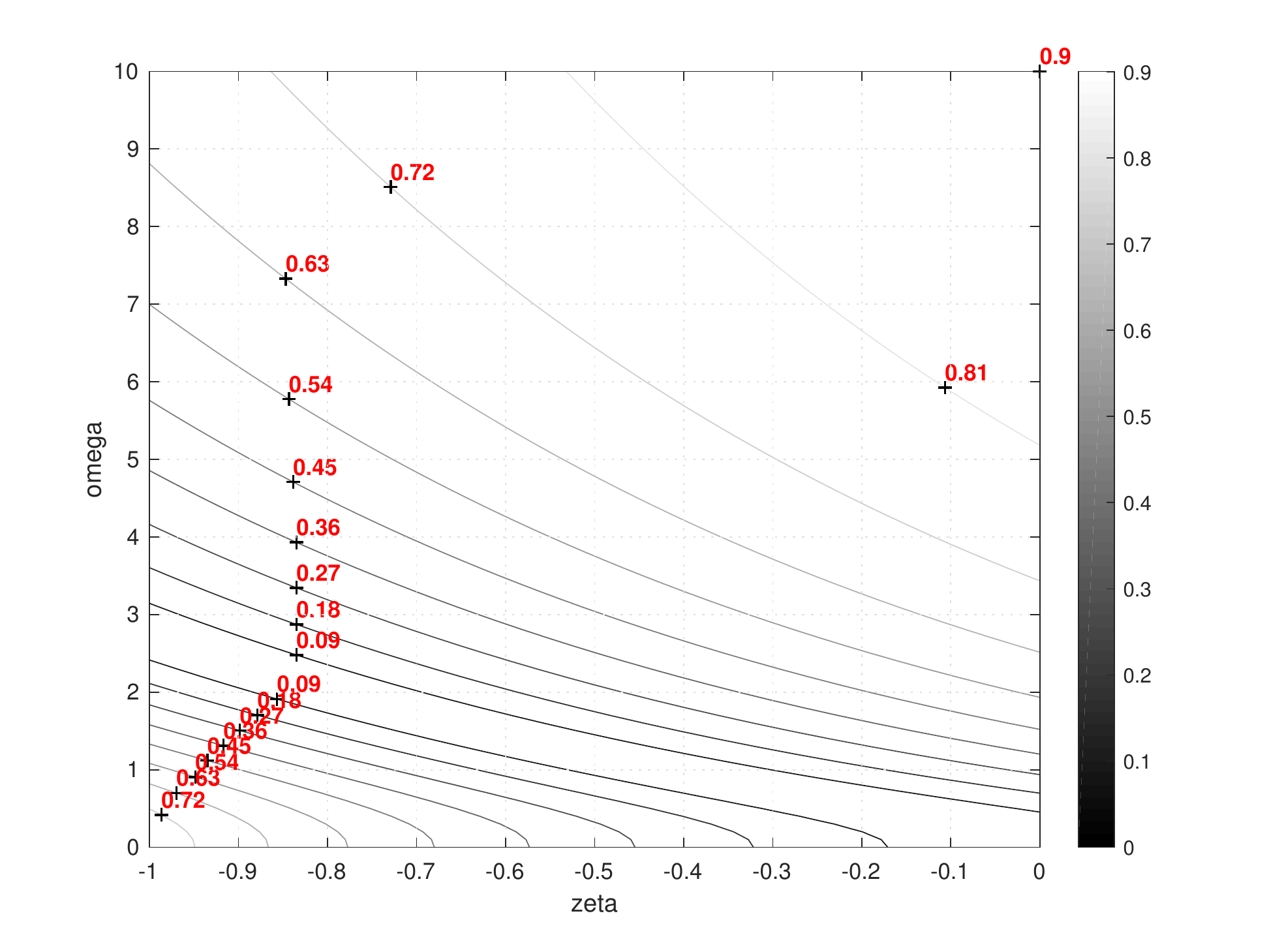}b)
	\includegraphics[width=0.45\textwidth]{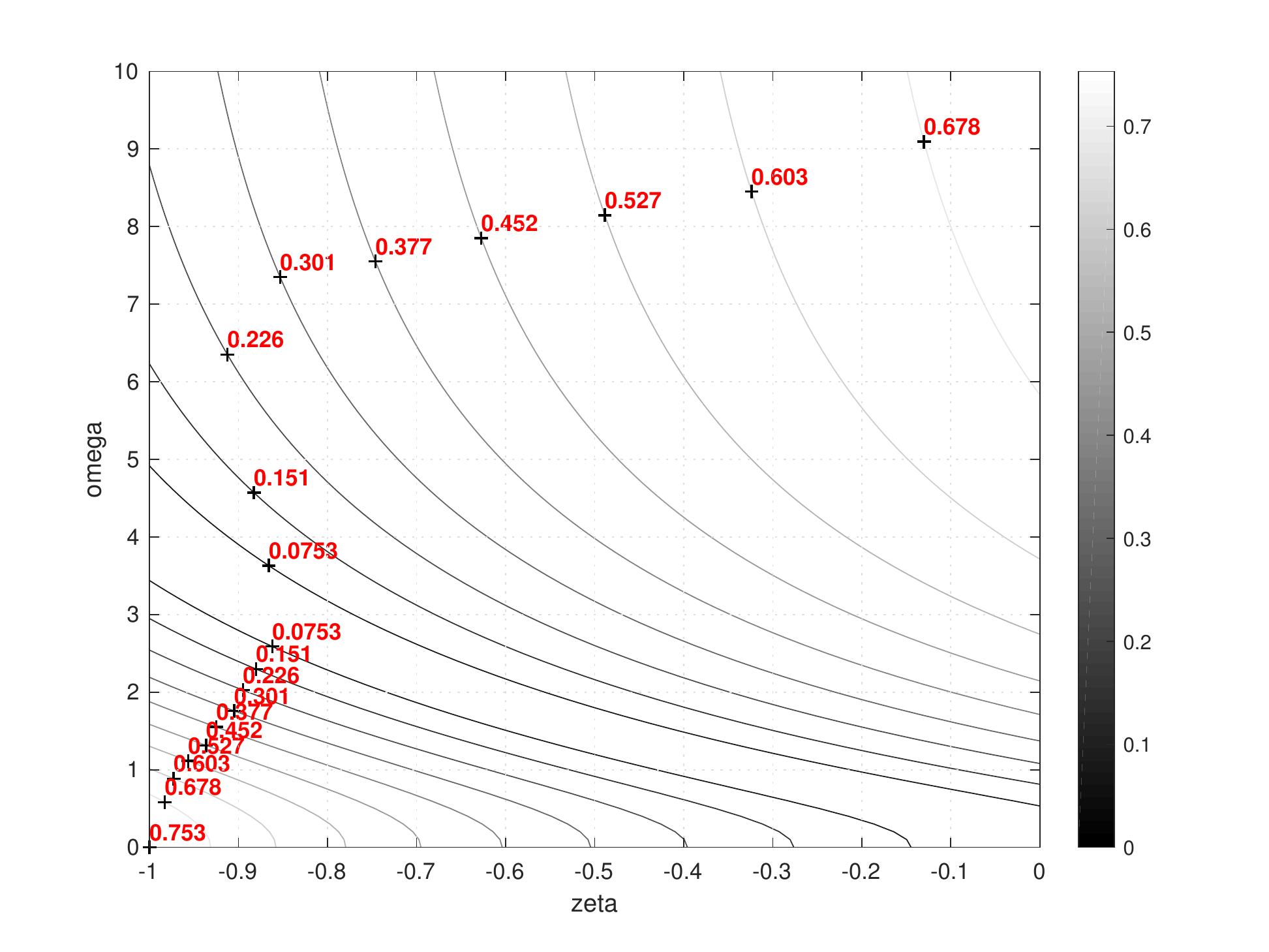}c)
	\includegraphics[width=0.45\textwidth]{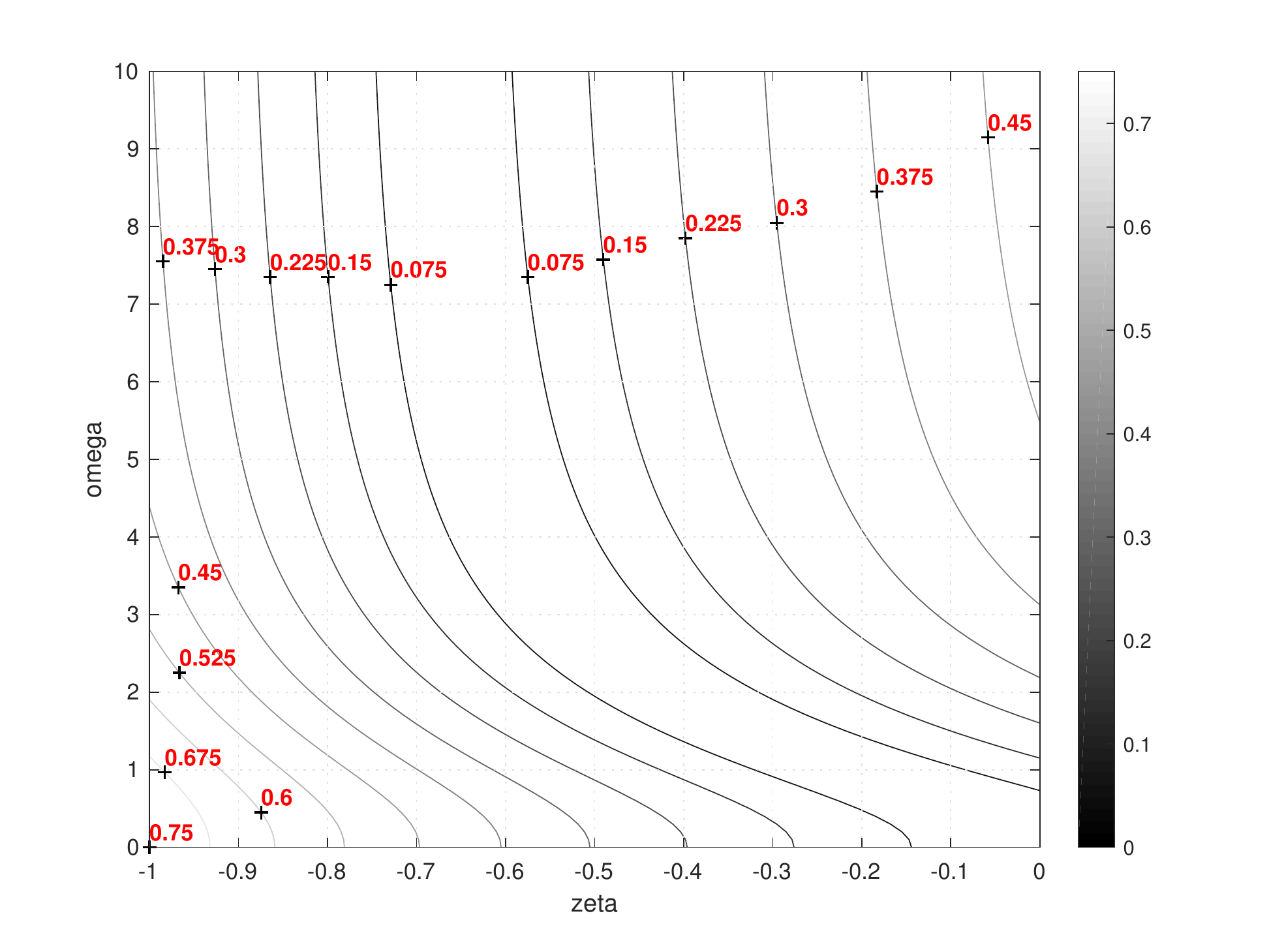}d)
	\includegraphics[width=0.45\textwidth]{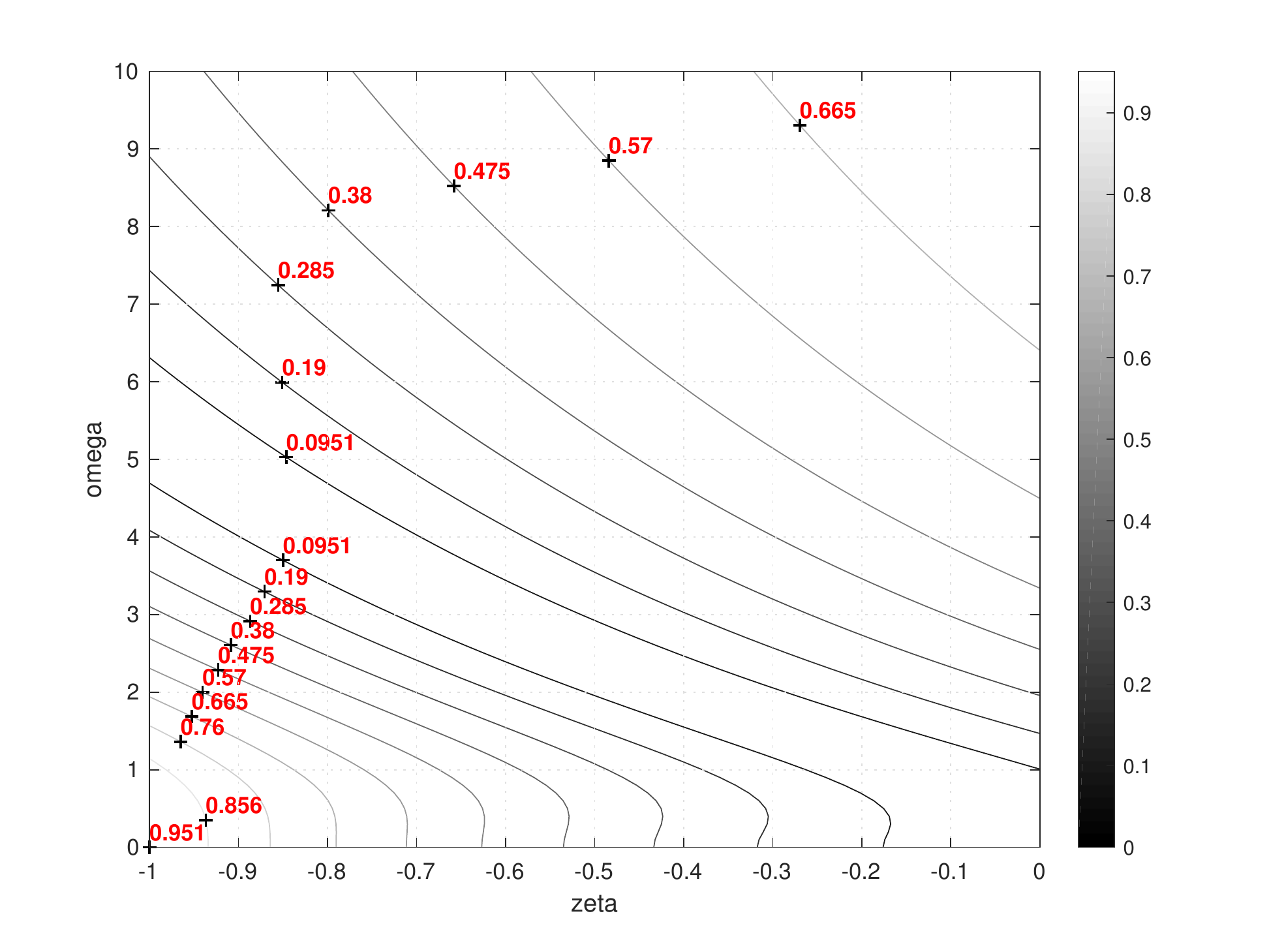}e)
	\includegraphics[width=0.45\textwidth]{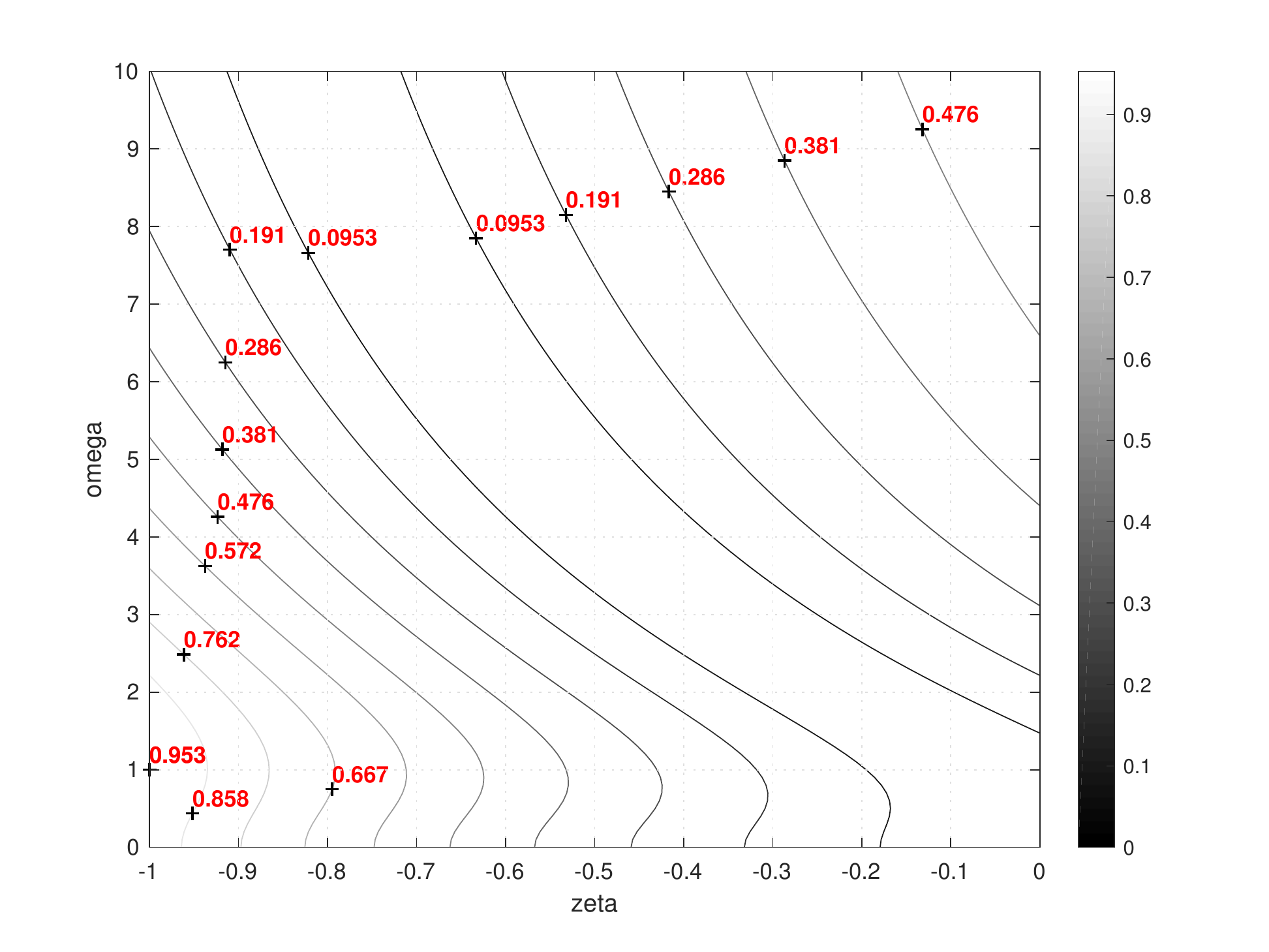}f)
	\end{center}
	\caption{Relative error in the spectral norm with respect to the exact evolution operator
	for a) TR-BDF2 method b) the Newmark method with $\rho_{\infty}=0 $   c) the Newmark method with $\rho_{\infty}=0.25 $ 
	 d) the Newmark method with $\rho_{\infty}=0.5 $  e) the CH-$\alpha$ method with $\rho_{\infty}=0 $ 
	  e) the CH-$\alpha$ method with $\rho_{\infty}=0.25.$ }
	\label{fig:relerrors_methods}
\end{figure}

The ratios of the spectral norms of the discrete evolution operators to that of the exact evolution operator
 are reported  for different methods in  Figure \ref{fig:ratio_methods} as a function of $(\zeta, \omega),$
 on the domain $(\zeta, \omega) \in [-1,0]\times [0,10]. $ Notice that these quantities are symmetric with respect to the
 $\omega=0 $ axis.
The corresponding
relative errors in the spectral norm  with respect to  the exact evolution operators are reported  
 in Figure \ref{fig:relerrors_methods}. In particular,  the TR-BDF2 method, the Newmark method with three different
 values of the damping parameters and the G($\alpha $) method introduced in \cite{chung:1993} (denoted by CH-$\alpha$)
 with two different
 values of the damping parameters
 are compared in this way.  
 
It can be observed that the TR-BDF2 method introduces significantly less numerical damping 
in the ranges corresponding to frequencies that need to be accurately resolved, while damping is  more effective on the highest frequencies.
Furthermore, the  TR-BDF2 method  also introduces a much smaller error than its counterparts in the same frequency ranges. These findings,
which could have been obtained in the undamped case with an analysis analogous to that presented in \cite{owren:1995},
are  also true if damping is added, thus confirming the advantages of the TR-BDF2 method for realistic applications including 
possibly stiff dissipative terms.
\newpage

\section{Numerical experiments}
\label{tests} \indent
A number of numerical experiments have been carried out in order to assess 
  the accuracy and efficiency of the TR-BDF2 method for applications to
structural mechanics and to compare it with that of the G($\alpha$) methods.
We first consider one  of the nonlinear numerical benchmarks discussed in \cite{erlicher:2002}.
Then, time discretizations of the wave equation are considered, in cases which are representative
of possible applications to structural mechanics and seismic engineering.
 For all the problems concerning the discretization of a wave equation, we have used  a finite element  spatial discretization and our computations have been performed with the free software FreeFem++  \cite {hecht:2012}, computing
 the errors in the space-time norms 
$$ \|u\|_{L^\infty(L^2)} =\displaystyle \max_n \|u(t^n, \cdot) \|_{L^2(\Omega)} $$
$$ \|u\|^2_{L^2(H^1)} = \displaystyle \sum_n \|u(t^N, \cdot) \|^2_{H^1(\Omega)}\Delta t $$
$$ \|u\|_{L^\infty(L^\infty)} =\displaystyle \max_n \|u(t^n, \cdot) \|_{L^\infty(\Omega)}, $$
where $t^n$ are the time levels used by the time  discretization on $[0,T].$

\subsection{Nonlinear system with 2 degrees of freedom}
 \label{2dofs}
\begin{figure}
  \begin{center}
	\includegraphics[width=0.45\textwidth]{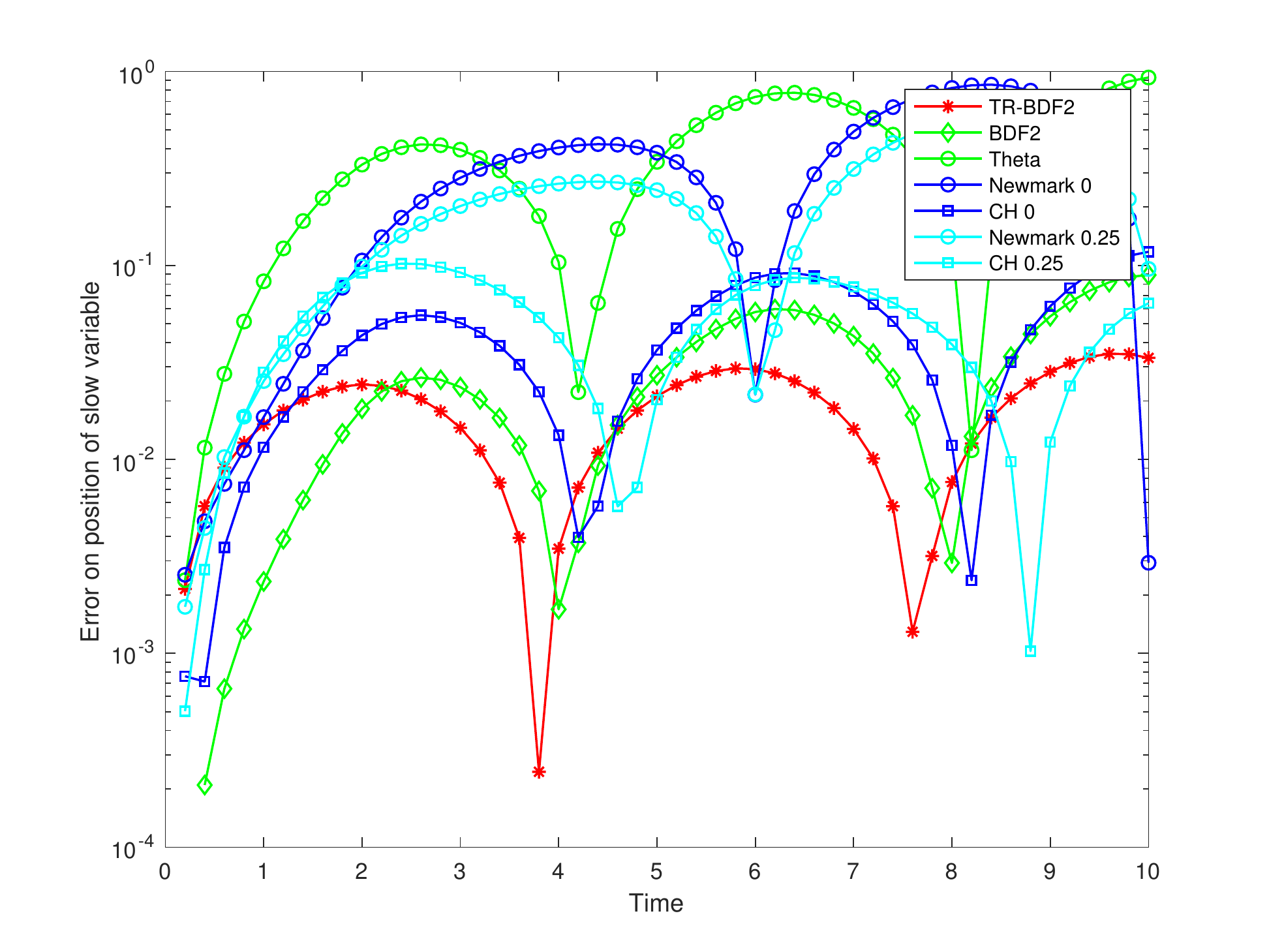}a)
	\includegraphics[width=0.45\textwidth]{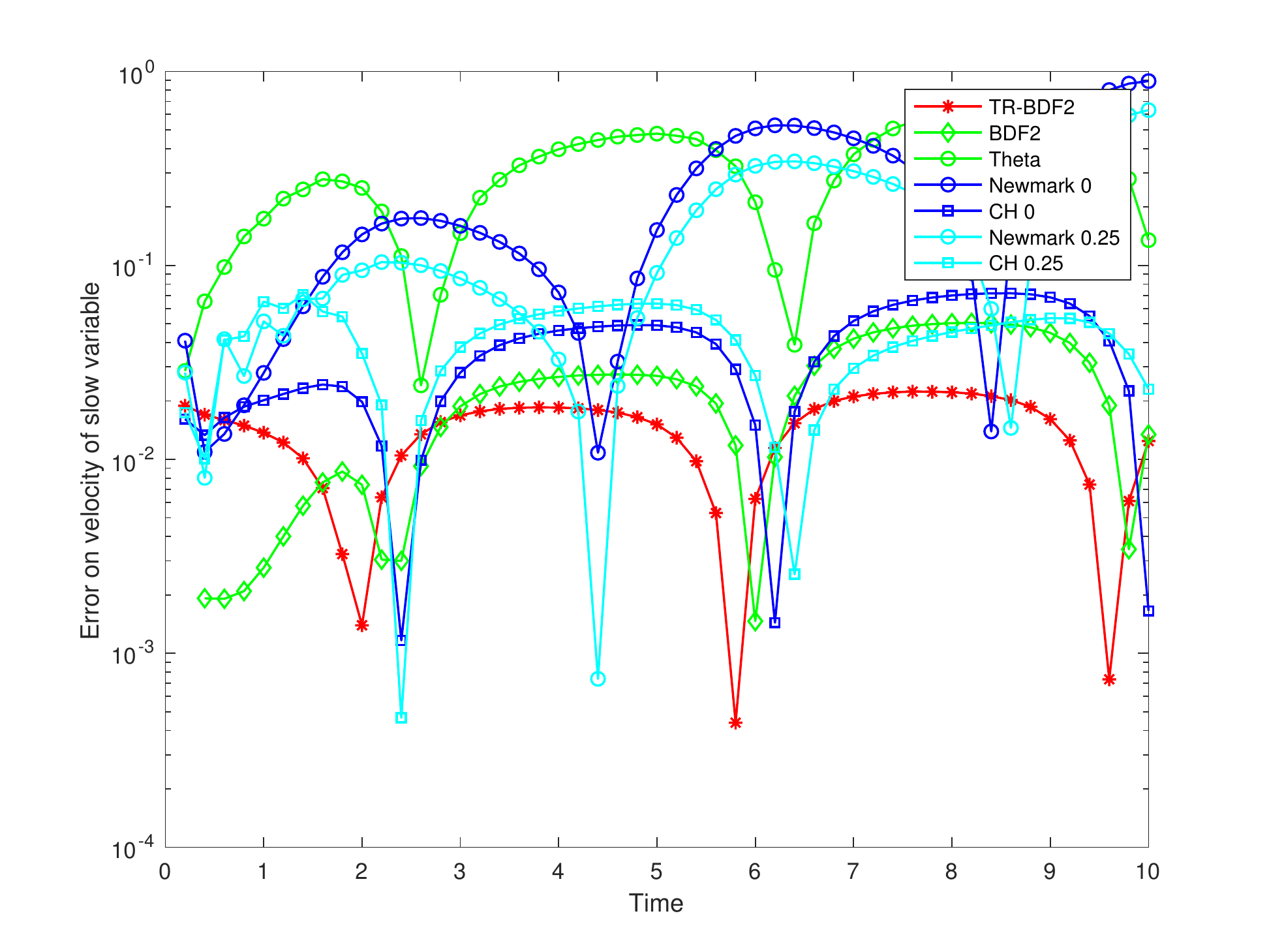}b)
	\includegraphics[width=0.45\textwidth]{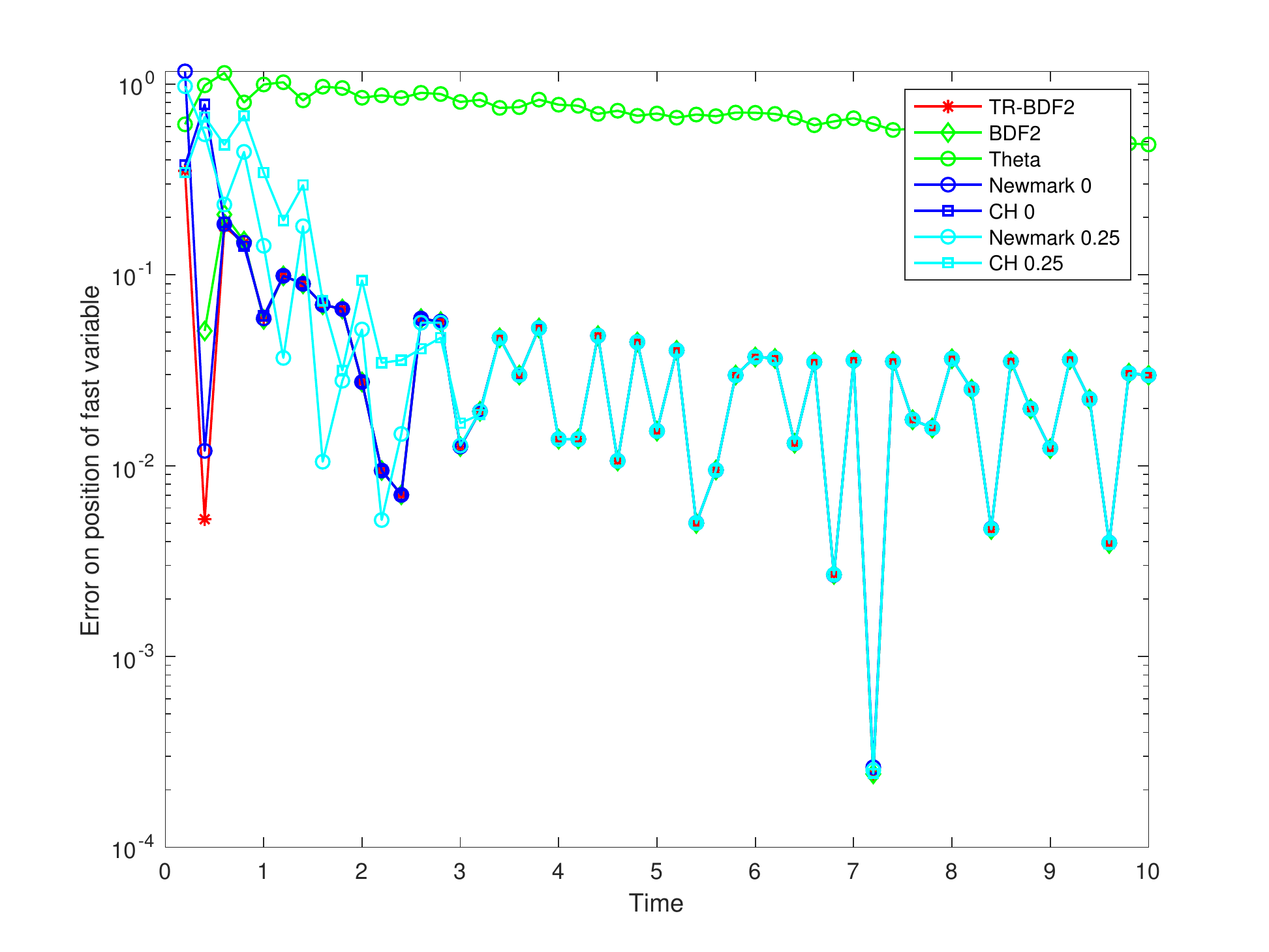}c)
	\includegraphics[width=0.45\textwidth]{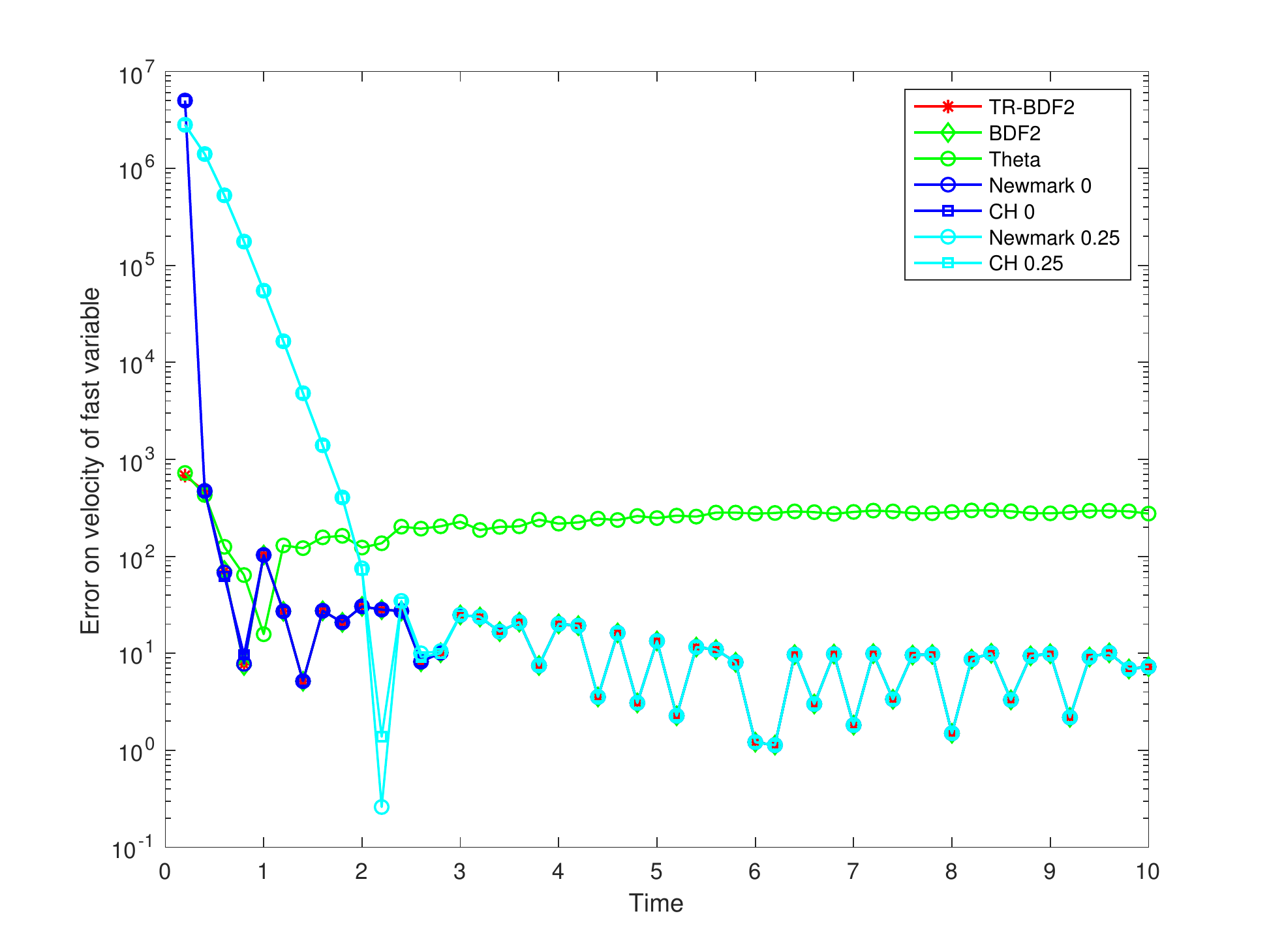}d) 
	\end{center}
	\caption{Absolute errors in the  test with 2 degrees of freedom on 
	  a) position for the slow degree of freedom b) velocity for the slow degrees of freedom   
	   c) position for the fast degree of freedom b) velocity for the fast degrees of freedom, as computed by TR-BDF2 and several other methods.}
	\label{fig:errors_2dofs}
\end{figure}

 In a first numerical experiment, we consider the strongly nonlinear system with two degrees of freedom
  \begin{equation}
 \label{eq:nonl2dofs}
\left[
  \begin{array}{c}
y_1\\
y_2\\
\end{array}
  \right]^{\prime\prime}
=  - \left[
  \begin{array}{c}
10^4y_1(1+10^4y_1^2)-\tanh{(y_2-y_1)}\\
\tanh{(y_2-y_1)}\\
\end{array}
  \right] 
\end{equation}
with the initial conditions ${\bf  y}(0) =[1,1.5]^T  \ \ \ \ {\bf  y}^{\prime}(0) = \mathbf{0}.$ The same system was used in
\cite{erlicher:2002} to compare the performance of  different G($\alpha$) methods.
 We compute a reference solution using the MATLAB {\tt ode15s} solver and we compare the performance of the TR-BDF2 method,
 of its parent methods, i.e. the off-centered Trapezoidal Rule or $\theta-$ method (setting $\theta=0.51$)  and the BDF2 method, and of
 4 different G($\alpha$) methods, more specifically the Newmark method with $\rho_{\infty}=0, 0.25 $ and the CH-$\alpha $ method
 with the same values of the dissipation parameter.
 As it can be seen in figure \ref{fig:errors_2dofs}, while the performance of the G($\alpha$) methods and of the TR-BDF2 method
 is essentially analogous on the the fast and strongly damped variable, the latter is in general the most accurate on the slow variables.

 \newpage

\subsection{Wave equation in 1 dimension with strongly varying coefficients}
 \label{wave1d}
\begin{figure}
  \centering
    \includegraphics[width=1\textwidth]{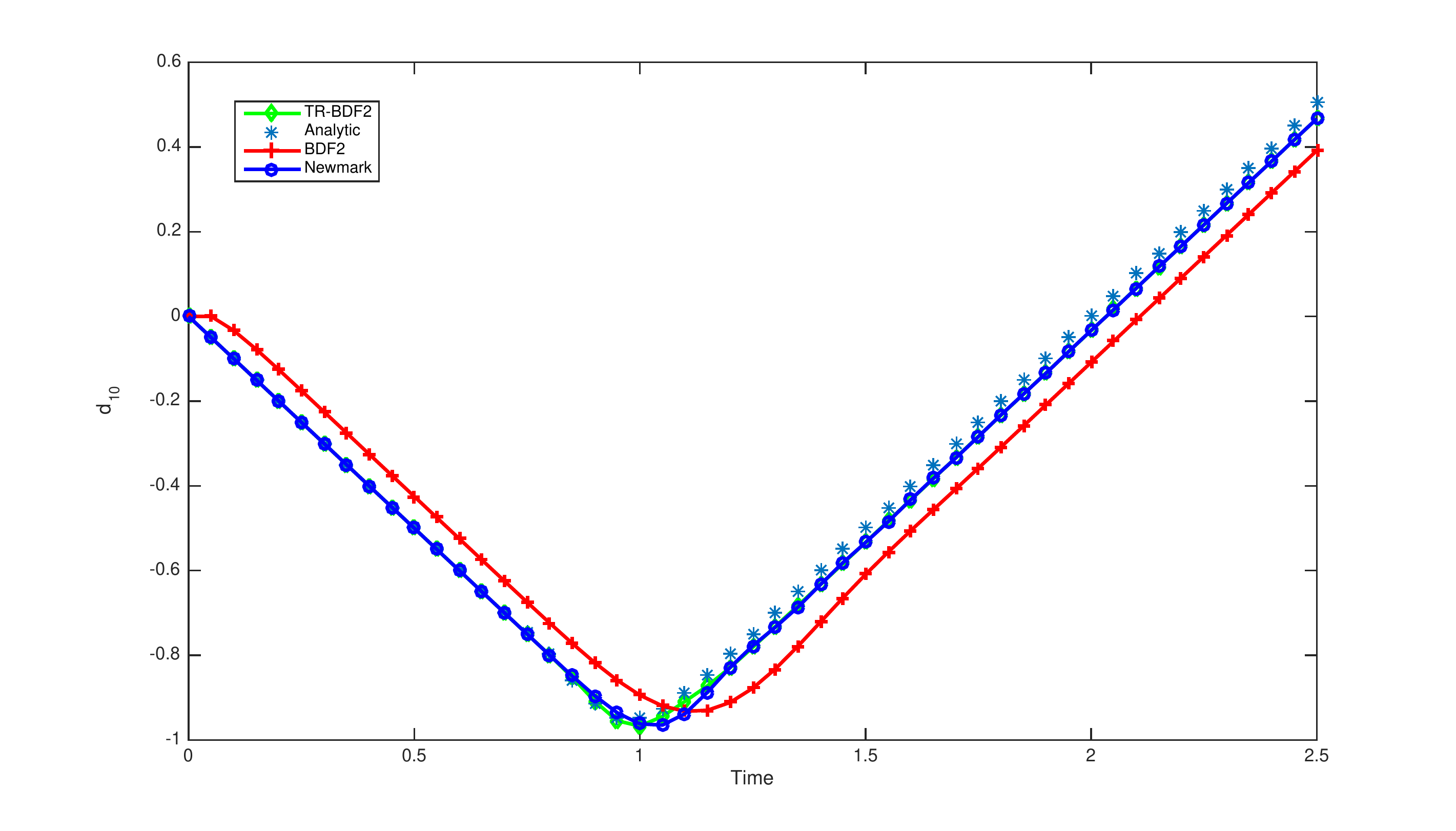}
  \caption{Displacement in x=10.5}
  \label{fig:d10}
\end{figure}

\begin{figure}
  \centering
    \includegraphics[width=1\textwidth]{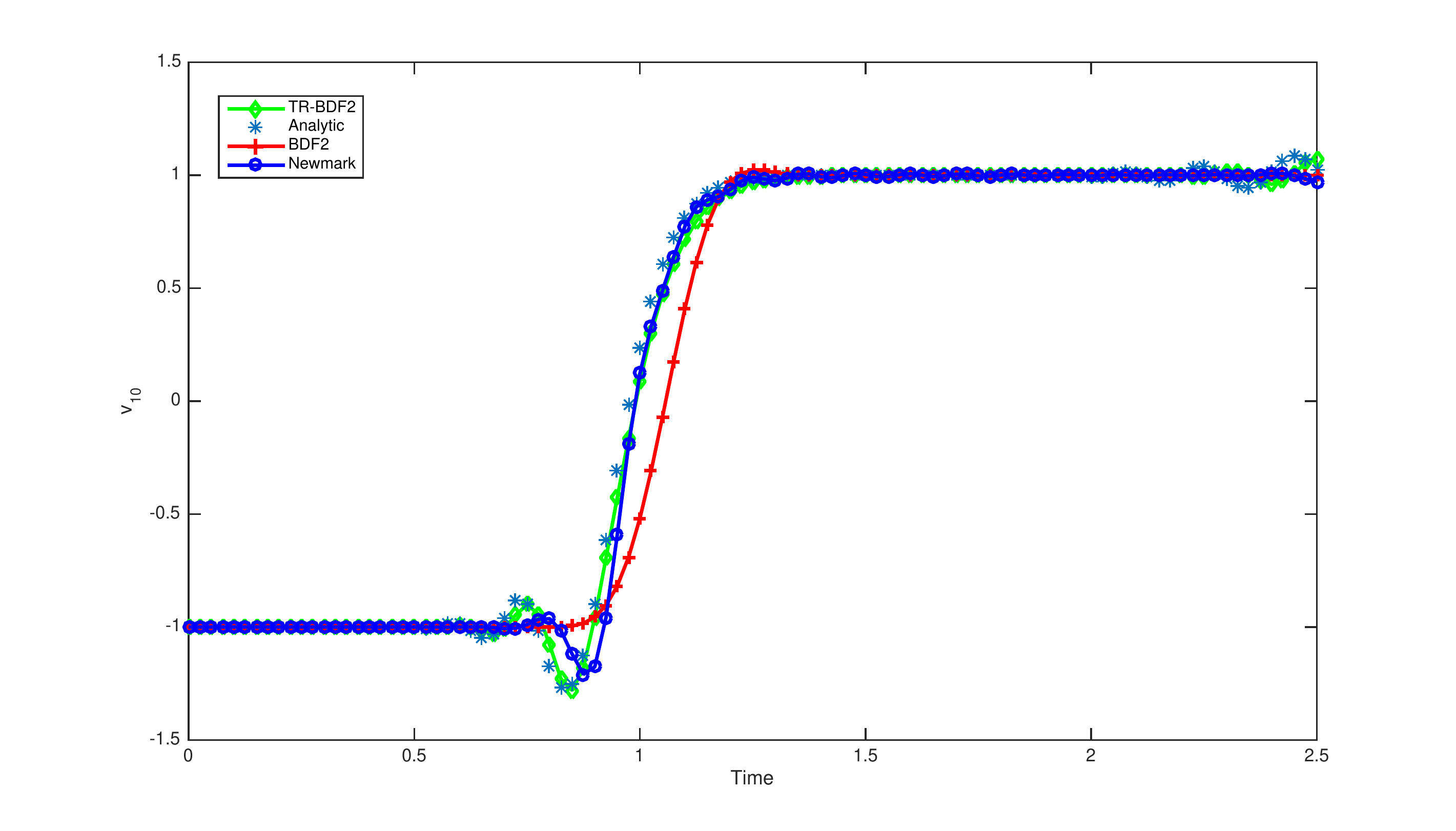}
  \caption{Velocity in x=10.5}
  \label{fig:v10}
\end{figure}

\begin{figure}
  \centering
    \includegraphics[width=1\textwidth]{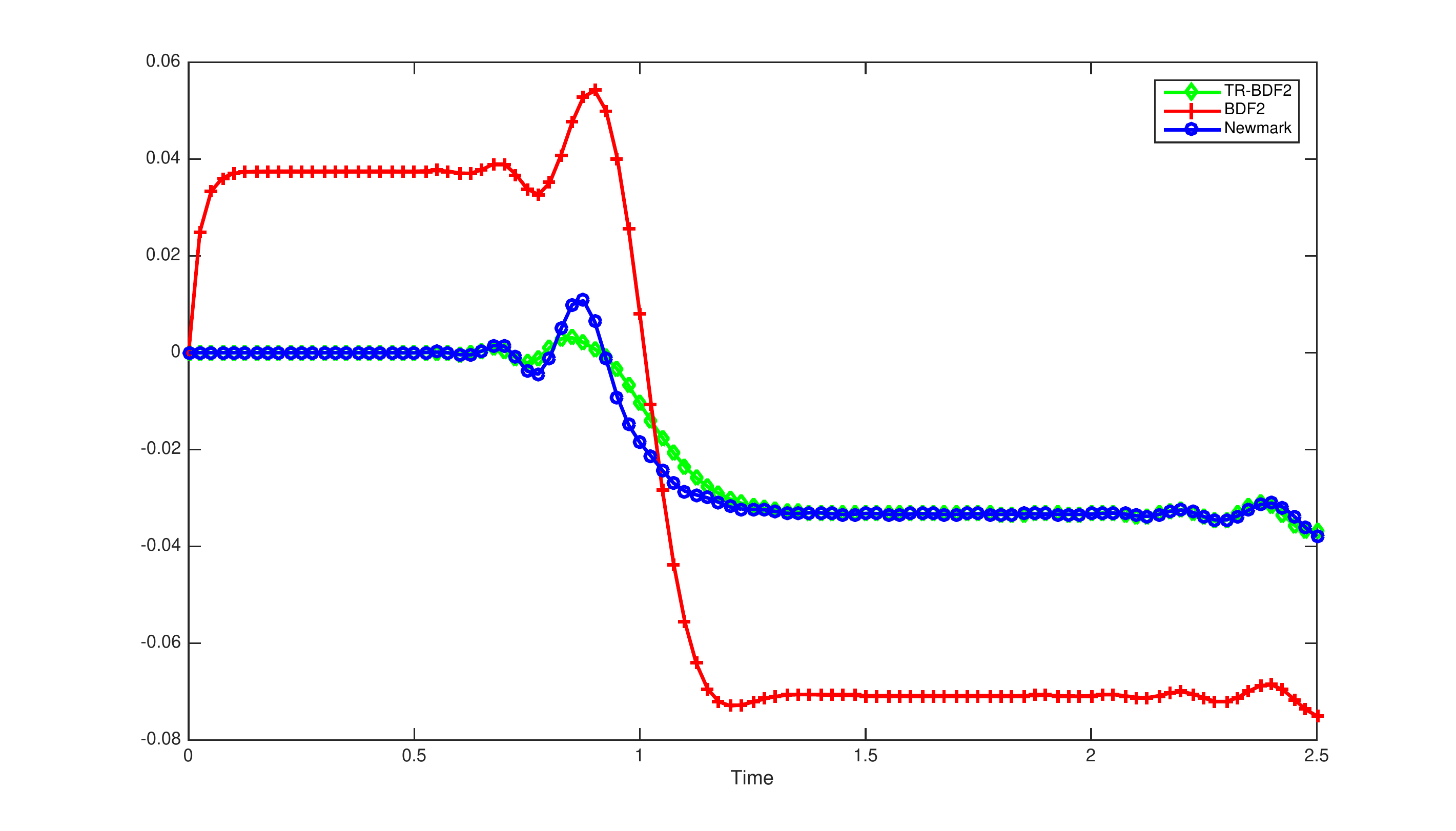}
  \caption{Displacement difference from Analytic solution in x=10.5}
  \label{fig:dd10}
\end{figure}

\begin{figure}
  \centering
    \includegraphics[width=1\textwidth]{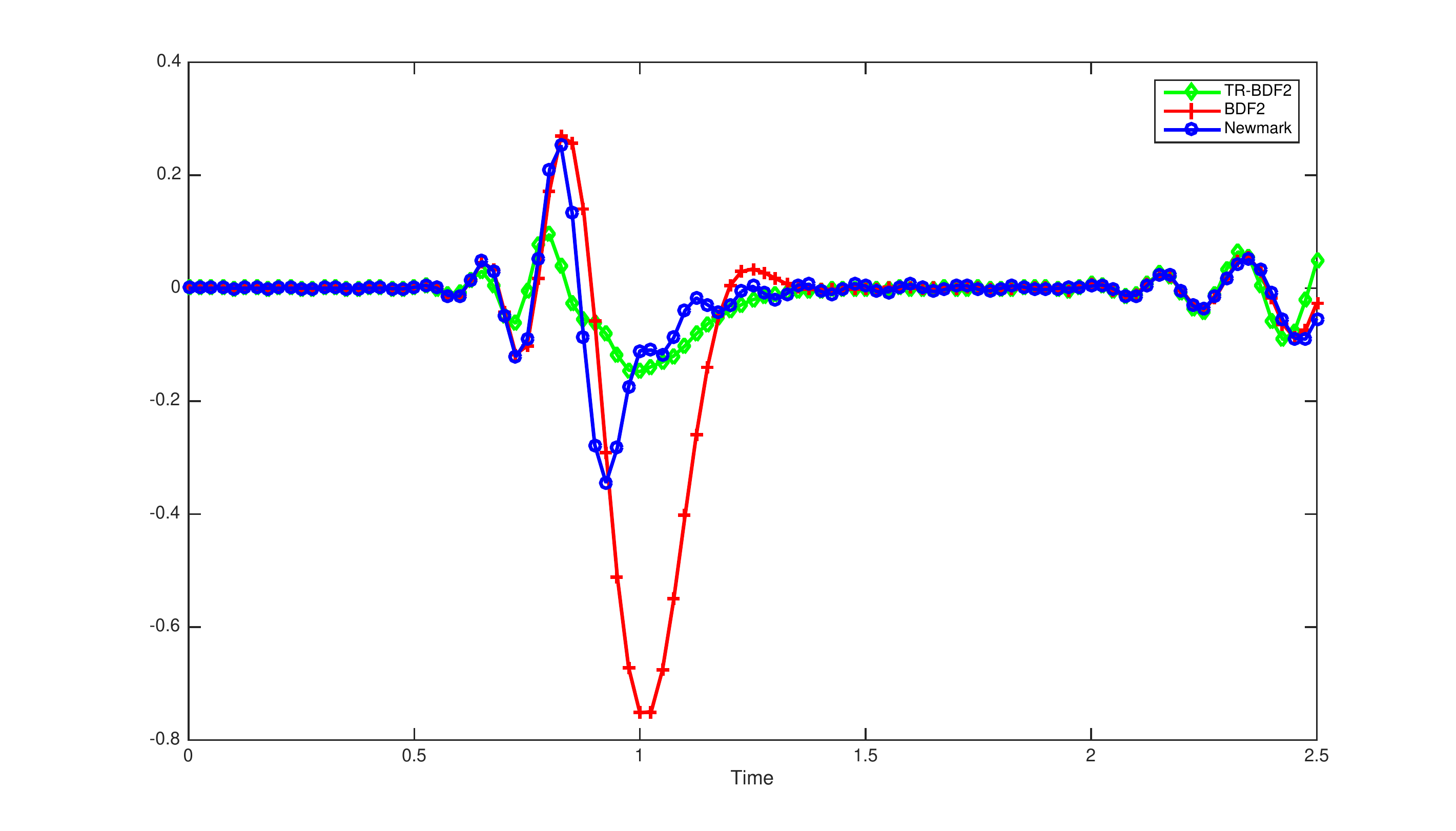}
  \caption{Velocity difference from Analytic solution in x=10.5}
  \label{fig:vv10}
\end{figure}

 We then consider the one dimensional wave propagation problem described by
 \begin{equation}
 \left\{
 \begin{array}{l}
 \displaystyle  \rho  \frac{\partial^2 u}{ \partial t^2}-\frac{\partial}{\partial x} \left( E(x) \frac{\partial u}{\partial x} \right)=0, \quad \mbox { \rm{in} } \,  [0,L] \times [0,T] \\ 
 \medskip
 u(0,t)=0, \quad \mbox { \rm{on} } \,  [0,T] \\
 \medskip
E(L)\displaystyle \frac{\partial u} {\partial x}(L,T)= 0, \quad \mbox { \rm{on} } \,  [0,T] \\
\medskip
u(x,0)=0, \quad \mbox { \rm{in} } \,  [0,L] \\
\medskip
\displaystyle \frac{\partial u}{\partial t}=-1, \quad \mbox { \rm{in} } \, [0,L]. 
 \end{array}
 \right.
 \label {hl78}
 \end{equation}
  These equations model the displacement $u$ of a clamped-free elastic rod with a unit constant cross-sectional area, 
   length $L, $   constant density $\rho$ and Young modulus given by $E(x).$
  The same problem has also been considered  in  \cite{hughes:1978}, \cite{piche:1995}, 
 \cite{erlicher:2002} to study the  performance of  G($\alpha$) methods
 in a stiff case in which strong oscillations and numerical overshoots can arise.   In particular, we assume $L=10.5$,  
 $\rho=0.01 $ and a Young modulus such that
 $$
 E(x)= \left \{  \begin{array}{ll} 10^7 & 0\leq x < 0.5 \\ 10^2 & 0.5 \leq x < 10 \\ 10^7 & 10 \leq x \leq 10.5,  \end{array} \right.
 $$ 
  which yields a rather  stiff problem.
 The system is discretized  using  $\mathbb{P}_1$ finite elements on
the spatial domain $[0,L] $ and either the TR-BDF2 method, the BDF2 method or the Newmark method
with  
$\beta=1/4 ,  $  $\gamma_N = 1/2$ for the time discretization. In all cases we have used a mesh with 21 nodes and $\Delta t=0.025 \rm s. $  
As a reference solution, we consider the one obtained by computing the exponential of the matrix of the 
spatial semi-discretization. 
In table \ref{tabla:com} we report the absolute errors between the considered methods and the reference solution,
computed at two different final times. It can be observed that the TR-BDF2 method  yields consistently the smallest errors
among the three second order methods under comparison.
\begin{table}[htbp]
\begin{center}
\begin{tabular}{| l | l | l | l | l |}
\hline
$T$ & Method & Error $L^\infty(L^2)$  &Error $L^2(H^1)$  & Error $L^\infty(L^\infty)$ \\ \hline \hline
1 & BDF2 &  6.63e-2&0.59 &  2.50e-2  \\ \hline
1& Newmark & 2.46e-2 &0.19 &  6.99e-3\\ \hline
1 &TR-BDF2 &  1.51e-2 &  6.00e-2&  2.67e-3 \\ \hline \hline
2.5 & BDF2 & 0.32 & 2.20 &  5.41e-2\\ \hline
2.5 & Newmark & 0.15 & 1.06 & 2.51e-2 \\ \hline
2.5 & TR-BDF2 & 0.14  & 0.92 &  7.73e-3  \\ \hline \hline
\end{tabular}
\caption{Absolute errors in the approximation of a stiff 1D wave equation}
\label{tabla:com}
\end{center}
\end{table}
     The time series of the displacements and velocities corresponding to the end of the rod ($x=L$) 
     are reported instead in Figures \ref{fig:d10}, \ref{fig:v10}, while the corresponding absolute differences between
     the computed solutions and the reference one are displayed in Figures \ref{fig:dd10} and \ref{fig:vv10}.
     Again, the smaller error produced by the TR-BDF2 method is apparent.
  

\subsection{Two-dimensional test with analytic solution}
\label{wave2d}
In order to verify the implementation used in the more advanced tests presented in Section \ref{elasticity2d},
we have considered the the problem
\begin{equation}
\left\{
\begin{array}{ll}
\displaystyle \frac{\partial^2 u}{\partial t^2} -2\Delta u =0 & {\mbox{ \rm{ in }}} \Omega \times [0,T], \\
u(x,y,0)=0 &  {\mbox{ \rm{ in }}} \Omega , \\
\displaystyle \frac{\partial u}{\partial t }(x,y,0)=2\pi {\rm{ sin}}(\pi x){\rm{sin}}(\pi y) &  {\mbox{ \rm{ in }}} \Omega, \\
u(x,y,t)=0 &  {\mbox{ \rm{ on }}} \partial  \Omega \times [0,T], 
\end{array}
\right.
\label{t1}
\end{equation}
where the domain $\Omega=[0,1]\times [0,1]$ and $T=1$. The solution for this problem is given by
$u(x,y,t)={\rm{ sin}}(2 \pi t) {\rm { sin}}(\pi x){\rm { sin}}(\pi y).$ This solution has been used to perform a convergence study
of the TR-BDF2 time discretization, coupled again to a $\mathbb{P}_1$  finite element discretization in space.
To estimate empirically the convergence order, we choose the maximum diameter $h$  of a mesh element equal
 in size to the time step $\Delta t$ and we compute the following expression for two  different values of $\Delta t$:
\begin{equation}r_{emp}= \frac{\log(err(\Delta t_2)/err(\Delta t_1))} {\log(\Delta t_2 /\Delta t_1 )}. \label{fo} \end{equation}
The results of the convergence test are displayed in  table \ref{tabla:er}.
\begin{table}[htbp]
\begin{center}
\begin{tabular}{| l | l | l | l | l |}
\hline
$h=\Delta t$ &  $L^\infty(L^2)$  error & $L^2(H^1)$ Error &  $r_{emp}$ in $L^\infty(L^2)$ & $r_{emp}$ in $L^2(H^1)$\\\hline \hline
0.1 & 1.26e-2 &4.08e-2 & & \\ \hline
0.05 &  2.78e-3 & 8.79e-3 & 2.17 & 2.21 \\ \hline
0.025 & 6.39e-4&  2.04e-3 &2.12  &2.11 \\ \hline
0.0125 & 1.52e-4 &4.91e-4 & 2.07  &2.05 \\ \hline
0.00625 & 3.71e-5 &1.21e-4 & 2.04 &2.03\\ \hline
0.003125 & 9.15e-6 &2.99e-5 & 2.02 & 2.01  \\ \hline
\end{tabular}
\caption{Convergence behaviour of the TR-BDF2 method in test \eqref{t1}, absolute errors 
and convergence rates in $L^\infty(L^2) $ and $L^2(H^1)$ norms. }
\label{tabla:er}
\end{center}
\end{table}
%

\begin{figure}
  \centering
    \includegraphics[width=1\textwidth]{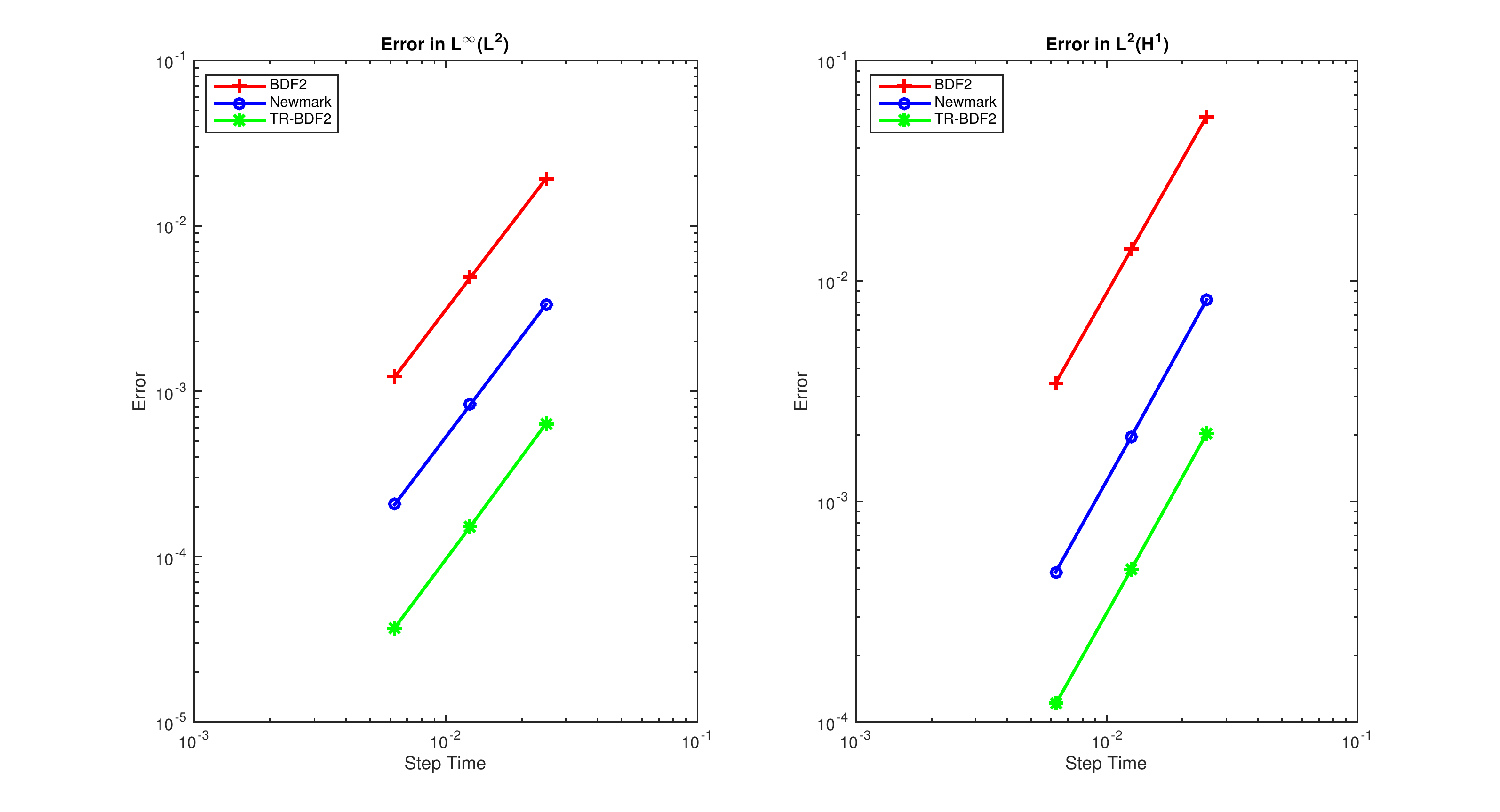}
  \caption{Comparison of convergence behaviour of different methods in test \eqref{t1}.}
  \label{fig:c1}
\end{figure}

We have also repeated the computation using the  BDF2 method and the  Newmark  method with 
$\beta= 1/4,  $  $\gamma_N = 1/2$. The results are shown in table \ref{tabla1:com} and figure \ref{fig:c1}.
While second order convergence is achieved by all methods, the results 
  clearly display the significantly smaller errors obtained by the TR-BDF2 method.  
  
\begin{table}[htbp]
\begin{center}
\begin{tabular}{| l | l | l | l | }
\hline
$h=\Delta t$ & Method & $L^\infty(L^2)$ error  &$L^2(H^1)$ error \\ \hline \hline
0.025 & BDF2 & 1.93e-2 & 5.58e-2\\ \hline
0.025 & Newmark & 3.35e-3 & 8.22e-3  \\ \hline
0.025 &TR-BDF2 &  6.39e-4 & 2.04e-3 \\ \hline \hline
0.0125 & BDF2 &  4.87e-3 &  1.38 e-2 \\ \hline
0.0125 & Newmark &  8.32e-4 &  1.96e-3  \\ \hline
0.0125 & TR-BDF2 & 1.52e-4 &  4.91e-4 \\ \hline \hline
0.00625 & BDF2 &  1.22e-4 &  3.44e-3 \\ \hline
0.00625 & Newmark & 2.06e-4 & 4.78e-4 \\ \hline
0.00625 & TR-BDF2 &3.71e-5 &  1.2e-4 \\ \hline
\end{tabular}
\caption{Convergence behaviour of the BDF2, Newmark and TR-BDF2 method in test \eqref{t1},  absolute errors 
  in $L^\infty(L^2) $ and $L^2(H^1)$ norms.  }
\label{tabla1:com}
\end{center}
\end{table}
 
 \newpage

  \subsection{Two dimensional elasticity}
 \label{elasticity2d}
 
 We consider the plane-strain equations of two dimensional elasticity,  presented in \cite{leveque:2004} as  the first order system:
 
 \begin{eqnarray}\label{eq:elastic_1stord}
 \frac{\partial \sigma^{x,x}}{\partial t}&=&    (\lambda +2\mu)  \frac{ \partial u}{\partial x}+ \lambda\frac{ \partial  v }{\partial y}\nonumber\\
  \frac{\partial \sigma^{x,y}}{\partial t}&=& \mu \frac{ \partial v }{\partial x}  +\mu \frac{ \partial  u}{\partial y}  \nonumber\\
  \frac{\partial \sigma^{y,y}}{\partial t}&=&   \lambda \frac{ \partial u }{\partial x} +  (\lambda +2\mu)\frac{ \partial  v}{\partial y}\nonumber\\
    \rho\frac{\partial u}{\partial t}&=&  \frac{  \partial \sigma^{x,x}}{\partial x} + \frac{ \partial \sigma^{x,y} }{\partial y}\nonumber\\
     \rho\frac{\partial v}{\partial t}&=& \frac{ \partial \sigma^{y,x} }{\partial x} + \frac{ \partial  \sigma^{y,y}}{\partial y}.
\end{eqnarray}
  Here, $\sigma^{i,j} $ denote the components of the (symmetric) stress tensor  $ \boldsymbol{\Sigma}, $ 
  $\mathbf{u}=(u,v) $ is the velocity vector with components $u,v $ in the
  $x,y $ directions, respectively, $\rho $ is the medium density and $\lambda, \mu $ are the medium Lam\'e coefficients.
  For the purpose of our work, we reformulate these equations as second order equations in terms of the
  displacements $ \boldsymbol{\delta}=(\delta^x,\delta^y).$ As a result, equations \eqref{eq:elastic_1stord} are equivalent to
  
  \begin{eqnarray}\label{eq:elastic_1stord2}
 \sigma^{x,x} &=&    (\lambda +2\mu)  \frac{ \partial \delta^x }{\partial x}+ \lambda\frac{ \partial  \delta^y }{\partial y}\nonumber\\
  \sigma^{x,y}&=& \mu \frac{ \partial \delta^y }{\partial x}  +\mu \frac{ \partial  \delta^x }{\partial y}  \nonumber\\
\sigma^{y,y}&=&   \lambda \frac{ \partial \delta^x  }{\partial x} +  (\lambda +2\mu)\frac{ \partial \delta^x   }{\partial y}\nonumber\\
    \rho\frac{\partial^2 \delta^x }{\partial t^2}&=&  \frac{  \partial \sigma^{x,x}}{\partial x} + \frac{ \partial \sigma^{x,y} }{\partial y}\nonumber\\
     \rho\frac{\partial^2 \delta^y}{\partial t^2}&=& \frac{ \partial \sigma^{y,x} }{\partial x} + \frac{ \partial  \sigma^{y,y}}{\partial y},
\end{eqnarray}
which, after substitution of the expression for the stresses into the momentum equation, can be rewritten in vector form as 
   \begin{equation}\label{eq:elastic_1stord_vec}
     \rho\frac{\partial^2 \boldsymbol{\delta}}{\partial t^2}={\rm div}\boldsymbol{\Sigma} =  \mu \Delta \boldsymbol{\delta} +(\lambda+\mu) \nabla {\rm div} \boldsymbol{\delta}.
   \end{equation}
   As discussed in \cite{leveque:2004}, this system has P and S-wave solutions with propagation speeds
   \begin{equation}
   c_P =\sqrt{\frac{\lambda +2\mu}{\rho} } \ \ \ \ c_s=\sqrt{\frac{\mu}{\rho} }.
   \end{equation}
   
\begin{figure}
 \begin{center}
	\includegraphics[width=0.45\textwidth]{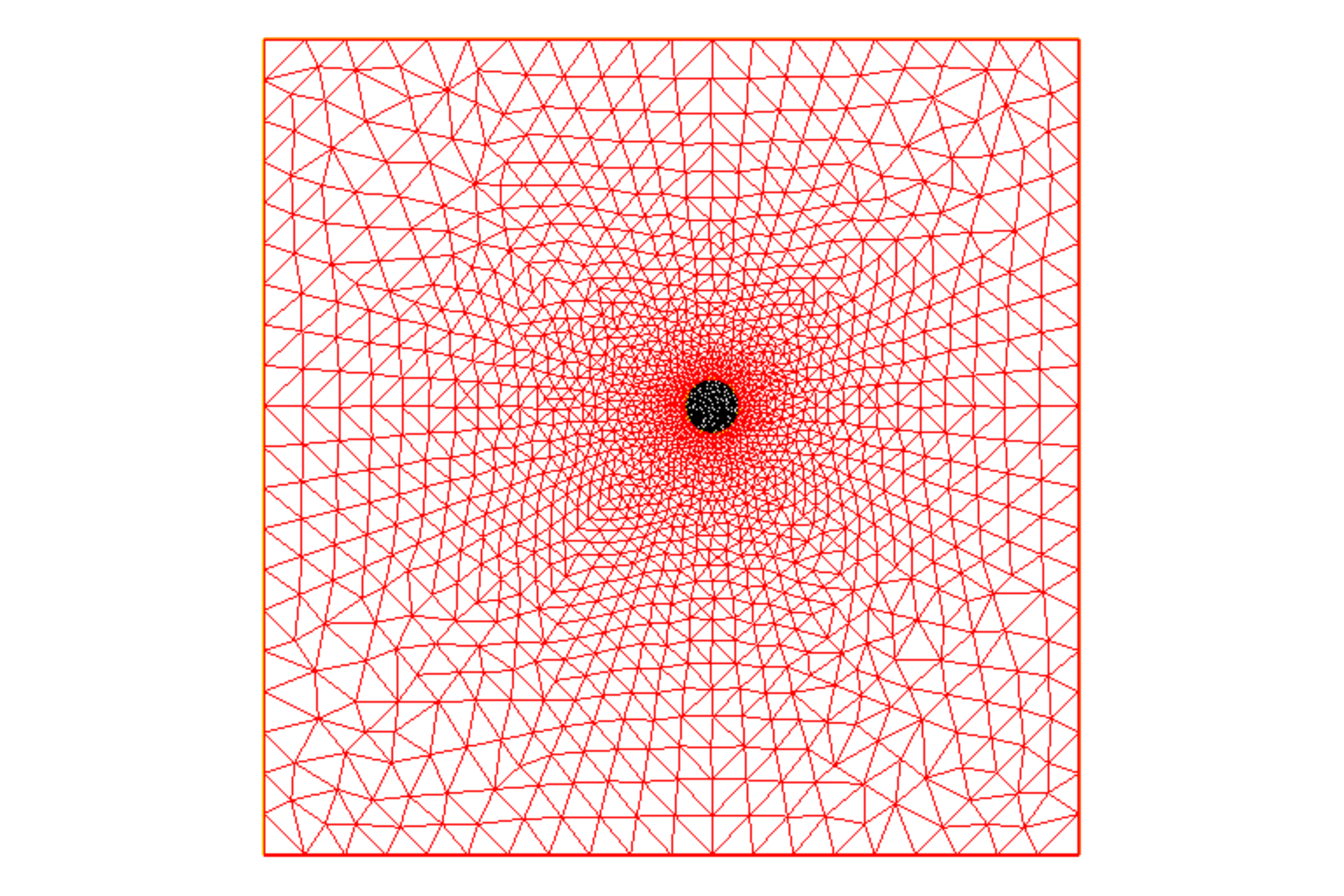}a)
	\includegraphics[width=0.45\textwidth]{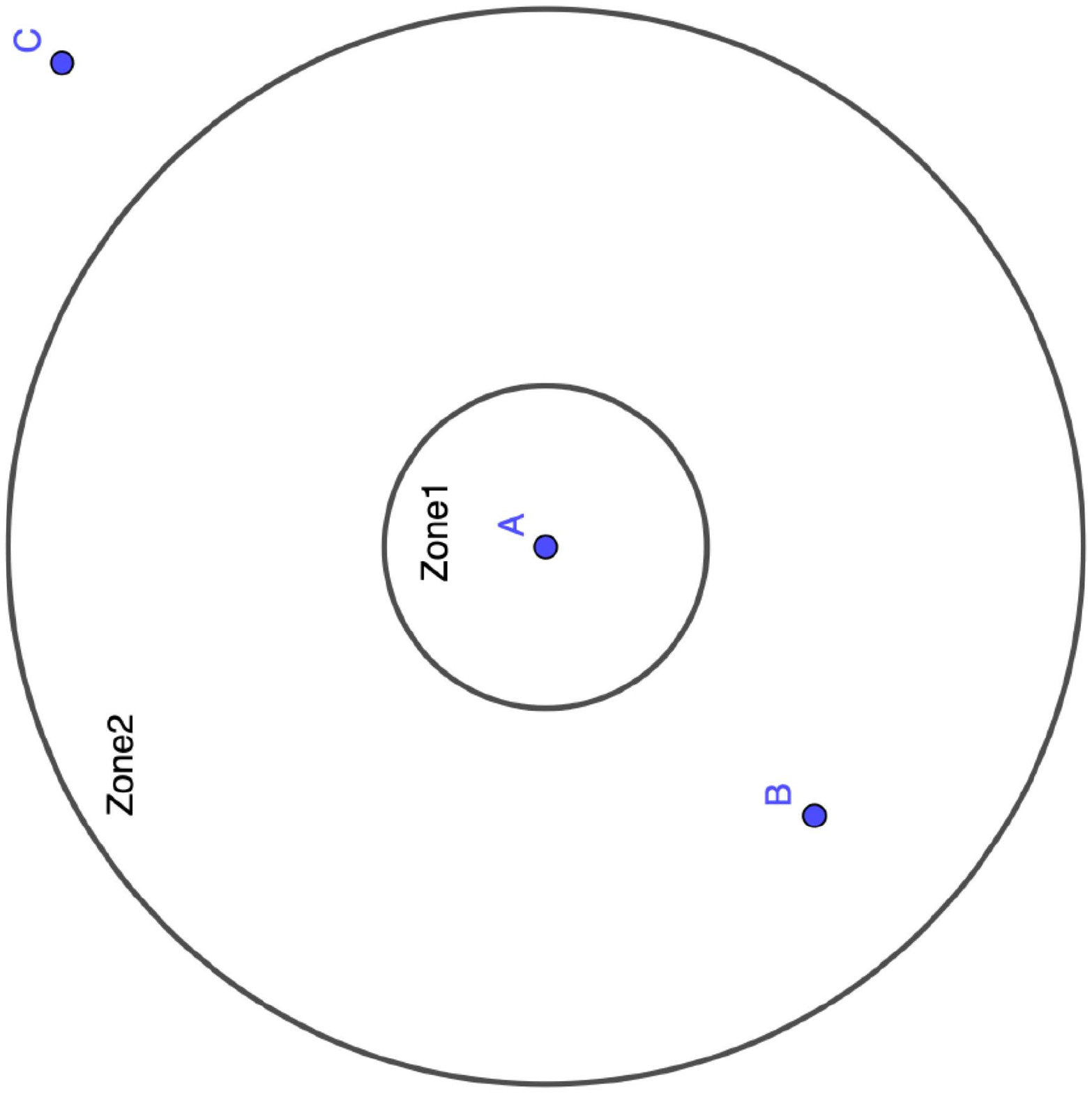}b)
	\end{center}
	\caption{Sketch of a) the computational domain and mesh  and b) the areas with more rigid material and initial perturbations.}
	\label{malla}
\end{figure}

   We consider the above equations on the domain  $\Omega=[0,3] \times [0,3] $ and on the time interval
   $[0,T], $ with $T=10^{-2} \ \rm s$. 
   We assume that the density is constant and equal one in all the domain, while
   the elastic constants are such that a more rigid inclusion is embedded   as shown in Figure \ref{malla}. 
   More specifically,  in the darker region, denoted as Zone 2 in Figure \ref{malla}b), whose boundary is described by the equation is $(x-x_0)^2+(y-y_0)^2 =1/100$, with $(x_0,y_0)=(1.65,1.65), $  we assume that  $\lambda=200, \, \mu=100,$ while in the region represented in red one has $\lambda=2$, $\mu=1. $ We then proceed to study the wave propagation resulting from imposing at the initial time 
    a   displacement   in the region denoted  as  Zone 1 in Figure \ref{malla}b),  whose boundary
    is described by the equation   $(x-x_0)^2+(y-y_0)^2 =9/10^4. $ 
    More specifically, we consider the function defined as $d(x,y)= \exp(-(x-x_0)^2/2-(y-y_0)^2/2) $  in Zone 1 and $d(x,y)=0 $ in the rest of the domain and we take as initial conditions
   $$ \delta^x=0, \, \delta^y= d(x,y), \, \frac{\partial \delta^x}{\partial t}=0, \,  \frac{\partial \delta^y}{\partial t}=v,$$
where   $v$ is defined by the finite difference approximation $ v(x,y)=d(x,y)/\Delta t, $  thus corresponding to the velocity that would be obtained if the impulsive displacement $ d(x,y) $ had been produced in a time interval of size $\Delta t $ starting from zero.
For the numerical discretization we consider $\mathbb{P}_1$ continuous finite elements on the locally refined unstructured mesh depicted in Figure \ref{malla} a), which is composed of 15036 vertices and 7559 elements with  an average diameter equal to    $h=2.39\times 10^{-2} $ in Zone 2 and equal to  $h=0.22 $ outside.  For the time discretization we use a time step of $\Delta t= 1.25\times 10^{-4}, $ which entails a   Courant number based on the  P-wave speed of $c_P\Delta t/h \approx 1.04 $ in Zone 1 and $c_P\Delta t/h \approx 0.01 $ in the rest of the domain, thus mimicking a situation in which   an unconditionally stable method is typically applied.

\begin{figure}
 \begin{center}
	\includegraphics[width=0.46\textwidth]{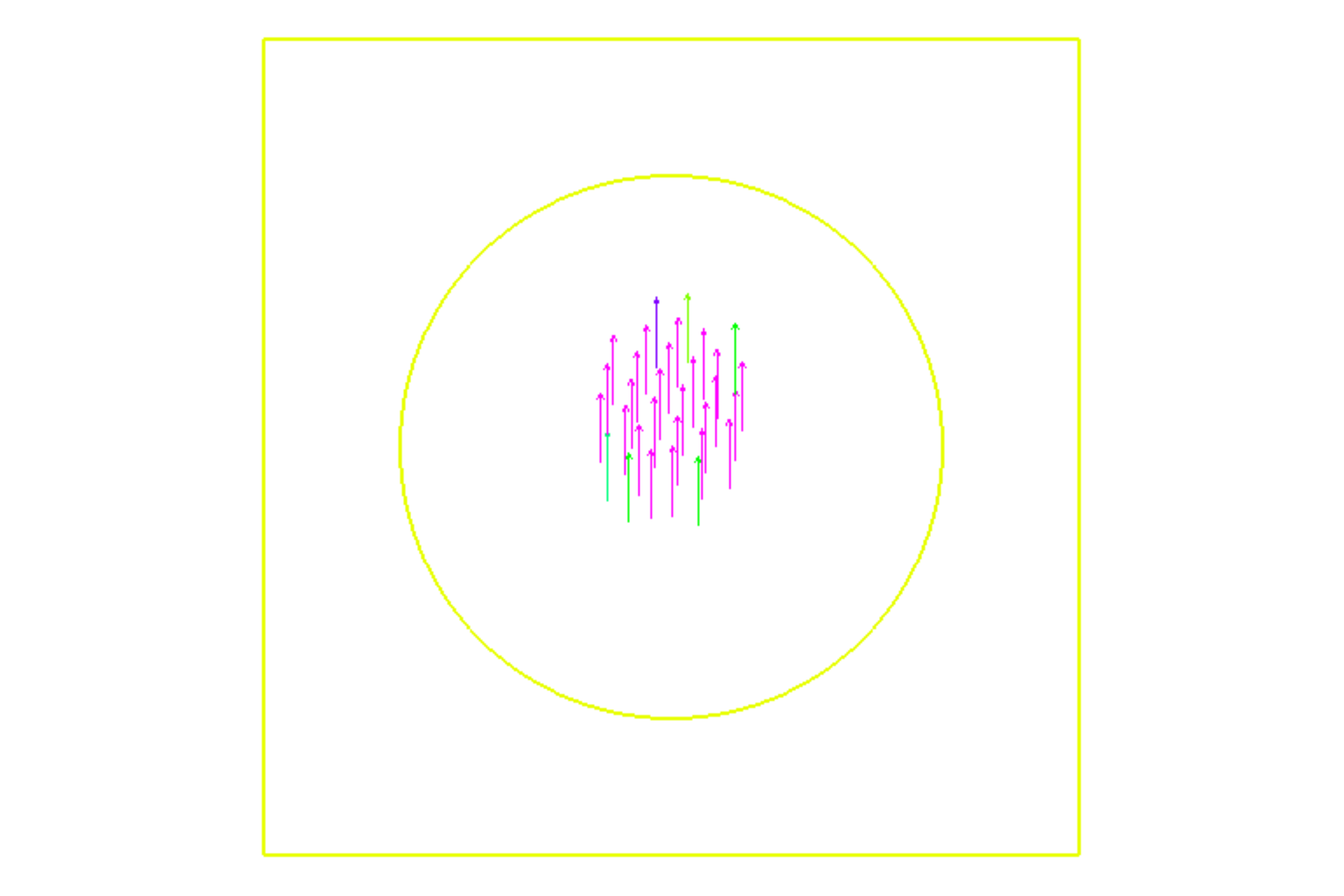}a)
	\includegraphics[width=0.46\textwidth]{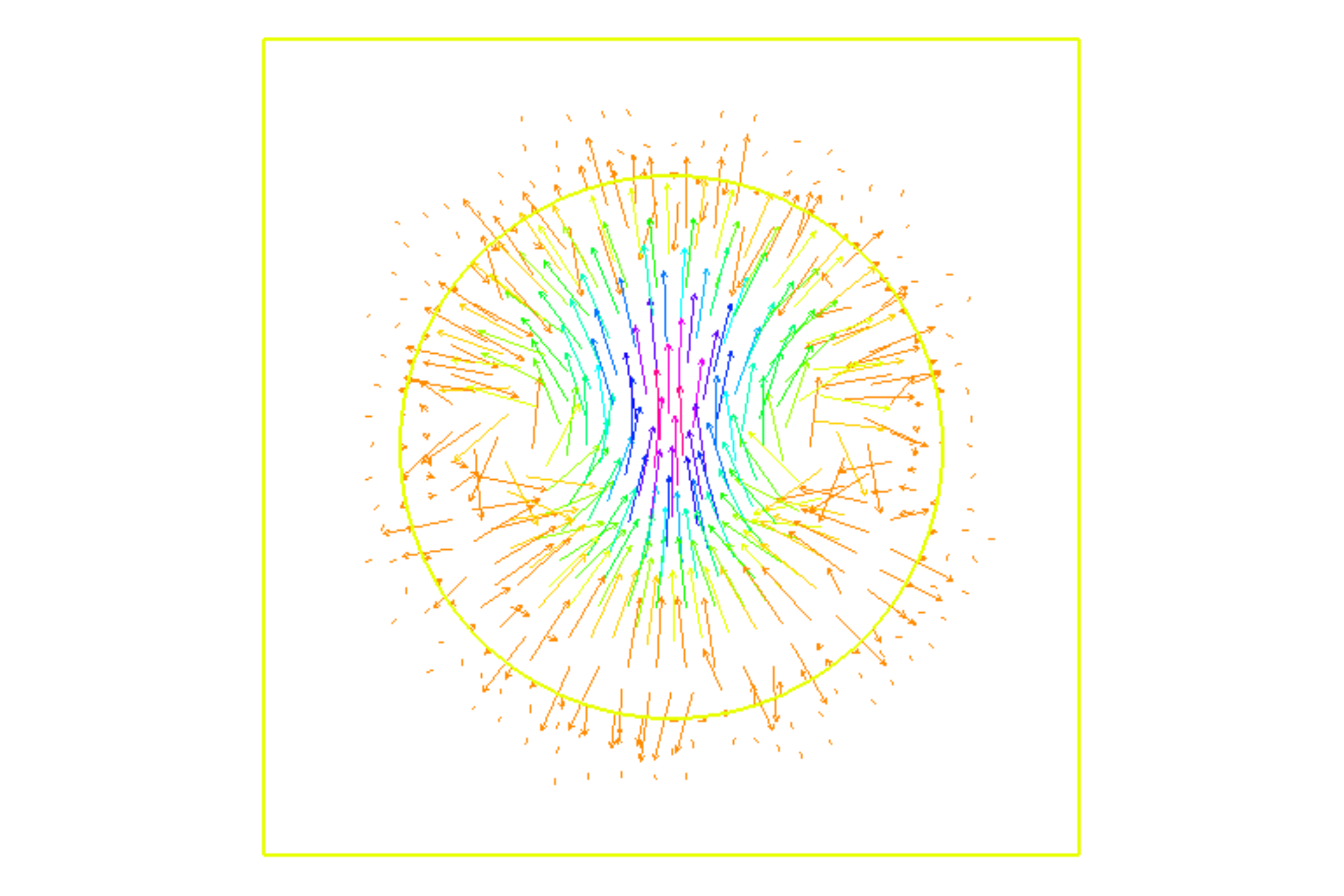}b)
	\includegraphics[width=0.46\textwidth]{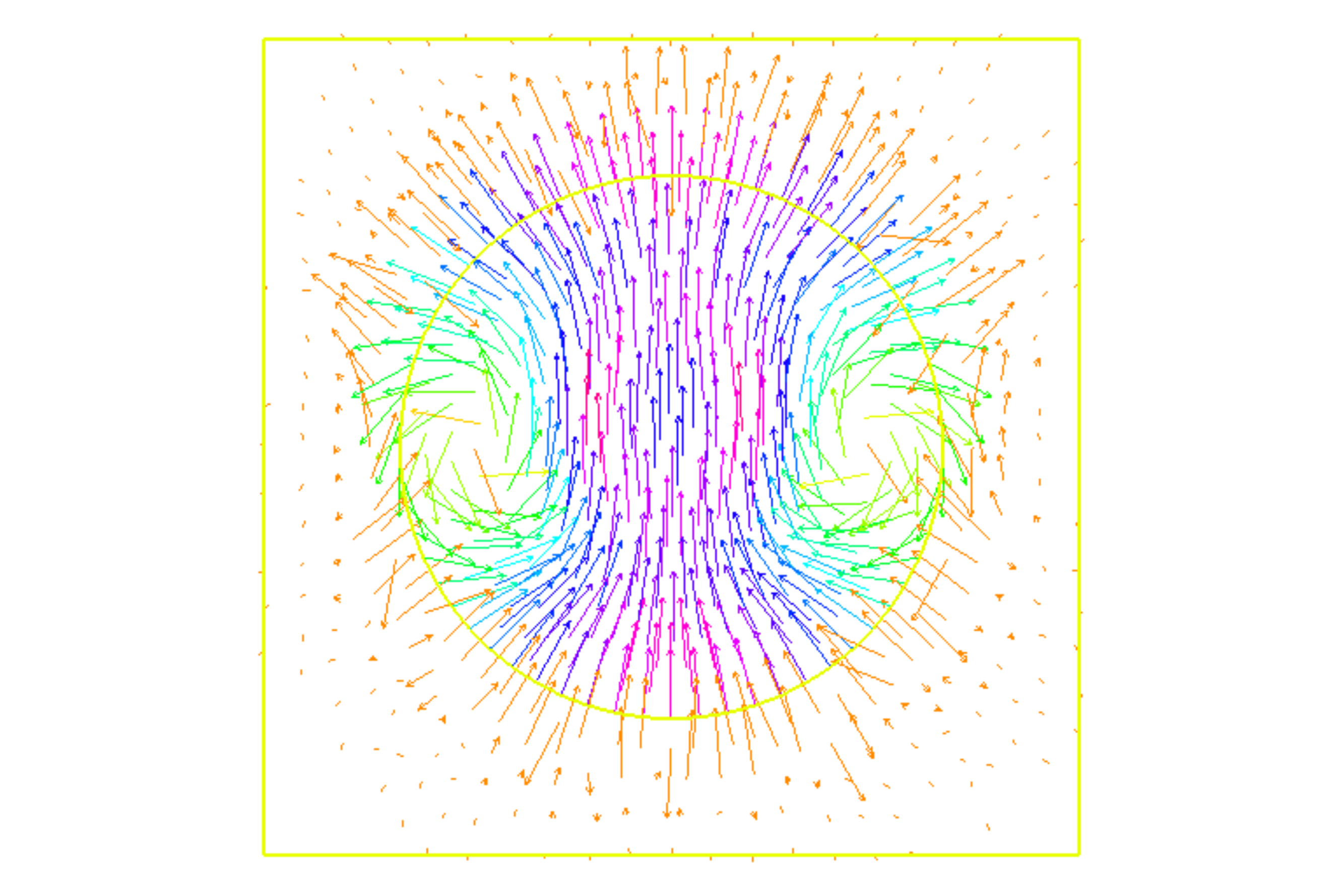}c)
	\includegraphics[width=0.46\textwidth]{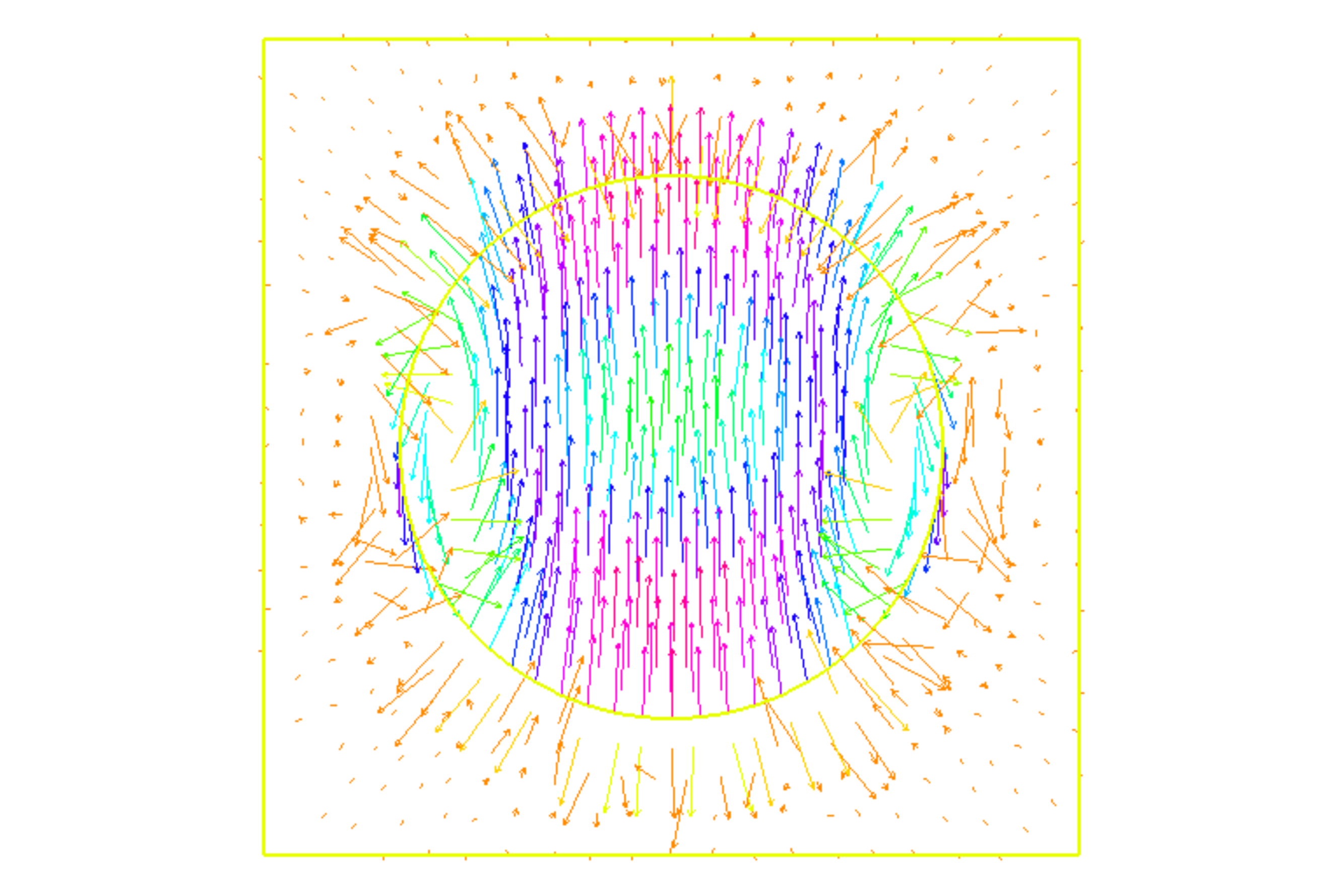}d)

	\end{center}
	\caption{Displacement vector field at times  a) $t=0, $ b) $t=2.5\times 10^{-3}, $ c)  $t=5\times 10^{-3}, $ d) $t=7.5\times 10^{-3}$.}
	\label{desplazamiento}
\end{figure}

\begin{figure}
 \begin{center}
	\includegraphics[width=0.45\textwidth]{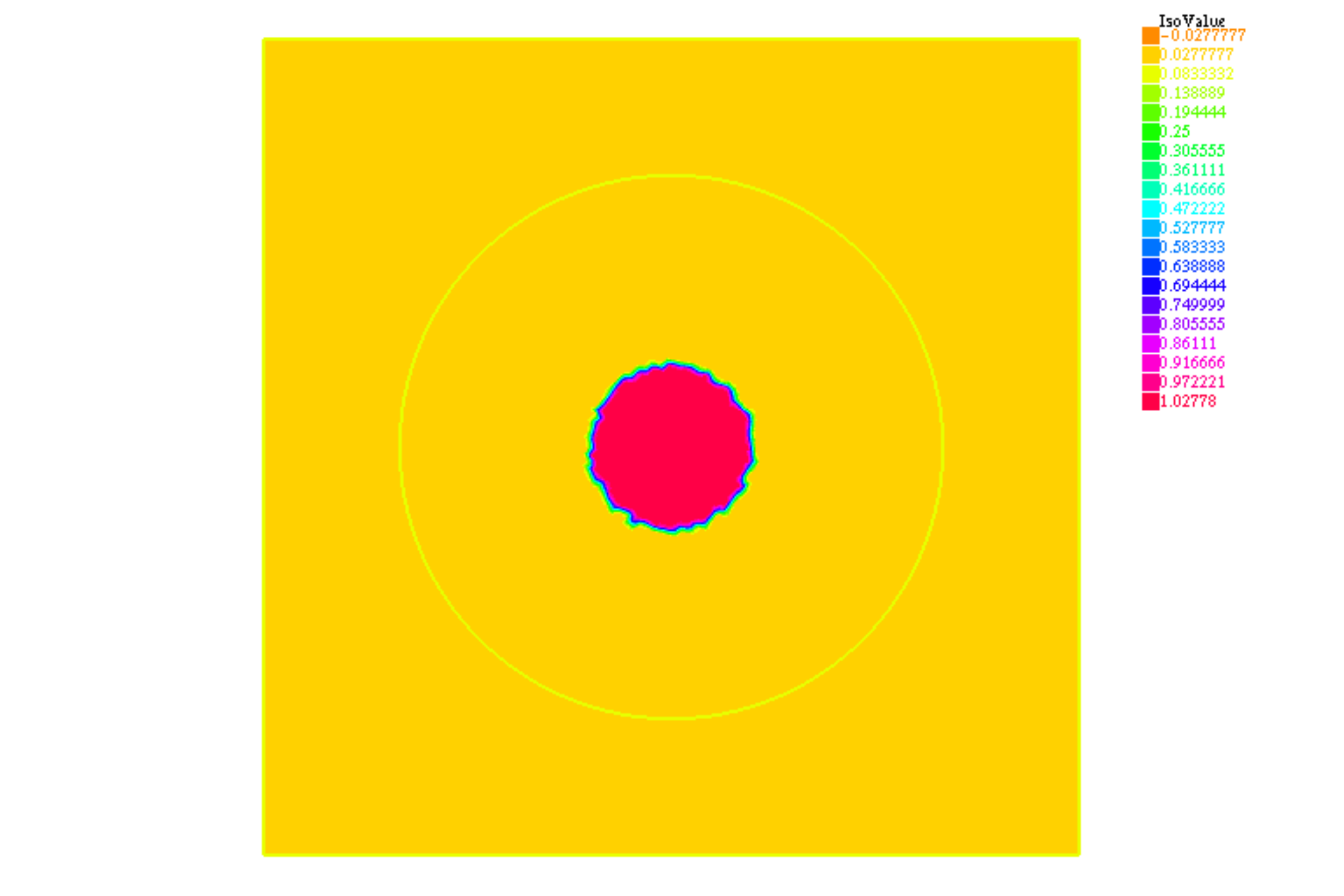}
	\includegraphics[width=0.45\textwidth]{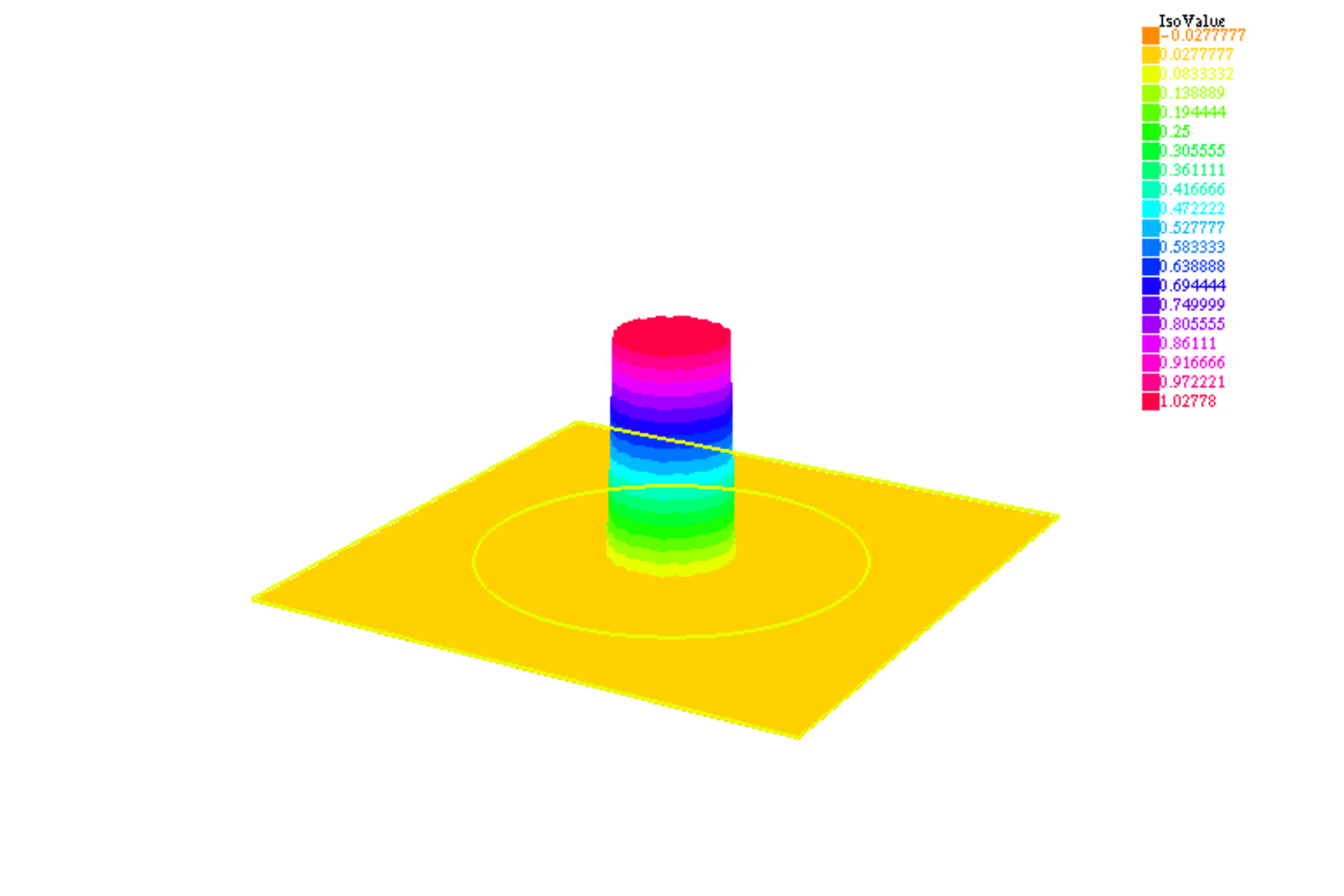} a)
	\includegraphics[width=0.45\textwidth]{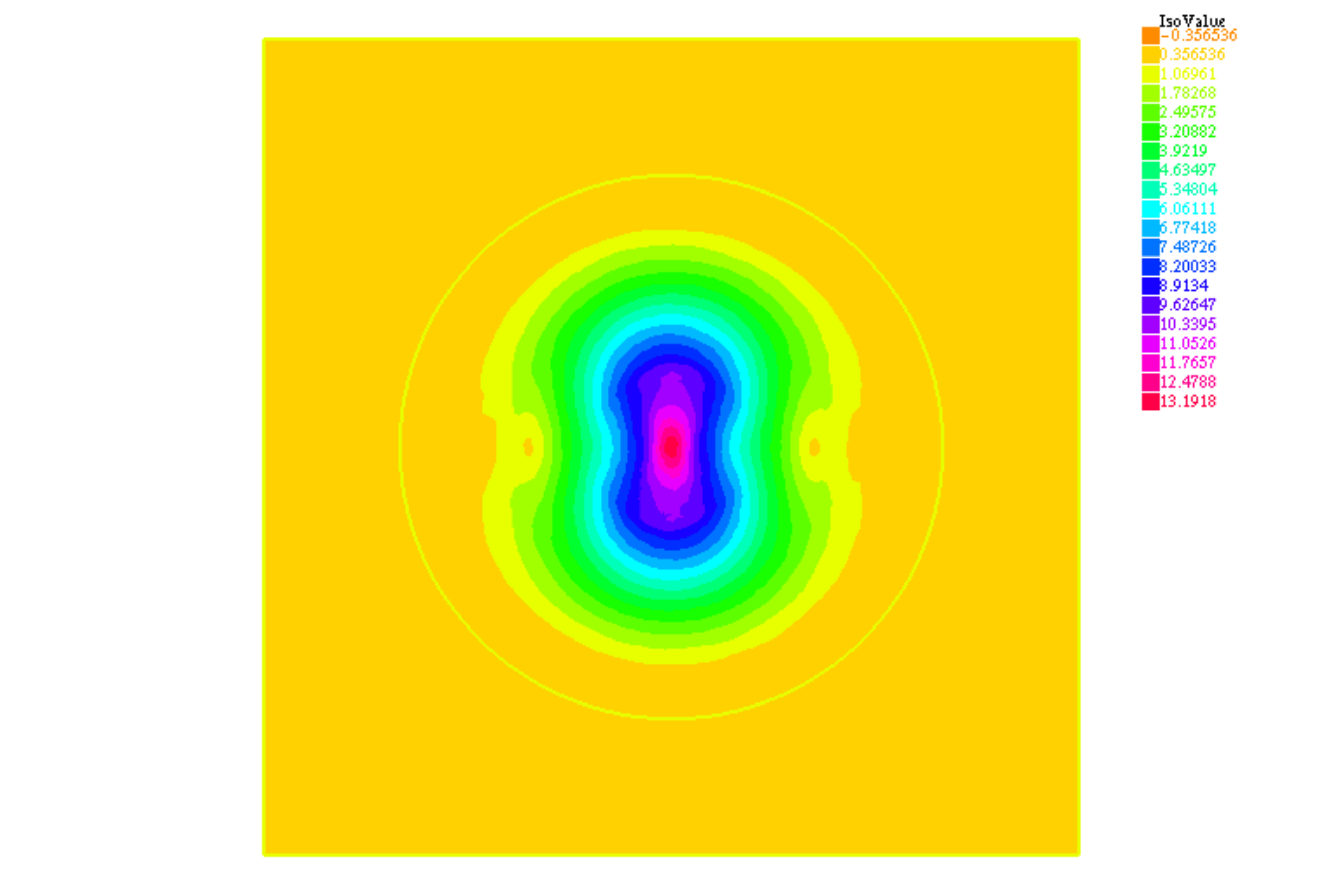}
	   \includegraphics[width=0.45\textwidth]{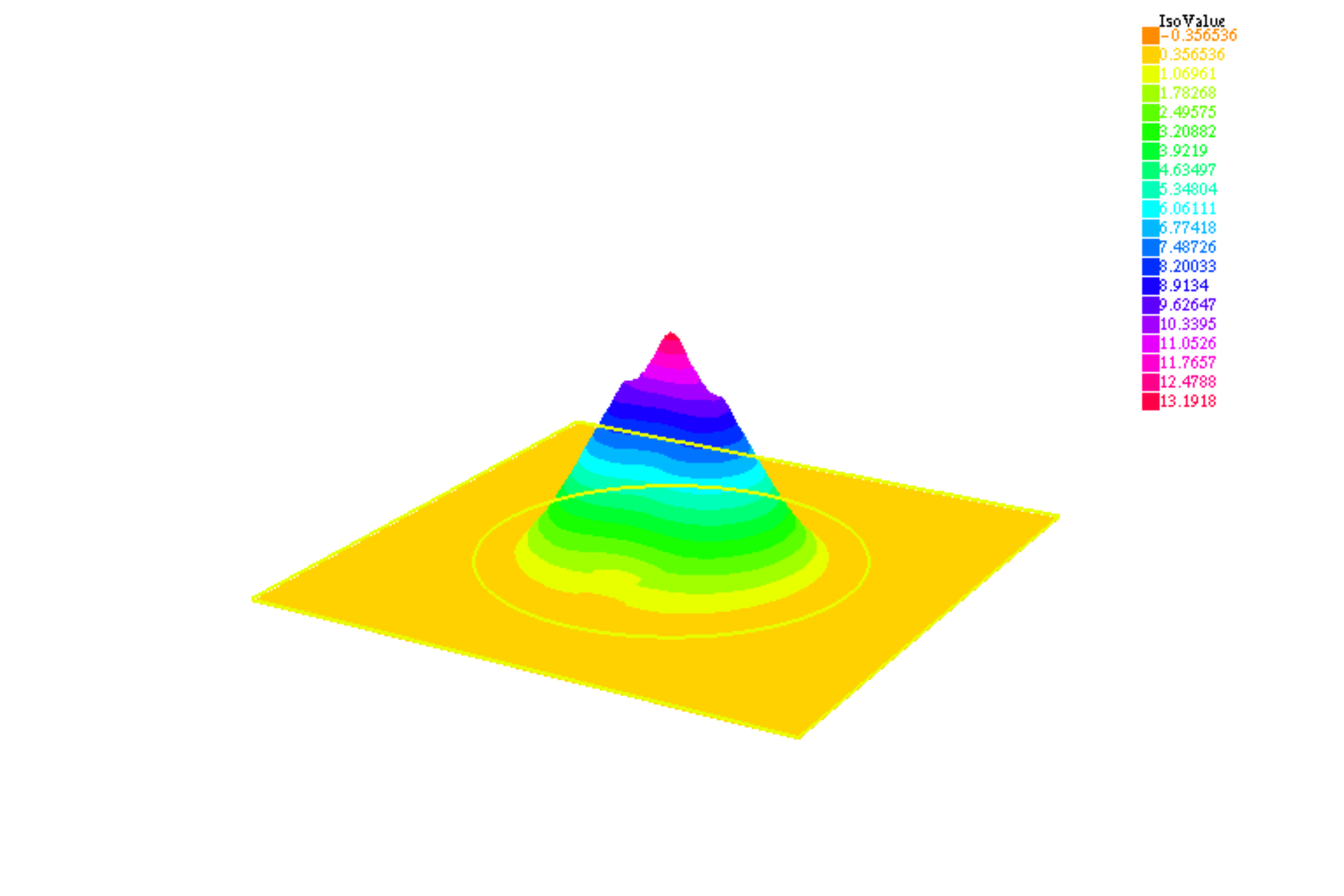} b)

	\includegraphics[width=0.45\textwidth]{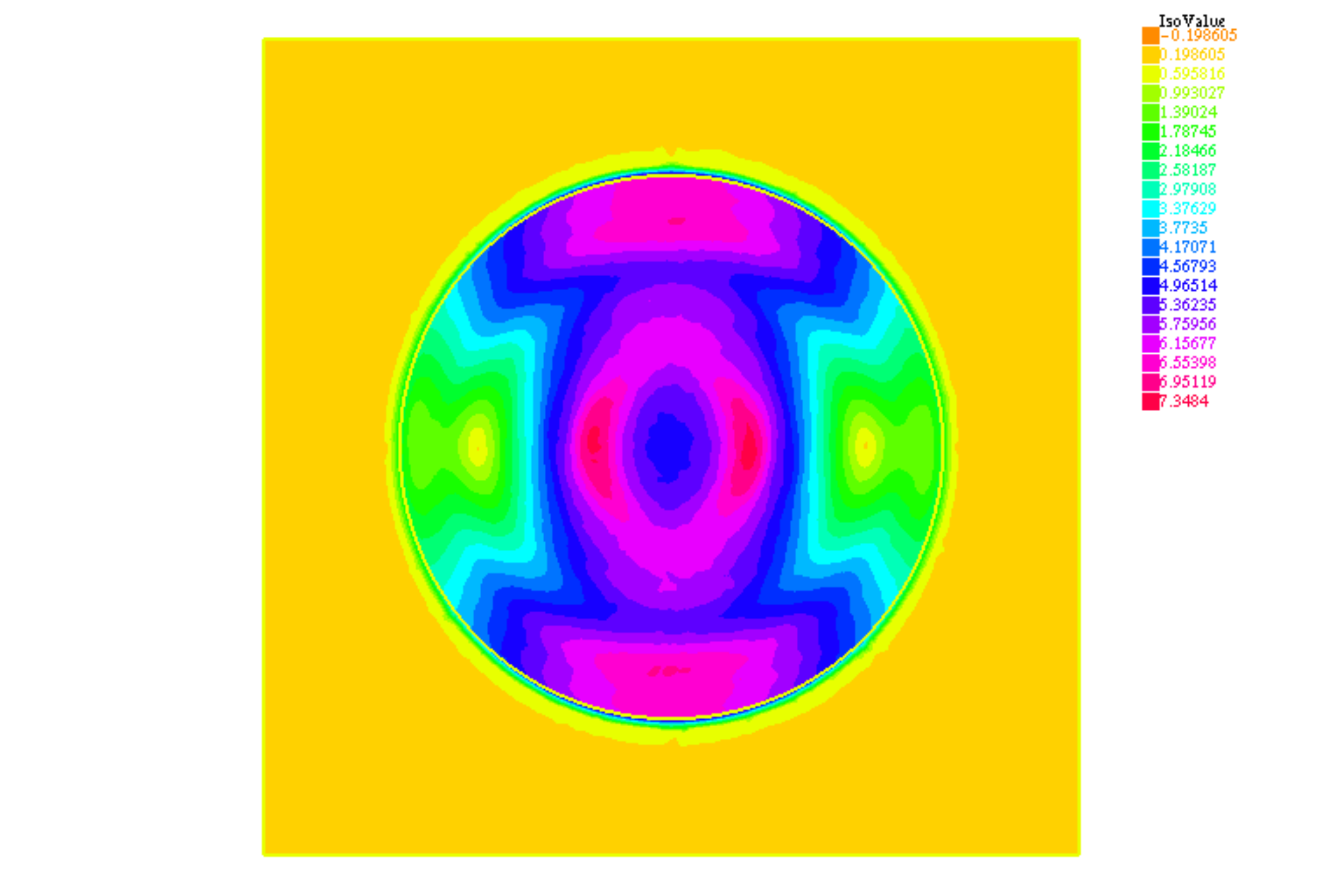}
	\includegraphics[width=0.45\textwidth]{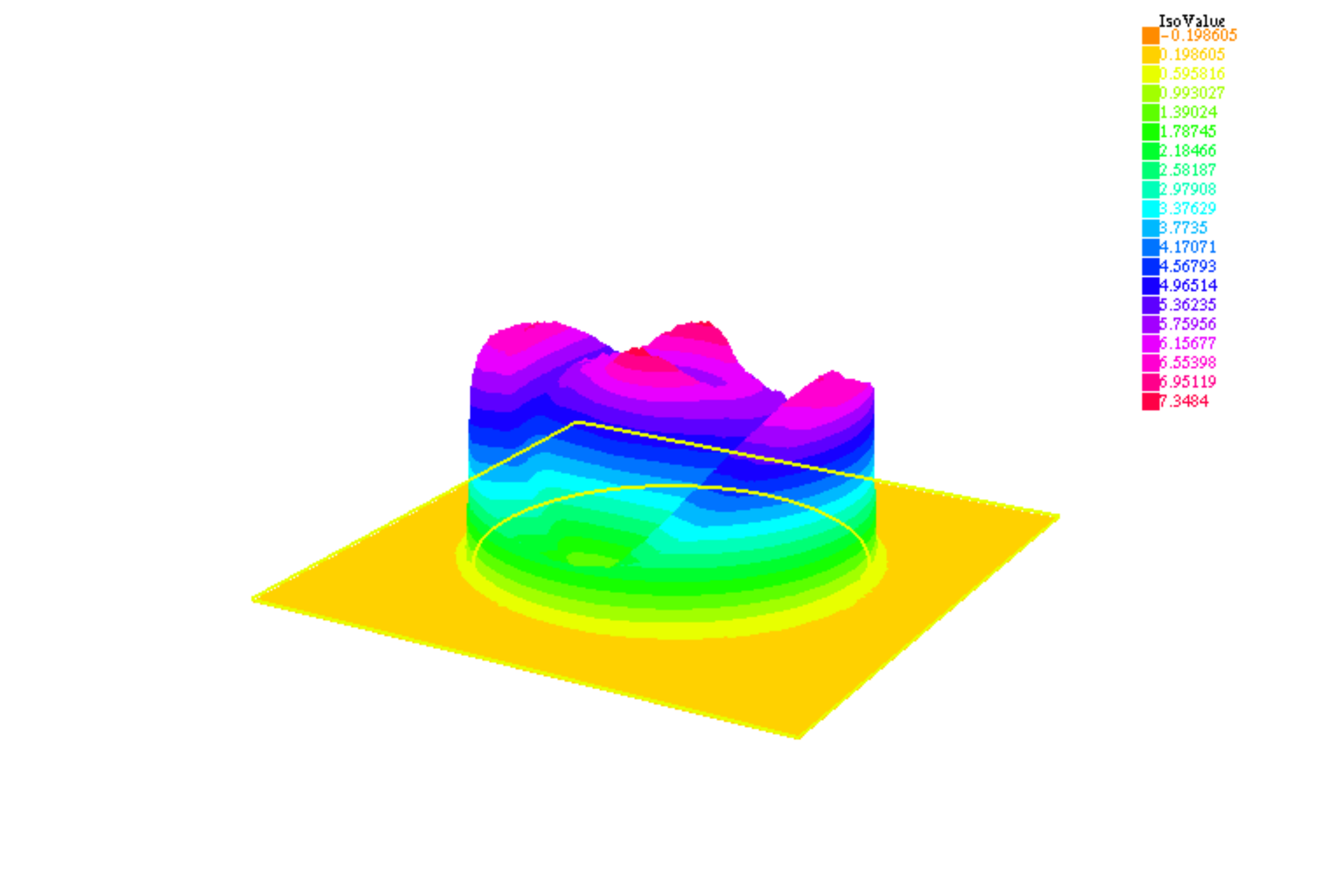} c)
	
	\includegraphics[width=0.45\textwidth]{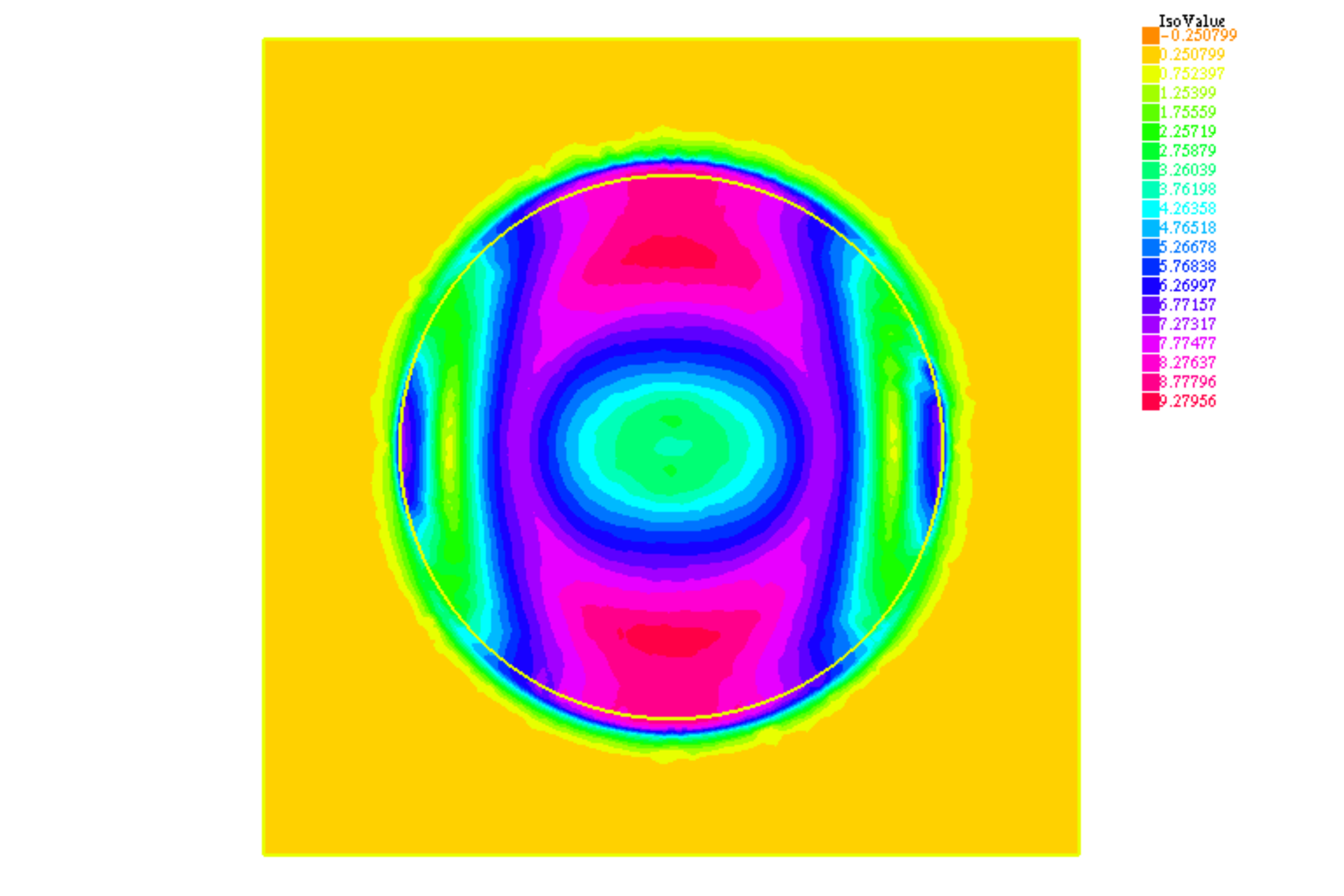}
	\includegraphics[width=0.45\textwidth]{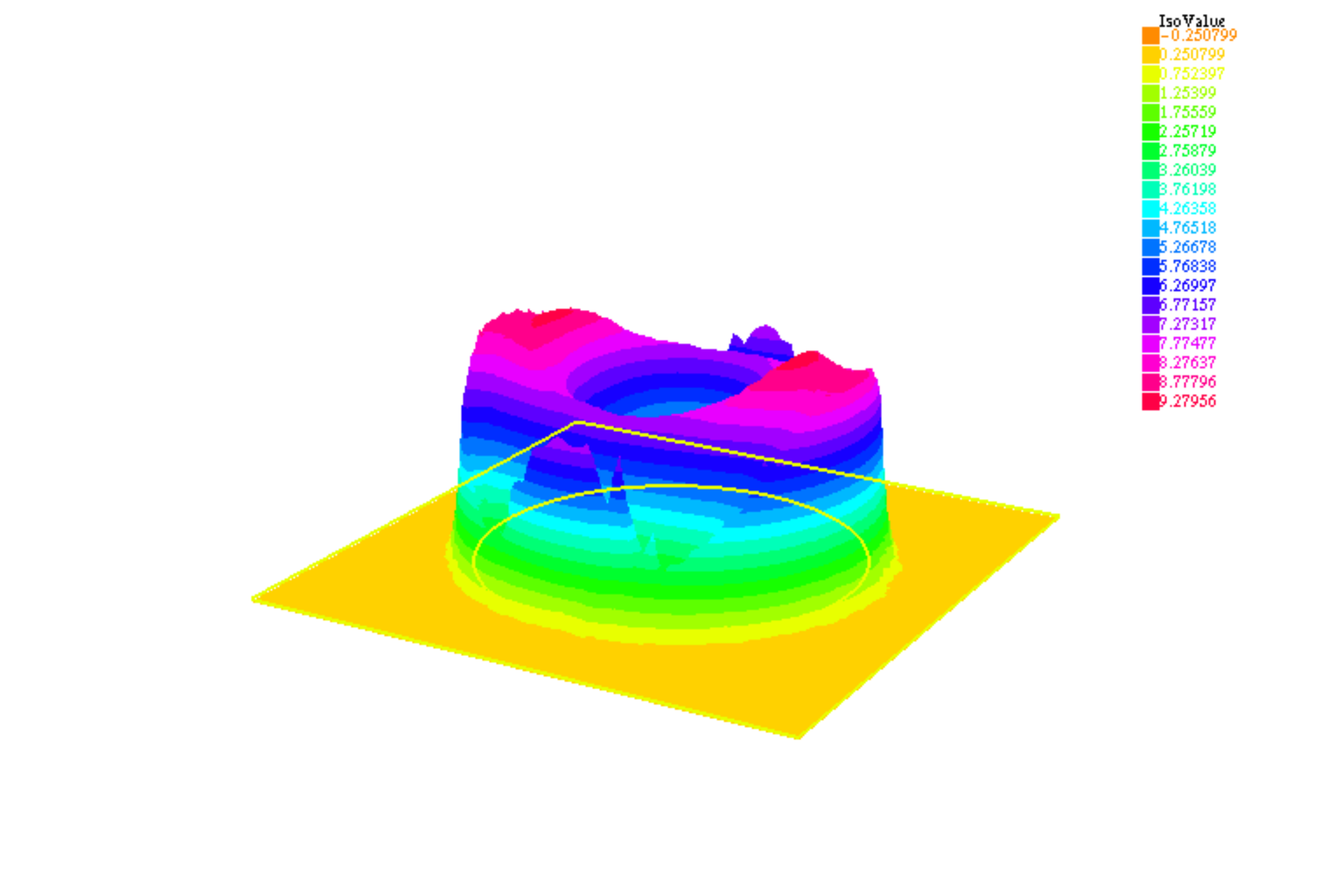} d)

\end{center}
	\caption{Absolute value of the displacement vector field at times a) $ t=0, $ b) $ t=2.5\times 10^{-3}, $ c) $t=5\times 10^{-3}, $ d) $t=7.5\times 10^{-3}$.} 
		\label{modulo}
\end{figure}

Since an analytic solution is not available, we consider as a reference solution that computed by the same space discretization
coupled to a fully implicit Gauss method of order four, see e.g. \cite{hairer:1993}. Therefore, the errors with respect to the reference
solution are an estimate of the time discretization error only.
We compare the implicit Euler method, the Crank-Nicolson method, the Newmark  method with $\beta=1/3$ and $\gamma_N=1/2$ and the
TR-BDF2 method. We include also the implicit Euler method as a representative of a robust, unconditionally stable, albeit first order method.
A pictorial view of the resulting wave propagation on the portion of the domain around Zone 2 at different instants as computed by the TR-BDF2 method is reported 
in Figures \ref{desplazamiento},\ref{modulo}.

 \newpage

 \newpage
  
  \newpage

\begin{figure}
 \begin{center}
	\includegraphics[width=0.45\textwidth]{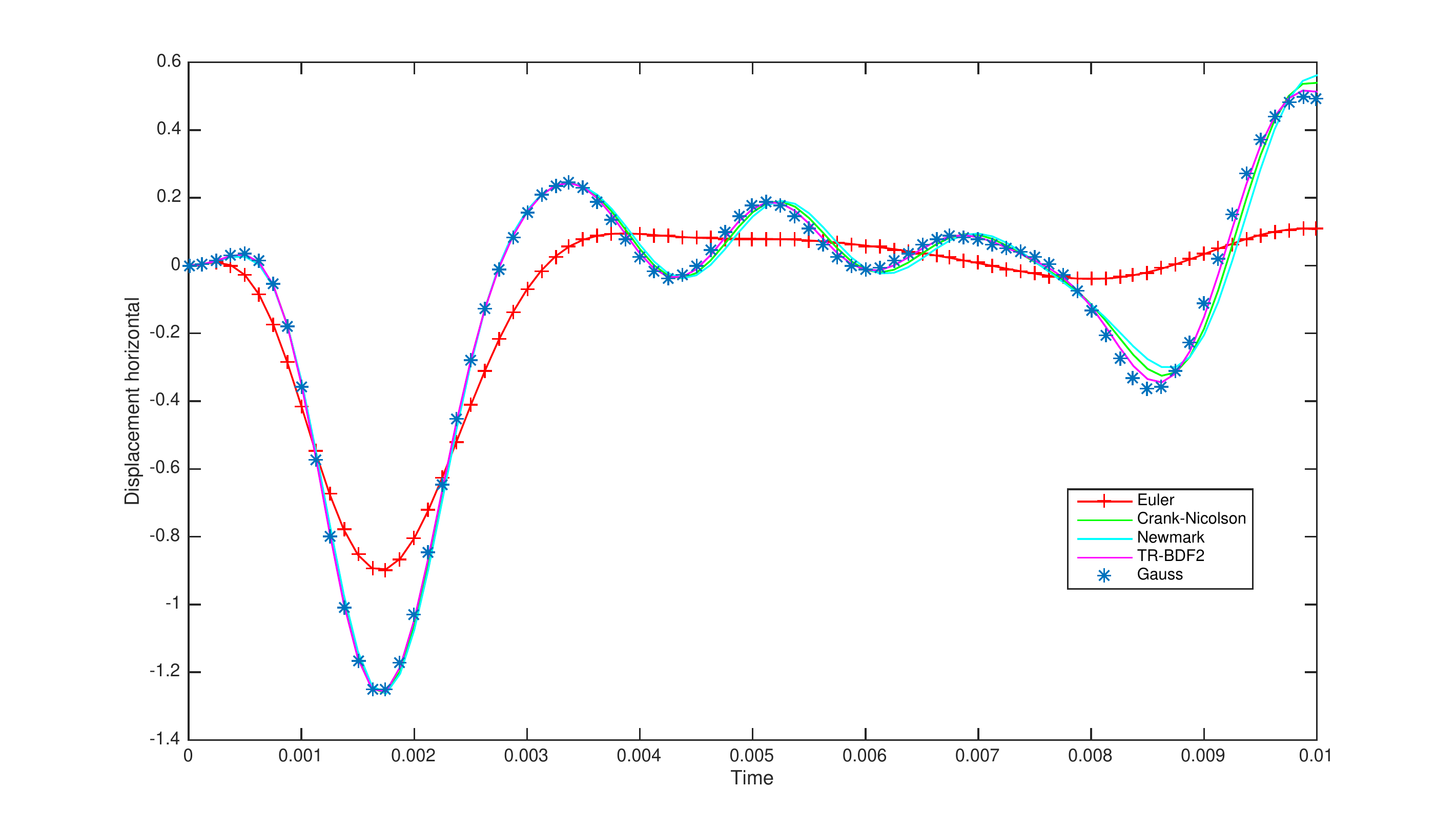}a)
	\includegraphics[width=0.45\textwidth]{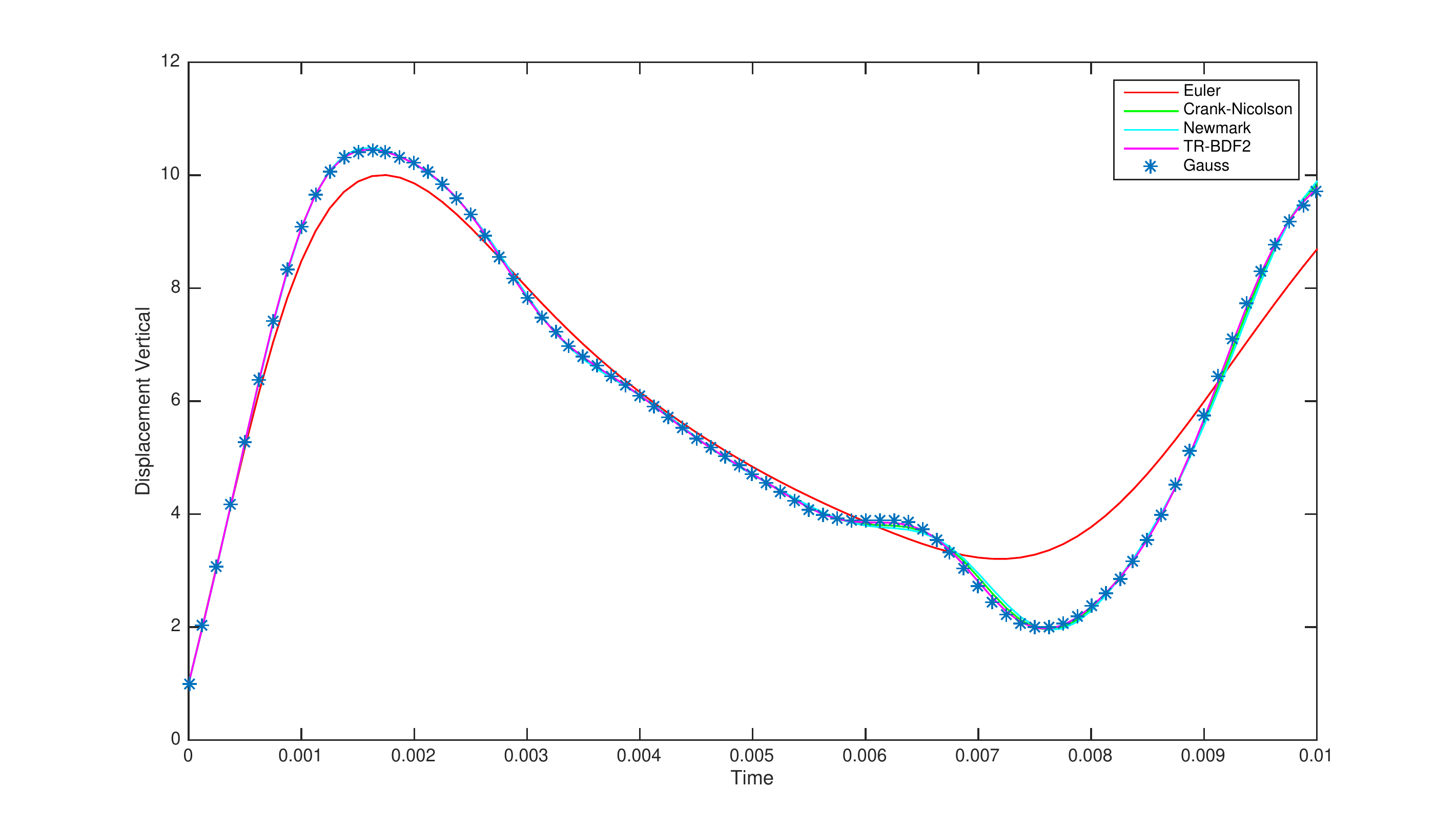}b)
	\includegraphics[width=0.45\textwidth]{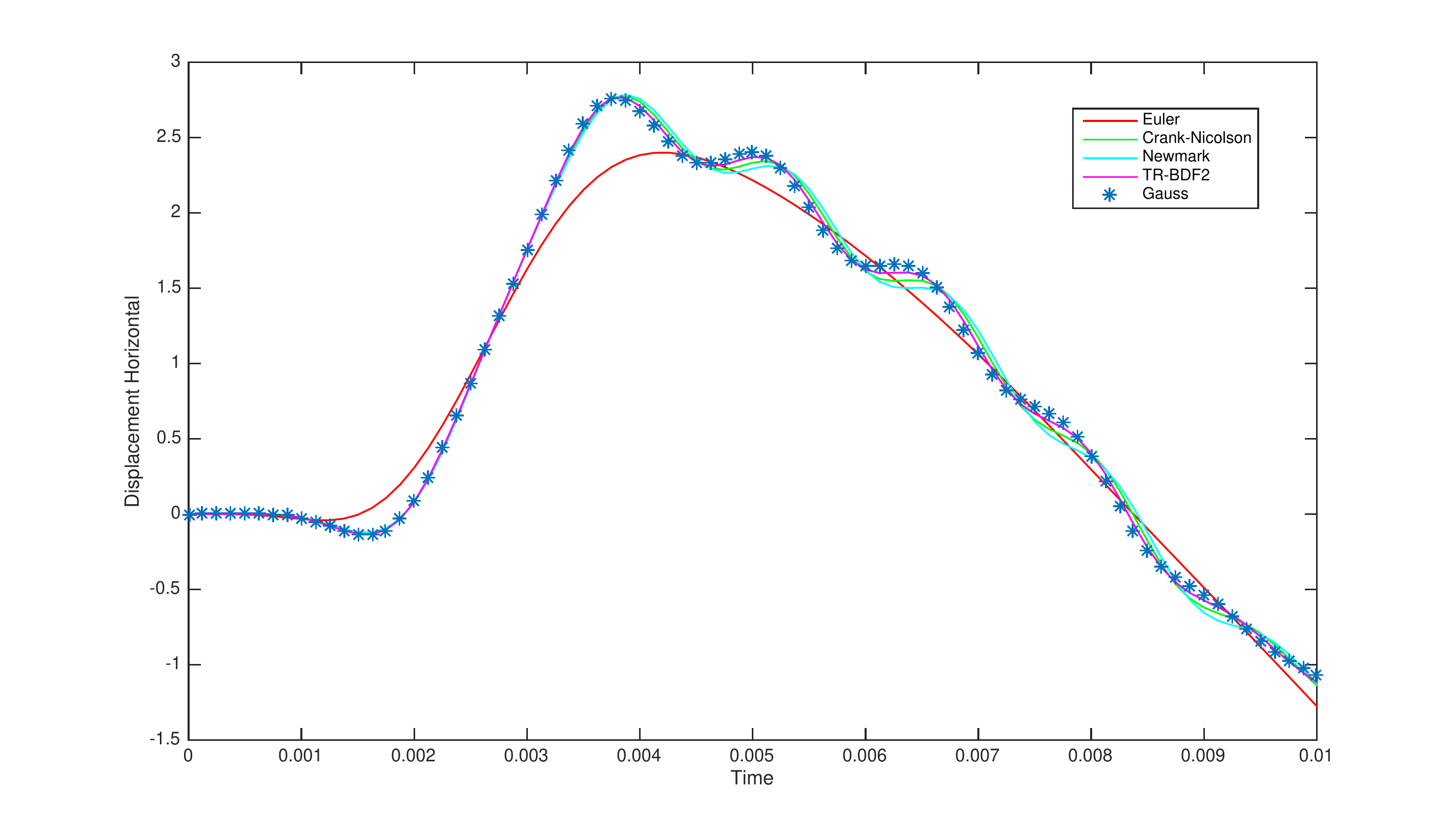}c)
	\includegraphics[width=0.45\textwidth]{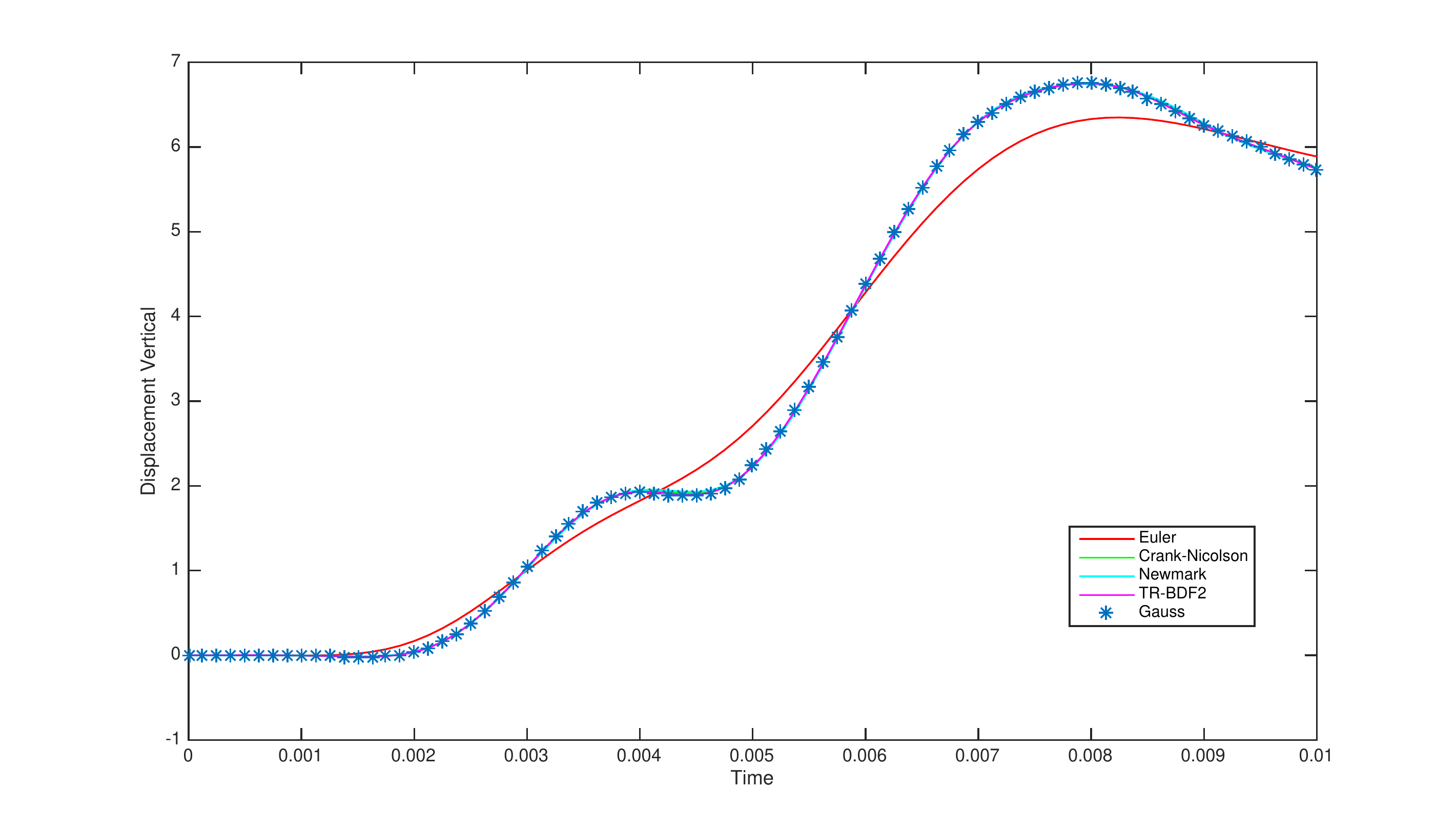}d)
	\includegraphics[width=0.45\textwidth]{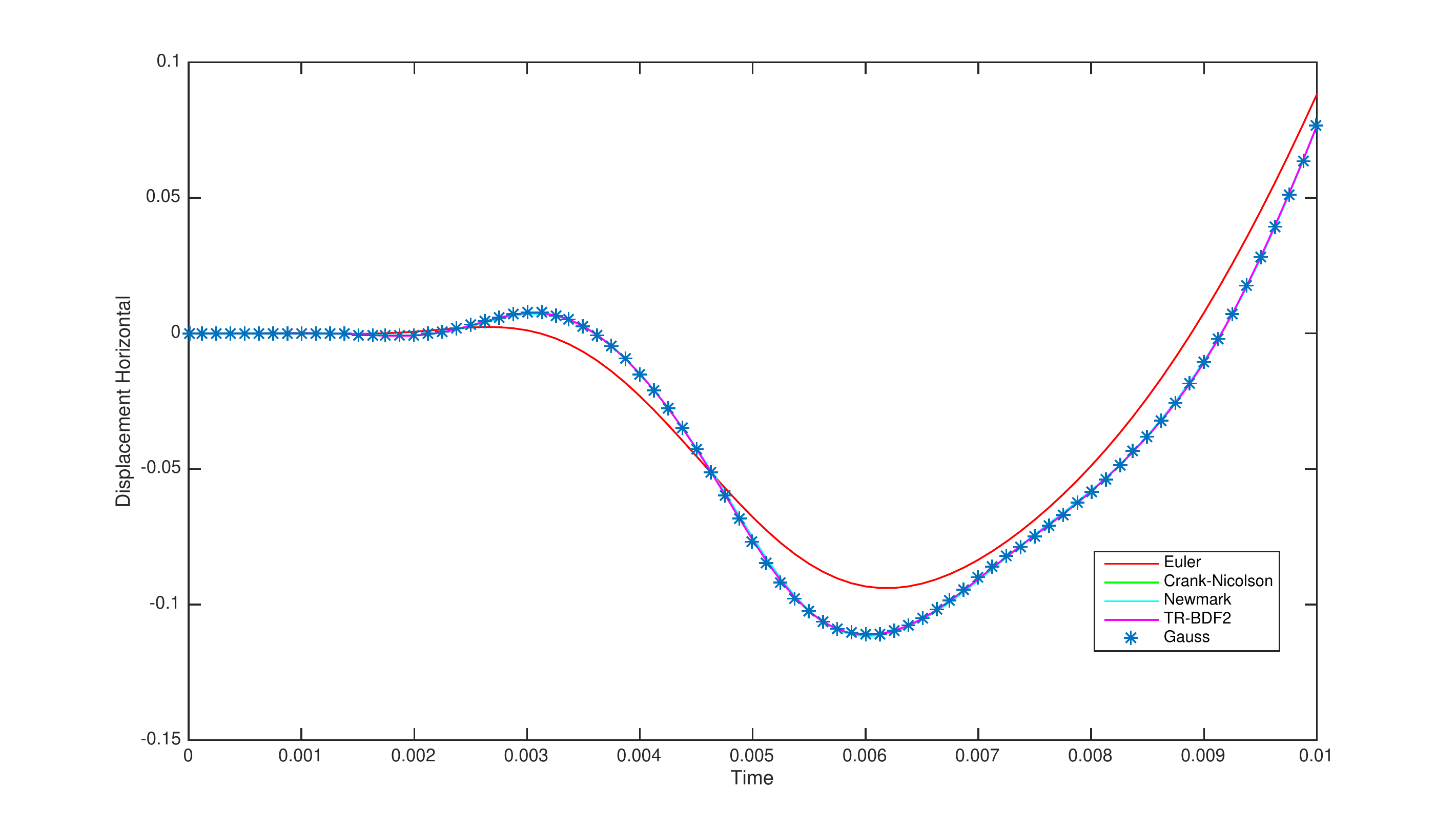}e)
	\includegraphics[width=0.45\textwidth]{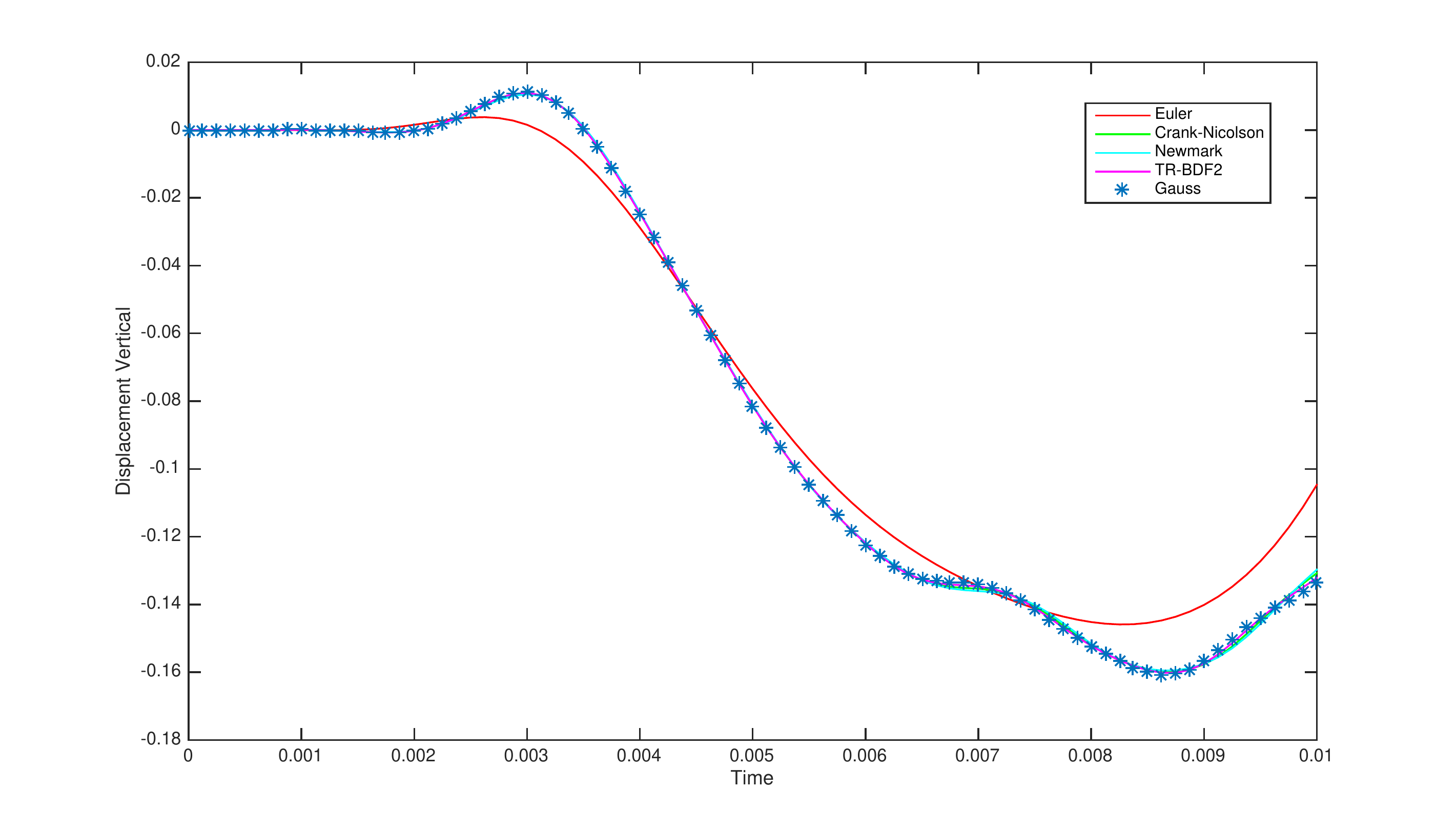}f)
	\end{center}
	\caption{Left column: displacement in the $x-$direction, right column: displacement  in the $y-$direction,
	computed in a)-b) Point A,   c)-d) Point B, e)-f) Point C.}
	\label{Puntossol}
\end{figure}

\begin{figure}
 \begin{center}
	\includegraphics[width=0.45\textwidth]{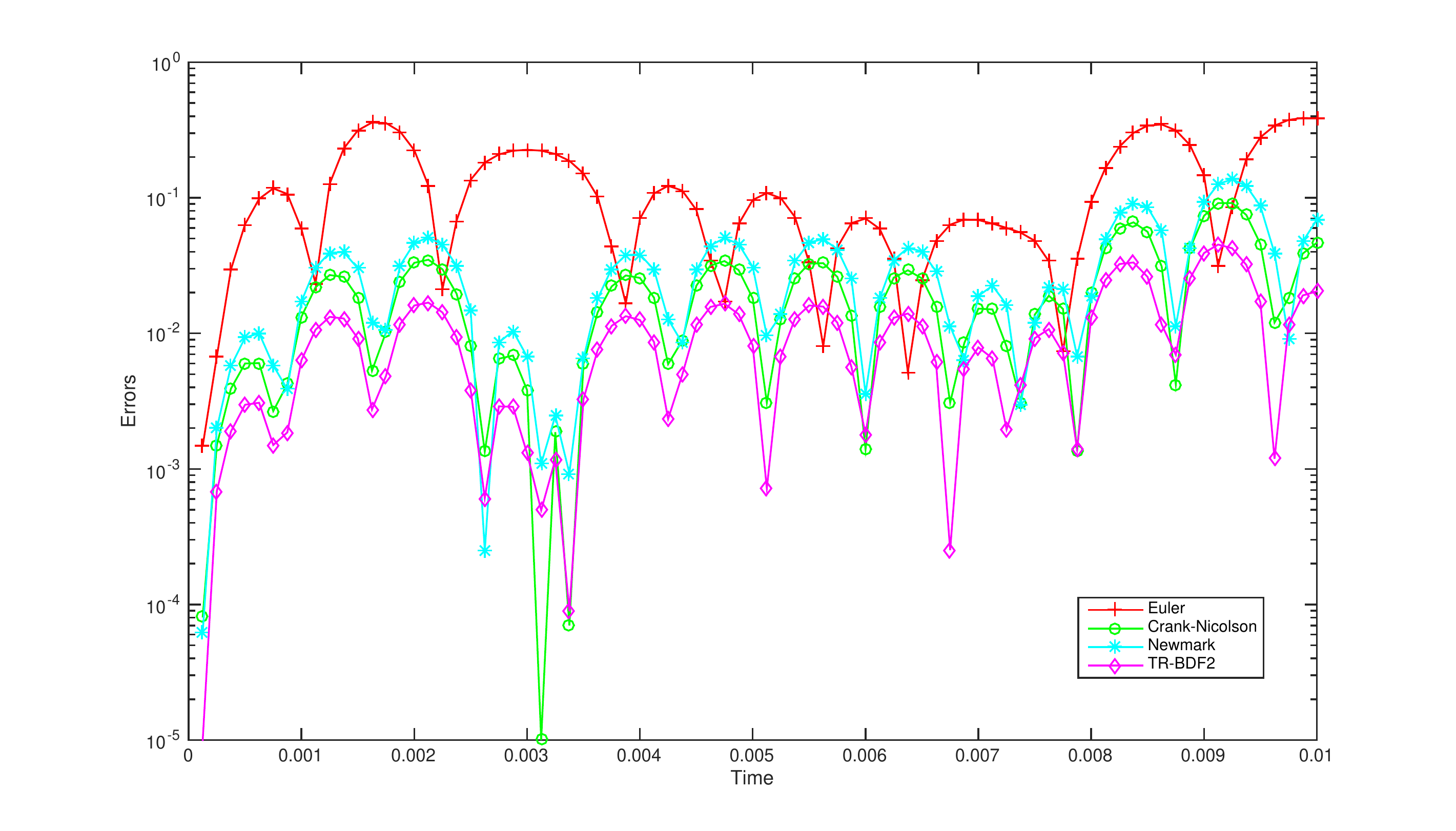}a)
	\includegraphics[width=0.45\textwidth]{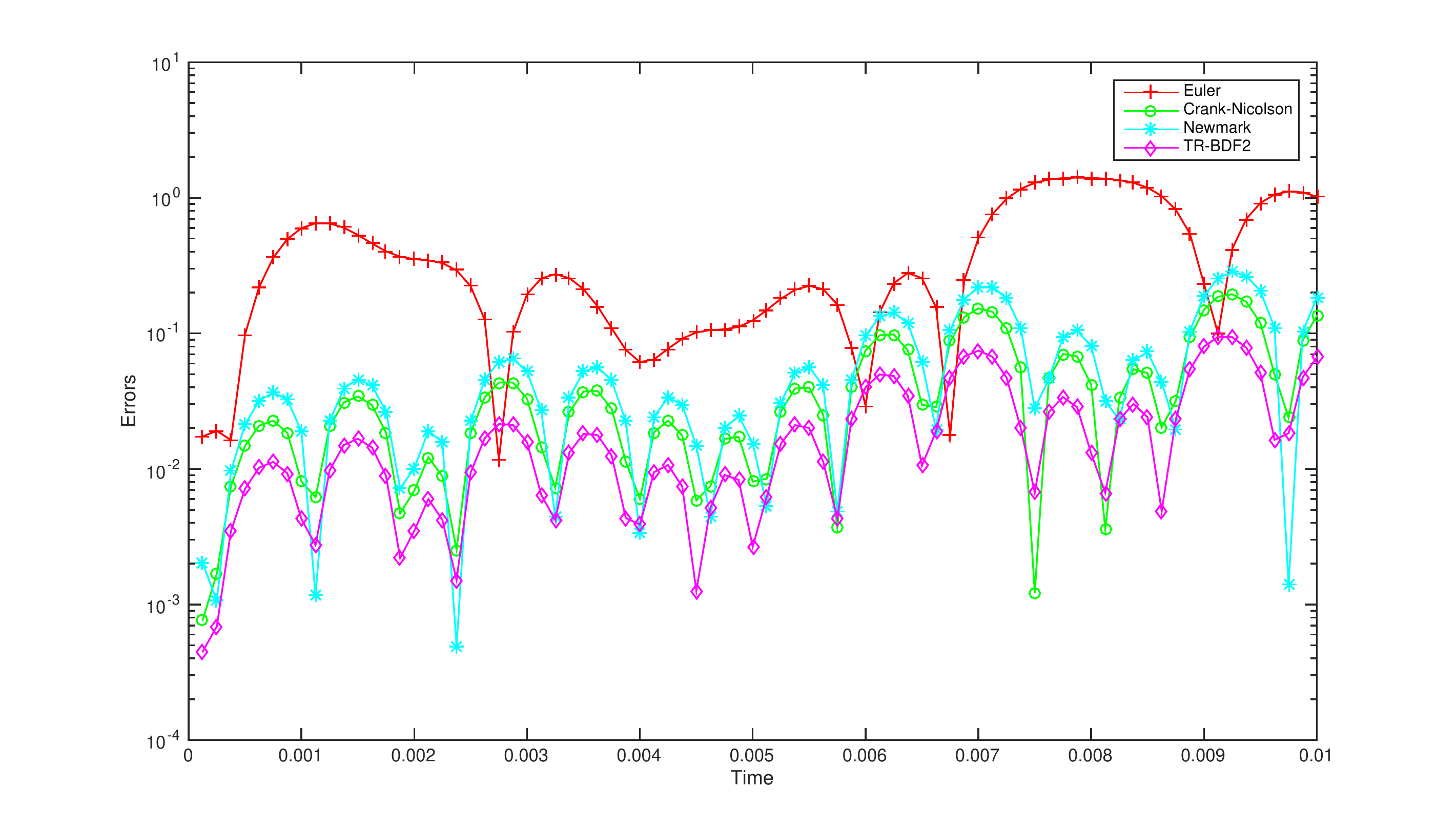}b)
	\includegraphics[width=0.45\textwidth]{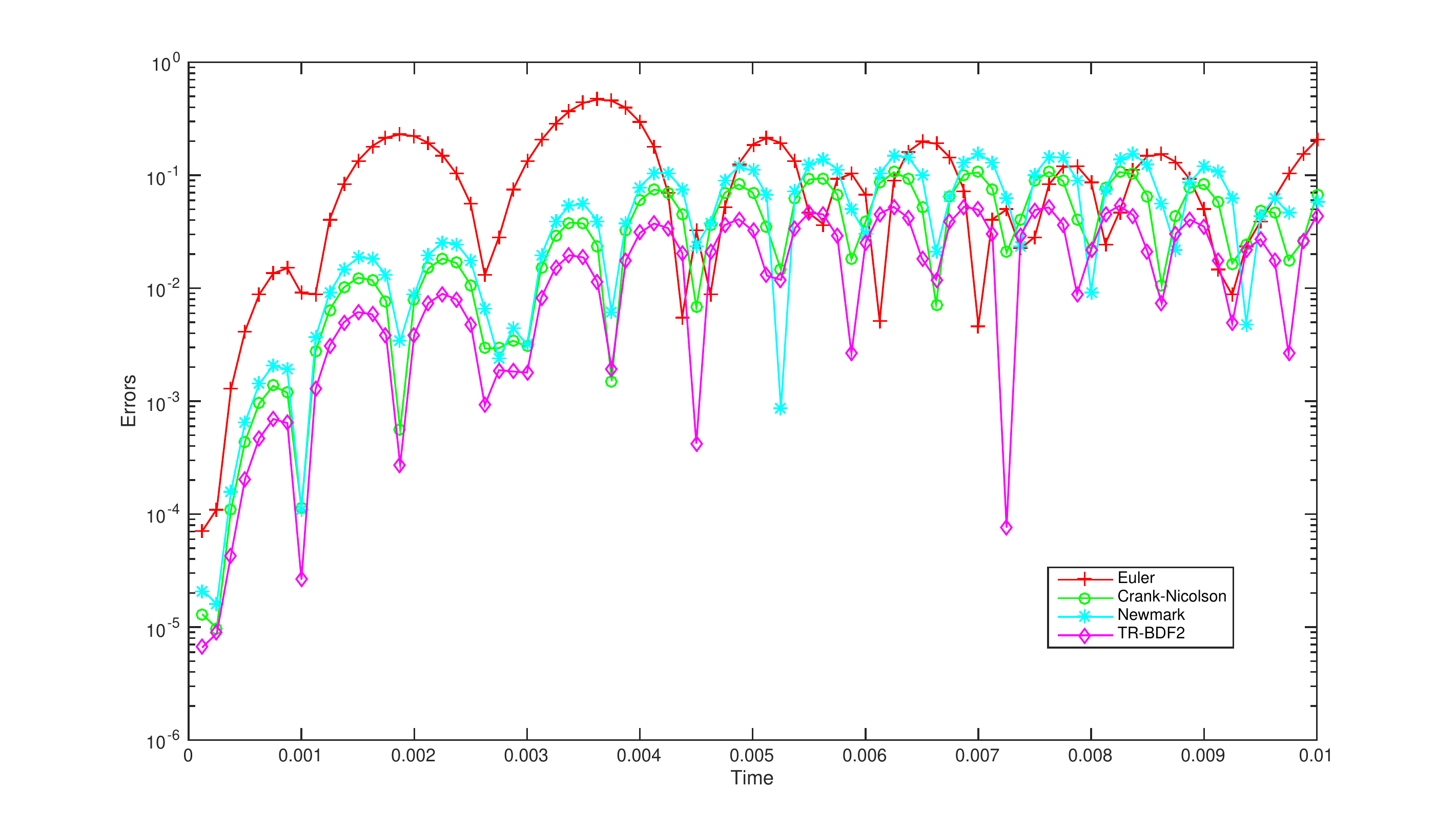}c)
	\includegraphics[width=0.45\textwidth]{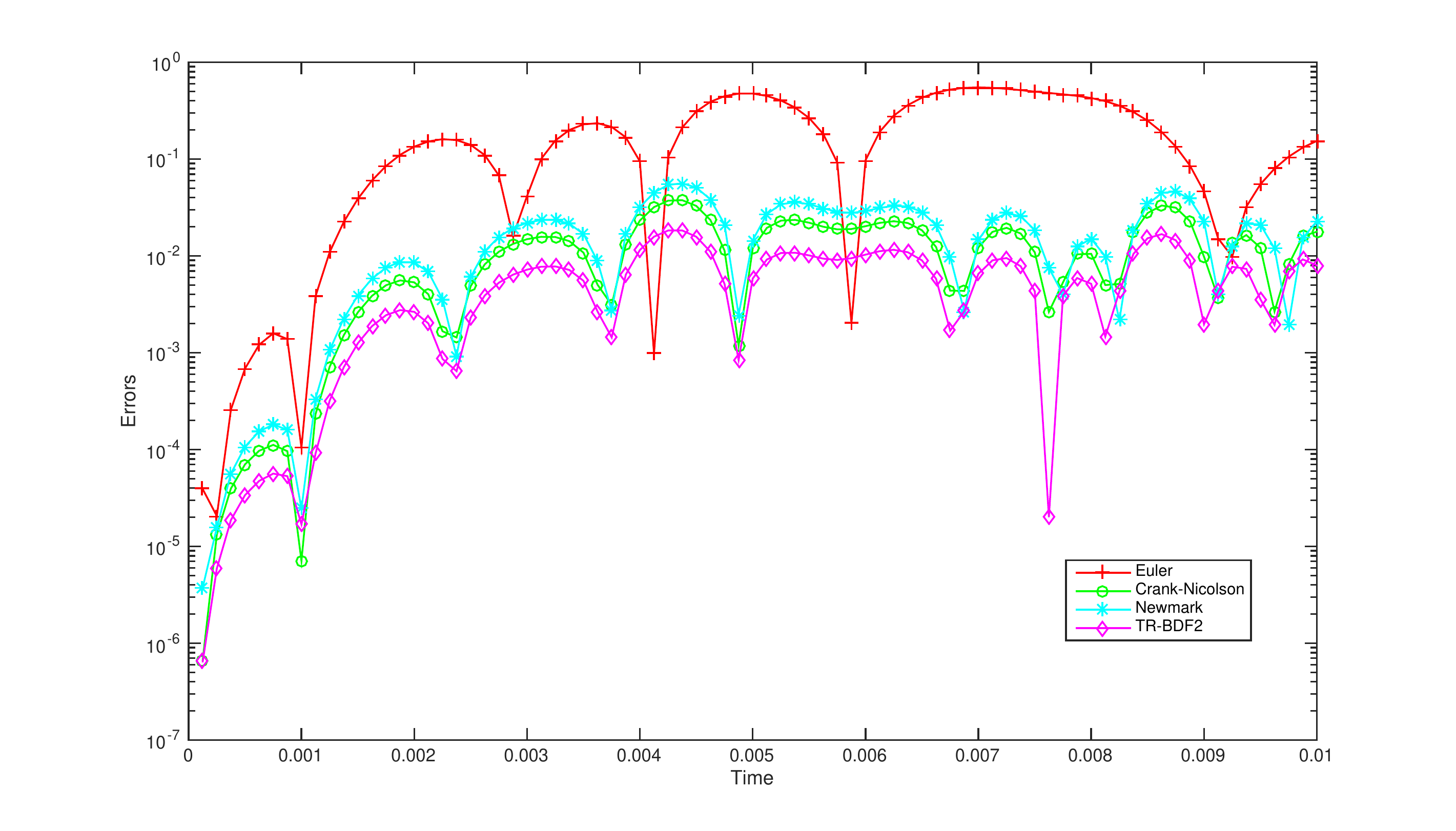}d)
	\includegraphics[width=0.45\textwidth]{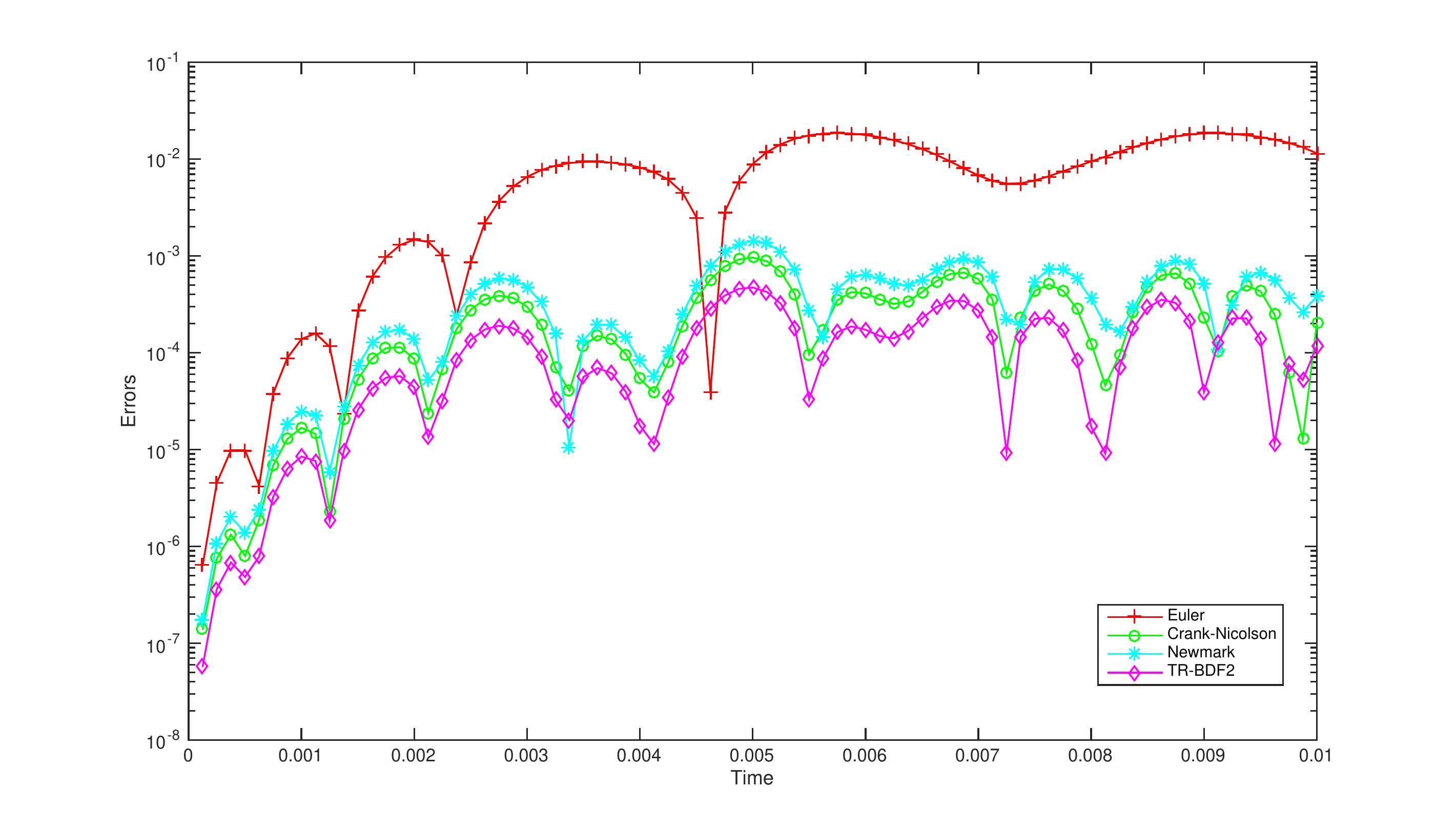}e)
	\includegraphics[width=0.45\textwidth]{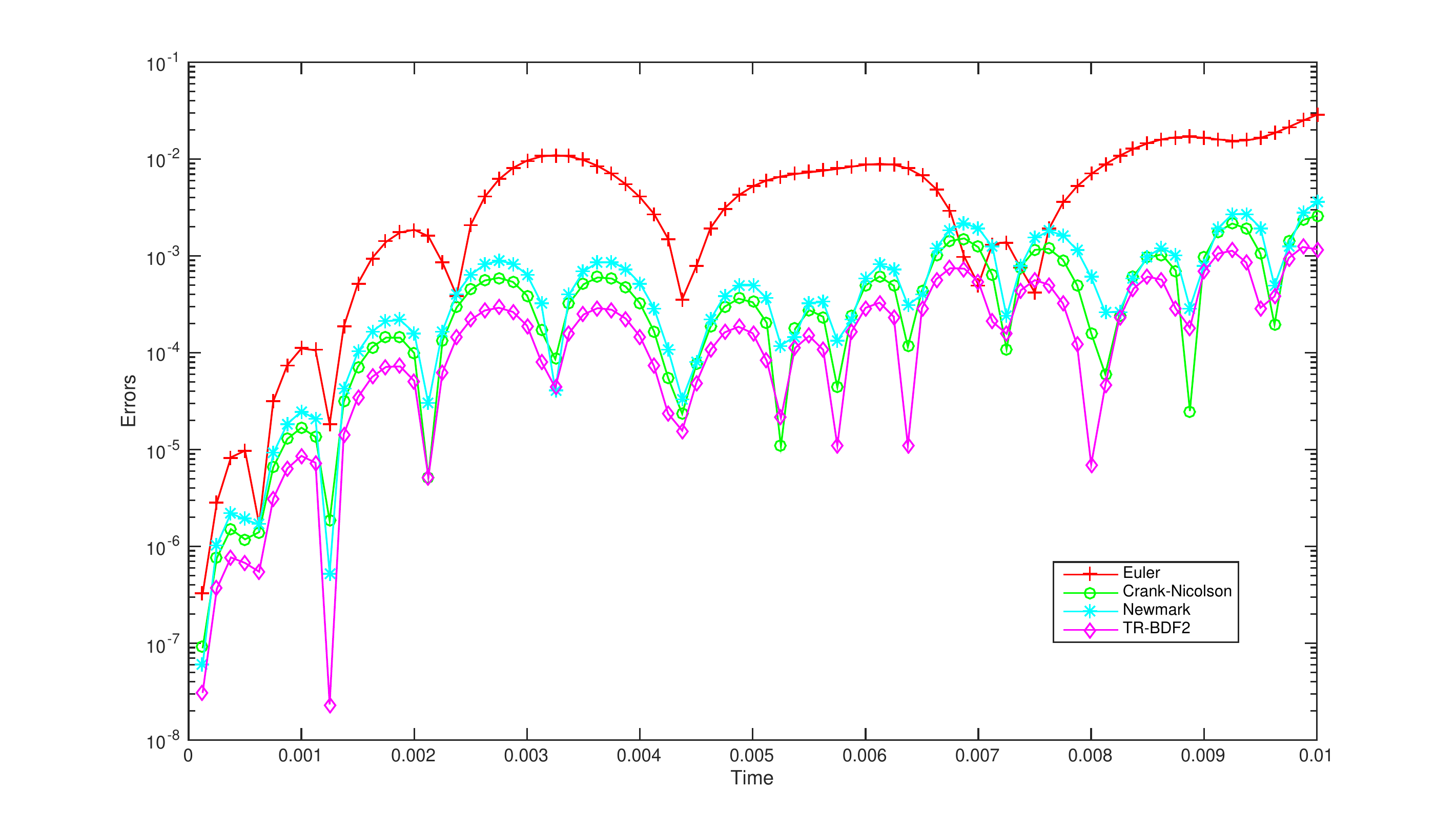}f)
	\end{center}
	\caption{Left column: absolute error on the displacement in the $x-$direction, right column: absolute error on the displacement  in the $y-$direction,
	computed in a)-b) Point A,   c)-d) Point B, e)-f) Point C.}
	\label{Puntosdif}
\end{figure}

The global relative errors with respect to the reference solution are then reported in
Tables \ref{tabla4:hor} and \ref{tabla5:ver} for the two displacement components, respectively.
It can be observed that the 
 TR-BDF2 consistently yields the smallest error values.
 \begin{table}[htbp]
\begin{center}
\begin{tabular}{| | l | l | l | }
\hline
Method & Error $L^\infty(L^2)$  &Error $L^2(H^1)$  \\ \hline \hline
Implicit Euler & 0.69 & 2.64e-2 \\ \hline
Newmark & 0.19 & 1.85e-3\\ \hline
Crank-Nicolson & 8.21e-2 & 1.15e-3  \\ \hline
TR-BDF2 & 4.40e-2 & 5.73e-4 \\ \hline  
\end{tabular}
\caption{Relative errors on the   displacement in the $x-$direction}
\label{tabla4:hor}
\end{center}
\end{table}
 \begin{table}[htbp]
\begin{center}
\begin{tabular}{| | l | l | l | }
\hline
Method & Error $L^\infty(L^2)$  &Error $L^2(H^1)$  \\ \hline \hline
Euler & 7.07e-2 & 1.71e-2 \\ \hline
Newmark & 8.35e-3& 5.29e-3\\ \hline
Crank-Nicolson & 5.89e-3 & 3.68e-4  \\ \hline
TR-BDF2 & 2.93e-3 & 1.82e-4\\ \hline  
\end{tabular}
\caption{Relative errors on the displacement in the $y-$direction}
\label{tabla5:ver}
\end{center}
\end{table}
The behaviour of the computed solutions over the whole time interval $[0,T] $ is also displayed in Figure \ref{Puntossol} 
at the three control points A,B and C  shown in Figure \ref{malla} b). Point A point is located at the center of Zone 1, where the initial displacement occurs,  point B is in  Zone 2 and C is the outer region but close to Zone 2. 
It can be observed that, apart from the solutions computed by the implicit Euler method, the others are all in good agreement.
 The corresponding absolute errors with respect to the reference solution are shown in Figure \ref{Puntosdif},
 while   the error norms at the three control points are reported    
in Tables \ref{tablaA}, \ref{tablaB} and \ref{tablaC}, respectively.  These results  highlight the fact that the TR-BDF2 method yields
more accurate solutions both in the higher and lower Courant number regions, with errors that are consistently at least 50\% smaller than
those of the other methods. This is especially important in view to the combination of these  methods with  higher order discretizations in space.
 
 \newpage

 \begin{table}[htbp]
\begin{center}
\begin{tabular}{ | l | c | c| c|c|}
\hline
Method &  $L^2$ error on $\delta^x$ & $L^\infty$  error on $\delta^x$  &  $L^2$ error on $\delta^y$   & $L^\infty$  error on 
$\delta^y$  \\ \hline \hline
Implicit Euler & 1.74e-2 & 0.39  & 6.31e-2 & 1.41\\ \hline
Newmark & 4.32e-3& 0.137& 9.39e-3&0.28 \\ \hline
Crank-Nicolson & 2.99e-3 & 9.16e-2 & 6.50e-3 & 0.19 \\ \hline
TR-BDF2 & 1.46e-3 & 4.52e-2 & 3.21e-3 & 9.43e-2\\ \hline \hline
\end{tabular}
\caption{Absolute errors on displacement in Point A}
\label{tablaA}
\end{center}
\end{table}

 \begin{table}[htbp]
\begin{center}
\begin{tabular}{ | l | c | c| c|c|}
\hline
Method &   $L^2$ error on $\delta^x$  &  $L^\infty$ error on $\delta^x$  &   $L^2$ error on   $\delta^y$ &  $L^\infty$ error on $\delta^y$  \\ \hline \hline
Implicit Euler & 1.58e-2 & 0.47 & 2.74e-2 & 0.55\\ \hline
Newmark & 7.50e-3& 0.15& 2.31e-3& 5.62e-2\\ \hline
Crank-Nicolson & 5.29e-3 & 0.11 & 1.58e-3 & 3.76e-2\\ \hline
TR-BDF2 & 2.64e-3 & 5.36e-2 & 7.79e-4 & 1.83e-2\\ \hline \hline
\end{tabular}
\caption{Absolute errors on displacement in Point B}
\label{tablaB}
\end{center}
\end{table}
      
 \begin{table}[htbp]
\begin{center}
\begin{tabular}{ | l | c | c| c|c|}
\hline
Method & $L^2$ error on $\delta^x$   & $L^\infty$ error on $\delta^x$  &   $L^2$ error on   $\delta^y$ & $L^\infty$ error on   $\delta^y$\\ \hline \hline
Implicit Euler & 1.04e-3 & 1.86e-2 & 9.14e-4 & 2.86e-2\\ \hline
Newmark & 5.36e-5& 1.43e-3& 1.03e-4&3.63e-3\\ \hline
Crank-Nicolson & 3.65e-5 & 9.63e-4 & 7.68e-5 & 2.60e-3\\ \hline
TR-BDF2 & 1.79e-5 & 4.72e-4 & 3.96e-5 & 1.24e-3\\ \hline \hline
\end{tabular}
\caption{Absolute errors on displacement in Point C}
\label{tablaC}
\end{center}
\end{table}

\section{Conclusions and future work}
\label{conclu} \indent

We have proposed a reformulation of  the  TR-BDF2 method that allows to apply it without overheads in terms of computational cost
or memory to large scale problems in structural mechanics. Our work 
extends the analysis and the comparison of \cite{owren:1995} and a similar reformulation proposed 
for two-stage Rosenbrock methods in \cite{piche:1995}.  In particular, we have presented
  a  reformulation of this method   that only implies the solution of nonlinear systems
  of the same size as the number of degrees of freedom necessary to describe the displacement variables. 
  Furthermore, the velocity  degrees of freedom
  do not have to be stored explicitly and can be recomputed whenever needed, thus avoiding the
  shortcomings of  naive implementations of solvers for first order ODE systems.
  An analysis of the dissipative behaviour of the  method, which
  was carried out considering also the damping term,  shows that  TR-BDF2 is superior  in terms of accuracy and efficiency
 to the classical  G($\alpha$) methods. This finding is confirmed by a number of numerical experiments on significant benchmarks 
 of increasing complexity, in which the TR-BDF2 method   consistently yields errors least 50\% smaller than
those of the other methods, both in high and low Courant number regions. This is especially important in view to the combination of 
robust, unconditionally stable time discretization  methods with  higher order discretizations in space.

 The future developments of this work include the development of a multirate version
 of the proposed method, based on the  self-adjusting multirate technique described in \cite{bonaventura:2020a},
 for application to structural mechanics problems with multiple time scales \cite{gravouil:2001} 
 and to monolithic treatment of fluid structure interaction problems \cite{heil:2008}, as well as the coupling of the proposed time discretization
 to high order, adaptive discontinuous finite element spatial discretizations, so as to achieve the performance improvements
 demonstrated in \cite{tumolo:2015} also in the wave propagation problems typical of structural and seismic engineering.

\section*{Acknowledgements}
L.B. would like to thank T. Chac\'on Rebollo, E. Fern\'andez Nieto and G. Narbona Reina for supporting several visits to Sevilla University,
which have allowed the inception and completion of this work.
This work has been supported by the Spanish Government  Project RTI2018-093521-B-C31. 
 
\bibliographystyle{plain}
\bibliography{struct}

\begin{thebibliography}{10}

\bibitem{bank:1985}
R.E. Bank, W.M. Coughran, W.~Fichtner, E.H. Grosse, D.J. Rose, and R.K. Smith.
\newblock {T}ransient {S}imulation of {S}ilicon {D}evices and {C}ircuits.
\newblock {\em IEEE {T}ransactions on {E}lectron {D}evices}, 32:1992--2007,
  1985.

\bibitem{bathe:1973}
K.J. Bathe and E.~L.Wilson.
\newblock Stability and accuracy analysis of direct integration methods.
\newblock {\em Earthquake Engineering and Structural Dynamics}, 1:283--291,
  1973.

\bibitem{bonaventura:2020a}
L.~Bonaventura, F.~Casella, L.~Delpopolo Carciopolo, and A.~Ranade.
\newblock A self adjusting multirate algorithm for robust time discretization
  of partial differential equations.
\newblock {\em Computers and Mathematics with Applications}, 79:2086--2098,
  2020.

\bibitem{bonaventura:2017}
L.~Bonaventura and A.~Della Rocca.
\newblock Unconditionally strong stability preserving extensions of the
  {TR-BDF2} method.
\newblock {\em Journal of Scientific Computing}, 70:859--895, 2017.

\bibitem{bursi:2008}
O.S. Bursi, A.~Gonzalez-Buelga, L.~Vulcan, S.A. Neild, and D.J. Wagg.
\newblock Novel coupling {R}osenbrock-based algorithms for real-time dynamic
  substructure testing.
\newblock {\em Earthquake Engineering \& Structural Dynamics}, 37:339--360,
  2008.

\bibitem{bursi:2011}
O.S. Bursi, C.~Jua, L.~Vulcan, S.A. Neild, and D.J. Wagg.
\newblock Rosenbrock-based algorithms and subcycling strategies for real-time
  nonlinear substructure testing.
\newblock {\em Earthquake Engineering \& Structural Dynamics}, 40:1--19, 2011.

\bibitem{chung:1993}
J.~Chung and G.M. Hulbert.
\newblock A time integration algorithm for structural dynamics with improved
  numerical dissipation: the generalized-$\alpha$ method.
\newblock {\em Journal of Applied Mechanics}, 60:371--375, 1993.

\bibitem{edwards:2011}
J.D. Edwards, J.E. Morel, and D.A. Knoll.
\newblock Nonlinear variants of the {TR/BDF2} method for thermal radiative
  diffusion.
\newblock {\em Journal of Computational Physics}, 230:1198--1214, 2011.

\bibitem{erlicher:2002}
S.~Erlicher, L.~Bonaventura, and O.~S. Bursi.
\newblock The analysis of the generalized-$\alpha$ method for non-linear
  dynamic problems.
\newblock {\em Computational Mechanics}, 28:83--104, 2002.

\bibitem{gravouil:2001}
A.~Gravouil and A.~Combescure.
\newblock Multi-time-step explicit--implicit method for non-linear structural
  dynamics.
\newblock {\em International Journal for Numerical Methods in Engineering},
  50:199--225, 2001.

\bibitem{hairer:1993}
E.~Hairer, H.P. Norsett, and G.~Wanner.
\newblock {\em Solving {O}rdinary {D}ifferential {E}quations. {I} {N}onstiff
  {P}roblems (2nd.\ Revised Edition)}.
\newblock Springer-Verlag, 1993.

\bibitem{hamkar:2012}
A.W. Hamkar, S.~Hartmann, and J.~Rang.
\newblock A stiffly accurate {R}osenbrock-type method of order 2 applied to
  {FE}-analyses in finite strain viscoelasticity.
\newblock {\em Applied Numerical Mathematics}, 62:1837--1848, 2012.

\bibitem{hartmann:2007}
S.~Hartmann and J.~Wensch.
\newblock Finite element analysis of viscoelastic structures using
  {R}osenbrock-type methods.
\newblock {\em Computational Mechanics}, 40:383--398, 2007.

\bibitem{hecht:2012}
F.~Hecht.
\newblock New development in freefem++.
\newblock {\em Journal of Numerical Mathematics}, 20:251--265, 2012.

\bibitem{heil:2008}
M.~Heil, A.L. Hazel, and J.~Boyle.
\newblock Solvers for large-displacement fluid--structure interaction problems:
  segregated versus monolithic approaches.
\newblock {\em Computational Mechanics}, 43:91--101, 2008.

\bibitem{hosea:1996}
M.E. Hosea and L.F. Shampine.
\newblock Analysis and implementation of {TR}-{BDF}2.
\newblock {\em Applied Numerical Mathematics}, 20:21--37, 1996.

\bibitem{hughes:1978}
T.~Hughes and W.Liu.
\newblock Implicit-explicit finite elements in transient analysis:
  implementation and numerical examples.
\newblock {\em Journal of Applied Mechanics, Transactions ASME},
  45(2):375--378, 1978.

\bibitem{lamarche:2009}
C.P. Lamarche, A.~Bonelli, O.S. Bursi, and R.~Tremblay.
\newblock A {R}osenbrock-{W} method for real-time dynamic substructuring and
  pseudo-dynamic testing.
\newblock {\em Earthquake Engineering \& Structural Dynamics}, 38:1071--1092,
  2009.

\bibitem{leveque:2004}
R.J. LeVeque.
\newblock {\em Finite Volume Methods for Hyperbolic Problems}.
\newblock Cambridge University Press, 2004.

\bibitem{meijaard:2003}
J.P. Meijaard.
\newblock Application of {R}unge--{K}utta--{R}osenbrock methods to the analysis
  of flexible multibody systems.
\newblock {\em Multibody System Dynamics}, 10:263--288, 2003.

\bibitem{newmark:1959}
N.~M. Newmark.
\newblock A method of computation for structural dynamics.
\newblock {\em Journal of the Engineering Mechanics Division}, 85:67--94, 1959.

\bibitem{owren:1995}
B.~Owren and H.~H. Simonsen.
\newblock Alternative integration methods for problems in structural dynamics.
\newblock {\em Computer Methods in Applied Mechanics and Engineering},
  122:1--10, 1995.

\bibitem{piche:1995}
R.~Pich\'e.
\newblock An {L}-stable {R}osenbrock method for step-by-step time integration
  in structural dynamics.
\newblock {\em Computer Methods in Applied Mechanics and Engineering},
  126:343--354, 1995.

\bibitem{tumolo:2016}
G.~Tumolo.
\newblock A mass conservative {TR-BDF2} semi-implicit semi-{L}agrangian {DG}
  discretization of the shallow water equations on general structured meshes of
  quadrilaterals.
\newblock {\em Communications in Applied and Industrial Mathematics},
  7:165--190, 2016.

\bibitem{tumolo:2015}
G.~Tumolo and L.~Bonaventura.
\newblock A semi-implicit, semi-{L}agrangian, {DG} framework for adaptive
  numerical weather prediction.
\newblock {\em Quarterly Journal of the Royal Meteorological Society},
  141:2582--2601, 2015.

\end{thebibliography}

\end{document}